\documentclass[10pt]{article}

\usepackage{amsmath}
\usepackage{amssymb}
\usepackage{latexsym}
\usepackage{amsthm}

\renewcommand{\theenumi}{\roman{enumi}}

\newcommand{\iterate}[1]{^{(#1)}}
\newcommand{\UpperSpinor}{\Tilde\chi}
\newcommand{\LowerSpinor}{\Tilde\eta}
\newcommand{\Pauli}{\chi_{P}}
\newcommand{\PauliData}{\chi_{P0}}
\newcommand{\FT}{\, \widehat{} \,\,}
\newcommand{\LPproj}{\Delta}
\newcommand{\eps}{^{\varepsilon}}
\newcommand{\ext}{_{\mathrm{ext}}}
\newcommand{\subeps}{_{\varepsilon}}
\newcommand{\init}{\vert_{t = 0}}
 
\newcommand{\abs}[1]{\left\vert #1 \right\vert}
\newcommand{\bigabs}[1]{\bigl\vert #1 \bigr\vert}

\newcommand{\norm}[1]{\left\Vert #1 \right\Vert}
\newcommand{\bignorm}[1]{\bigl\Vert #1 \bigr\Vert}

\newcommand{\Sobnorm}[2]{\norm{#1}_{H^{#2}}}
\newcommand{\bigSobnorm}[2]{\bignorm{#1}_{H^{#2}}}
\newcommand{\Sobdotnorm}[2]{\norm{#1}_{\dot{H}^{#2}}}
\newcommand{\bigSobdotnorm}[2]{\bignorm{#1}_{\dot{H}^{#2}}}
\newcommand{\Lpnorm}[2]{\norm{#1}_{L^{#2}}}

\newcommand{\Lxpnorm}[2]{\norm{#1}_{L_{x}^{#2}}}
\newcommand{\bigLxpnorm}[2]{\bignorm{#1}_{L_{x}^{#2}}}

\newcommand{\twonorm}[2]{\norm{#1}_{L^2#2}}
\newcommand{\bigtwonorm}[2]{\bignorm{#1}_{L^2#2}}
\newcommand{\inftynorm}[2]{\norm{#1}_{L^\infty#2}}

\newcommand{\mixednorm}[3]{\norm{#1}_{L_{t}^{#2}L_{x}^{#3}}}
\newcommand{\mixednormlocal}[3]{\norm{#1}_{L_{t}^{#2}L_{x}^{#3}(S_{T})}}
\newcommand{\LpHs}[3]{\norm{#1}_{L_{t}^{#2}H^{#3}}}

\newcommand{\LpHslocal}[3]{\norm{#1}_{L_{t}^{#2}H^{#3}(S_{T})}}
\newcommand{\bigLpHslocal}[3]{\bignorm{#1}_{L_{t}^{#2}H^{#3}(S_{T})}}

\newcommand{\mixed}[2]{L_{t}^{#1}L_{x}^{#2}}
\newcommand{\bigmixednorm}[3]{\bignorm{#1}_{L_{t}^{#2}L_{x}^{#3}}}
\newcommand{\bigmixednormlocal}[3]{\bignorm{#1}_{L_{t}^{#2}L_{x}^{#3}(S_{T})}}

\newcommand{\bigenergylocal}[1]{\bignorm{#1}_{L_{t}^{\infty}H^{1}(S_{T})}}

\newcommand{\WaveSobnorm}[3]{\norm{#1}_{H^{#2,#3}\subeps}}

\newcommand{\SecondWaveSobnorm}[3]{\norm{#1}_{\scrH\subeps^{#2,#3}}}

\newcommand{\SecondWaveSobnormArg}[4]{\norm{#1}_{\scrH\subeps^{#2,#3}(#4)}}

\newcommand{\WaveSobnormArg}[4]{\norm{#1}_{H\subeps^{#2,#3}(#4)}}
\newcommand{\bigWaveSobnormArg}[4]{\bignorm{#1}_{H\subeps^{#2,#3}(#4)}}
\newcommand{\WaveSobdotnorm}[3]{\norm{#1}_{\dot{H}\subeps^{#2,#3}}}

\newcommand{\Anorm}[3]{\norm{#1}_{\dot{\mathcal H}\subeps^{#2,#3}}}

\newcommand{\AnormArg}[4]{\norm{#1}_{\dot{\mathcal H}\subeps^{#2,#3}(#4)}}
\newcommand{\bigAnormArg}[4]{\bignorm{#1}_{\dot{\mathcal H}\subeps^{#2,#3}(#4)}}
\newcommand{\Xnorm}[3]{\norm{#1}_{X^{#2,#3}_{\tau = \pm h\subeps(\xi)}}}
\newcommand{\bigXnorm}[3]{\bignorm{#1}_{X^{#2,#3}_{\tau = \pm h\subeps(\xi)}}}
\newcommand{\XnormPlus}[3]{\norm{#1}_{X^{#2,#3}_{\tau = + h\subeps(\xi)}}}
\newcommand{\bigXnormPlus}[3]{\bignorm{#1}_{X^{#2,#3}_{\tau = + h\subeps(\xi)}}}
\newcommand{\XnormMinus}[3]{\norm{#1}_{X^{#2,#3}_{\tau = - h\subeps(\xi)}}}
\newcommand{\bigXnormMinus}[3]{\bignorm{#1}_{X^{#2,#3}_{\tau = - h\subeps(\xi)}}}
\newcommand{\XnormArg}[4]{\norm{#1}_{X^{#2,#3}_{\tau = \pm h\subeps(\xi)}(#4)}}
\newcommand{\bigXnormArg}[4]{\bignorm{#1}_{X^{#2,#3}_{\tau = \pm h\subeps(\xi)}(#4)}}

\newcommand{\A}{\mathbf{A}}

\newcommand{\B}{\mathbf{B}}
\newcommand{\C}{\mathbb{C}}
\newcommand{\D}{\abs{\nabla}}
\newcommand{\E}{\mathbf{E}}
\newcommand{\J}{\mathbf{J}}

\newcommand{\scrH}{\mathcal{H}}

\newcommand{\N}{\mathbb{N}}
\newcommand{\Proj}{\mathcal{P}}
\newcommand{\R}{\mathbb{R}}

\newcommand{\Z}{\mathbb{Z}}

\newcommand{\innerprod}[2]{\left\langle \, #1 , #2 \, \right\rangle}

\newcommand{\angles}[1]{\langle #1 \rangle}

\newcommand{\low}{_{\mathrm{low}}}
\newcommand{\high}{_{\mathrm{high}}}
\newcommand{\fixedpar}{\Lambda}
\newcommand{\curl}{\nabla \times}

\DeclareMathOperator{\diag}{diag}

\DeclareMathOperator{\im}{Im}
\DeclareMathOperator{\re}{Re}
\DeclareMathOperator{\dv}{div}

\newtheorem{theorem}{Theorem}[section]
\newtheorem{proposition}[theorem]{Proposition}
\newtheorem{lemma}[theorem]{Lemma}
\newtheorem{corollary}[theorem]{Corollary}
\newtheorem*{TransferPrinciple}{Transfer Principle}

\theoremstyle{definition}
\newtheorem{definition}[theorem]{Definition}

\theoremstyle{remark}
\newtheorem{remark}[theorem]{Remark}

\numberwithin{equation}{section}

\title{On the asymptotic analysis of the Dirac-Maxwell system in the 
nonrelativistic limit}

\author{Philippe Bechouche, Norbert J. Mauser and Sigmund Selberg\\
Wolfgang Pauli Institute, c/o Inst. f. Math., \\
Universit\"at Wien, \\
Strudlhofgasse 4, A-1090 Wien}

\date{}

\begin{document}

\maketitle

{
 \renewcommand{\thefootnote}{}
}

\begin{abstract}
  We deal with the ``nonrelativistic limit'', i.e.\ the limit $c \to \infty$, 
  where $c$ is the speed of light, of the nonlinear PDE system obtained by
  coupling the Dirac
  equation for a 4-spinor to the Maxwell equations for the self-consistent
  field created by the ``moving charge'' of the spinor.
  This limit, sometimes also called ``Post-Newtonian'' limit, yields a
  Schr\"odinger-Poisson system, where the spin and the magnetic field no longer appear.
  However, our splitting of the 4-spinor into two 2-spinors preserves the symmetry of 
  ``electrons'' and ``positrons''; the latter obeying a Schr\"odinger equation
  with ``negative mass'' in the limit.
  We rigorously prove that in the nonrelativistic limit solutions of the 
  Dirac-Maxwell 
  system on $\R^{1+3}$ converge in the energy space $C([0,T];H^{1})$ to
  solutions of a Schr\"odinger-Poisson system, under appropriate (convergence)
  conditions on the initial data.
  
  We also prove that the time interval of existence of local solutions of 
  Dirac-Maxwell is bounded from below by \emph{log}$(c)$. In fact, for this result we only 
  require uniform $H^{1}$ bounds on the initial data, not convergence.

  Our key technique is ``null form estimates'', extending the work of 
  Klainerman and Machedon and our previous work on the nonrelativistic
  limit of the Klein-Gordon-Maxwell system.
\end{abstract}
\section{Introduction}

In this paper we study the behavior of solutions to the 
Dirac-Maxwell (abbr.\ \emph{DM}) system in the limit $c \to \infty$, 
where $c$ is the speed of light. Coupled to the Coulomb gauge 
condition, this system has the form
\begin{equation}\label{DMC}
  \left( i \gamma^{\mu} \partial_{\mu} - M + g \gamma^{\mu} A_{\mu} 
  \right) \psi = 0,
  \qquad
  \partial^{\nu} F_{\mu\nu} = J_{\mu} / c,
  \qquad
  \partial^{j} A_{j} = 0.
\end{equation}
Here the unknowns are the spinor field $\psi(t,x) \in \C^{4}$, 
regarded as a column vector, and the electromagnetic potential 
$A_{\mu}(t,x) \in \R$, $\mu = 0,1,2,3$. Further, $F_{\mu\nu} = \partial_{\mu} 
A_{\nu} - \partial_{\nu} A_{\mu}$ is the electromagnetic field 
tensor, and
$$
  J^{\mu} = c \innerprod{\gamma^{0}\gamma^{\mu} \psi}{\psi}_{\C^{4}}
$$
is the 4-current density.
On the Minkowski spacetime $\R^{1+3}$ we use relativistic coordinates
$x^{0} = ct \in \R$, $x = ( x^{1}, x^{2}, x^{3} ) \in \R^{3}$.
$\partial_{\mu}$ stands for $\tfrac{\partial}{\partial x^{\mu}}$.
Thus, $\partial_{0} = \tfrac{1}{c} \partial_{t}$, where
$\partial_{t} = \tfrac{\partial}{\partial t}$. We also write
$\nabla = (\partial_{1},\partial_{2},\partial_{3})$,
$\Delta = \partial_{1}^{2} + \partial_{2}^{2} +\partial_{3}^{2}$
and $\D^{s} = (-\Delta)^{s/2}$ for $s \in \R$.
Indices are raised and lowered using the metric $(\eta_{\mu\nu}) = 
\diag(-1,1,1,1)$. The Einstein summation convention is in effect. 
Thus, repeated greek indices $\mu,\nu,\dots$ are summed over $0,1,2,3$,
and repeated roman indices $j,k,\dots$ over $1,2,3$. For example, 
$\Delta = \partial_{j} \partial^{j}$. We denote by 
$\innerprod{\cdot}{\cdot}_{\C^{n}}$ the standard inner product on $\C^{n}$.

The physical constants are $M = m_{0} c / \hbar$, $g = e/\hbar c$,
where $m_{0}$ is the spinor's rest mass, $\hbar$ is the Planck 
constant and $e$ is the unit charge. By $\gamma^{\mu}$, $\mu = 
0,1,2,3$, we denote the $4 \times 4$ Dirac matrices, given in $2 
\times 2$ block form by
$$
  \gamma^{0} = \begin{pmatrix}
    I & 0  \\
    0 & -I
  \end{pmatrix},
  \qquad
  \gamma^{j} = \begin{pmatrix}
    0 & \sigma^{j}  \\
    -\sigma^{j} & 0
  \end{pmatrix},
$$
where the Pauli matrices $\sigma^{j}$ are given by
$$
  \sigma^{1} = 
  \begin{pmatrix}
    0 & 1  \\
    1 & 0
  \end{pmatrix},
  \qquad
  \sigma^{2} = 
  \begin{pmatrix}
    0 & -i  \\
    i & 0
  \end{pmatrix},
  \qquad
  \sigma^{3} = 
  \begin{pmatrix}
    1 & 0  \\
    0 & -1
  \end{pmatrix}.
$$
The following related matrices occur frequently:
$$
  \alpha^{j} := \gamma^{0}\gamma^{j}
  = \begin{pmatrix}
    0 & \sigma^{j}  \\
    \sigma^{j} & 0
  \end{pmatrix},
  \qquad
  S^{m} := i \gamma^{k} \gamma^{l} = 
  \begin{pmatrix}
    \sigma^{m} & 0  \\
    0 & \sigma^{m}
  \end{pmatrix},
$$
where $(k,l,m)$ is any cyclic permutation of $(1,2,3)$. Note the 
identities
\begin{equation}\label{alphaIdentities}
  \alpha^{j} \alpha^{k} = - \alpha^{k} \alpha^{j} + 2 \delta^{jk} I
  = \delta^{jk} I + i \epsilon^{jkl} S_{l}.
\end{equation}

The first equation in \eqref{DMC} is the Dirac equation. Multiplying 
it on the left by $\gamma^{0}$ and taking the imaginary part of its 
$\C^{4}$ inner product with $\psi$ yields the conservation law 
$\partial_{\mu} J^{\mu} = 0$.
Thus, the ``charge'' is conserved:
\begin{equation}\label{ChargeConservation}
  \int \innerprod{\psi}{\psi}_{\C^{4}} \, dx
  = \twonorm{\psi(t)}{}^{2} = \text{\emph{const.}}
\end{equation}

The second equation in \eqref{DMC} is the Maxwell equation. We split 
$A_{\mu}$ into its temporal part $A_{0}$, the electric potential, and 
its spatial part $\A = (A^{1},A^{2},A^{3})$, the magnetic potential.
Hence the electric field is given by $\E = \nabla A_{0} - 
\partial_{0} \A$ and the magnetic field by $\B = \curl \A$, and the 
second equation in \eqref{DMC} is seen to be equivalent to the 
Maxwell system in classical form, with charge density $\rho = J^{0}/c$
and current density $(J^{k})_{k=1,2,3}$.

The third equation in \eqref{DMC} is the Coulomb gauge 
condition $\dv \A = 0$. The reason for this choice of gauge condition will be 
explained later. It is equivalent to $\Proj \A = \A$, where 
$\Proj$ is the projection onto divergence free vector fields 
in $\R_{x}^{3}$. The second and third equations in \eqref{DMC} are 
then seen to be equivalent to
$$
  \Delta A_{0} = J^{0}/c, \qquad
  \left( c^{-2} \partial_{t}^{2} - \Delta \right)
  \A = c^{-1} \Proj (J^{k})_{k=1,2,3}
$$
provided the initial data of 
$\A$ are divergence free. Thus, when properly rescaled
(see \cite{BMP}, \cite{MM}), the system 
\eqref{DMC} is conveniently expressed in terms of a small
dimensionless parameter $$\varepsilon \simeq \frac{1}{c}$$ as follows:
\begin{subequations}\label{DMCscaled}
\begin{align}
  \label{DiracEq}
  i \partial_{t} \psi\eps &= - i \varepsilon^{-1} \alpha^{j} 
  \partial_{j} \psi\eps + \varepsilon^{-2} \gamma^{0} \psi\eps
  - A_{j}\eps \alpha^{j} \psi\eps - A_{0}\eps\psi\eps,
  \\
  \label{A0Eq}
  \Delta A_{0}\eps &= \rho\eps,
  \\
  \label{AEq}
  \square\subeps \A\eps &= \varepsilon \Proj \J\eps,
\end{align}
\end{subequations}
where we have put in superscripts to emphasize the dependence on
$\varepsilon$. Here
$$
  \square\subeps = \varepsilon^{2} \partial_{t}^{2} - \Delta
$$
and
\begin{equation}\label{ChargeCurrent}
  \rho\eps = \innerprod{\psi\eps}{\psi\eps}_{\C^{4}},
  \qquad
  \J\eps = \varepsilon^{-1} \left\{  
  \innerprod{\alpha^{k}\psi\eps}{\psi\eps}_{\C^{4}} \right\}_{k=1,2,3}.
\end{equation}
We consider the Cauchy problem for \eqref{DMCscaled} with ``finite 
energy'' initial data
\begin{equation}\label{Data}
  \psi\eps \init = \psi_{0}\eps \in H^{1},
  \qquad
  (\A\eps, \partial_{t} \A\eps) \init
  = (\mathbf a_{0}\eps, \mathbf a_{1}\eps)
  \in
  \Proj \dot H^{1} \times \Proj L^{2}.
\end{equation}

We prove three types of results for this system as $\varepsilon \to 0$.
First, local well-posedness (abbr.\ l.w.p.) with a logarithmic lower bound
on the existence time.
Second, convergence in the nonrelativistic limit if the initial datum 
of $\psi$ converges. Third, we prove some more precise results on the 
asymptotic behavior of the Dirac spinor under various smallness assumptions on
its ``positron part''. These results are described in detail in the next three subsections. 

\subsection{Local existence}

There are two issues here: (i) l.w.p.\ for $\varepsilon$ fixed, and (ii) the nature of the
$\varepsilon$-dependence of the local existence time as $\varepsilon \to 0$.

Concerning (i), the main difficulty  is that one cannot directly estimate
the bilinear term $A_{j} \alpha^{j} \psi$ in the Dirac equation, due to the failure of
the endpoint Strichartz estimate for the wave equation in $1+3$ 
dimensions. The crucial fact proved here is that when
the Dirac equation is squared,
the bilinear terms resulting from this dangerous term
can all be expressed in terms of null bilinear forms, provided
the Coulomb gauge condition is used, and this enables us to prove l.w.p.\ of
\emph{DM} in the energy space \eqref{Data}, a result entirely analogous to
that of Klainerman and Machedon \cite{KlMa94} for the Klein-Gordon-Maxwell
(\emph{KGM}) system. (The square of the Dirac eq.\ is similar to the Klein-Gordon
eq., but contains some additional bilinear terms due to the presence of spin.)

Bournaveas \cite{Bourn} proved l.w.p.\ of \emph{DM} in the space
$(\psi(t),\A(t)) \in H^{1/2+\delta} \times H^{1+\delta}$ for $\delta > 0$, but this
result does not take into account the null structure; in fact, by using
the special structure of the equations and the so-called Wave-Sobolev spaces,
the result can be improved to\footnote{Even this is not optimal (the scale invariant space
is $L^2 \times H^{1/2}$), and in view of the recent work of Machedon and Sterbenz
\cite{MaSt} on \emph{KGM} one may indeed hope to do better.}
$(\psi(t),\A(t)) \in H^{s} \times H^{1}$ for $1/4 < s \le 1$; this is proved in
an upcoming paper by the third author. It is worth pointing out that these results
are all independent w.r.t.\ $s$, since the regularity of $\A(t)$ is kept fixed.
Thus, e.g., the l.w.p.\ in $H^{1/2} \times H^1$ does not imply l.w.p.\ in
$H^1 \times H^1$ (or vice versa).

The question of global existence and uniqueness for \emph{DM} remains largely 
open\footnote{This in contrast to the situation for \emph{KGM}; see \cite{KlMa94}.
The crucial point is that \emph{KGM} has a positive Hamiltonian, unlike \emph{DM}.}
(but see Georgiev \cite{Ge} for a small data result), however, we prove---and this brings
us to the second issue mentioned above---that as $\varepsilon \to 0$ the
local existence time goes to infinity, subject to the initial assumptions
\begin{equation}\label{DataBound}
  \Sobnorm{\psi_{0}\eps}{1} = O(1),
  \quad
  \Sobdotnorm{\mathbf a_{0}\eps}{1}
  + \varepsilon \twonorm{\mathbf a_{1}\eps}{} = 
  O\left(\frac{1}{\varepsilon^{\fixedpar}}\right)
  \quad \text{as} \quad \varepsilon \longrightarrow 0,
\end{equation}
where
\begin{equation}\label{GrowthParameter}
  0 < \fixedpar < \frac{1}{2}
\end{equation}
will be kept fixed throughout the paper.
(The upper bound $1/2$ is explained by the 
factor $\varepsilon^{1/2}$ appearing in the $L^{2}$ bilinear
estimates discussed in Sect.\ \ref{BilinearSection}.)

\begin{theorem}\label{LWPThm} \textbf{($H^{1}$ l.w.p.\ of DM.)}
The initial value problem (\ref{DMCscaled}), (\ref{Data})
is locally well posed for fixed $\varepsilon$,
with an existence time $T\subeps > 0$
depending only on $\varepsilon$ and the size of the norms of the data. Moreover, if
\eqref{DataBound} holds, then
\begin{equation}\label{Lifespan}
  T\subeps \ge c_{0} \log \frac{1}{\varepsilon}
  \quad \text{as} \quad
  \varepsilon \longrightarrow 0,
\end{equation}
where $c_{0} > 0$ is a universal constant, and we have
\begin{equation}\label{EvolutionBound}
  \Sobnorm{\psi\eps(t)}{1}  = O(1),
  \quad
  \Sobdotnorm{\A\eps(t)}{1}
  + \varepsilon \twonorm{\partial_{t} \A\eps(t)}{}
  = O\left(\frac{1}{\varepsilon^{\fixedpar}}\right)
\end{equation}
uniformly in every finite time interval as $\varepsilon \to 0$.
\end{theorem}

In order to control the evolution as $\varepsilon \to 0$, it is
crucial to have estimates which are sufficiently strong w.r.t.\
powers of $\varepsilon$. To this end we employ analytical techniques 
used in our earlier paper \cite{BMS}, where the nonrelativistic limit
of \emph{KGM} was considered. The analysis of \emph{DM} 
is more involved, however, due to the additional terms that come up in
``squared Dirac'' compared to the usual Klein-Gordon (\emph{KG}) 
equation. In particular, we prove some new bilinear spacetime estimates which are
needed to control these extra terms.

Once we have obtained closed estimates for the system---sufficiently
strong w.r.t.\ powers of $\varepsilon$---we use a bootstrap argument
to prove existence in a short time interval
\emph{depending only on the $L^{2}$ norm of $\psi_{0}\eps$}, provided 
$\varepsilon$ is sufficiently small, depending on the size of 
\eqref{DataBound}.
On account of the conservation of charge \eqref{ChargeConservation}
for the Dirac equation we can then iterate this argument to obtain 
the long time result.

We stress the fact that no convergence assumption is made on the data in the 
above theorem---all we need is the uniform bound \eqref{DataBound}.
However, if we \emph{do} assume that $\psi_{0}\eps$
converges in $H^{1}$, then we can pass to the nonrelativistic limit,
which we discuss next.

\subsection{Nonrelativistic limit}

The nonrelativistic limit of the linear Dirac equation with a given
time-dependent electromagnetic potential was treated in \cite{BMP}
(earlier papers, see e.g.\ \cite{CC}, dealt only with the static case,
i.e.\ time-independent potential).
There are also some results on the nonlinear Dirac and Klein-Gordon equations
in the literature, see e.g.\ \cite{Naj1}, but for the coupled nonlinear
Dirac-Maxwell and Klein-Gordon-Maxwell systems there are no results previous
to our work (i.e.\ the present paper as well as \cite{BMS, BMSproc})
and the completely independent work of Masmoudi and Nakanishi \cite{MN}.

The most marked difference between
our work and that of Masmoudi and Nakanishi is that our estimates are strong enough to give
uniform (w.r.t.\ $\varepsilon$) bounds for the solutions of \emph{DM}
assuming only the initial boundedness condition \eqref{DataBound}---no
convergence assumption is necessary. This, of course, is crucial as far as
proving Theorems \ref{LWPThm} and \ref{SeminonrelLimitThm} is concerned.
By contrast, in \cite{MN} the convergence
assumption is essential  because uniform estimates are obtained only on an
arbitrarily small time interval, and in order to push the result to a larger time
interval they must use bounds on solutions of the limiting Schr\"odinger-Poisson
system.

Let us now state our result.
We split the Dirac spinor into its upper and lower components:
\begin{equation}\label{UpperLower}
  \psi\eps = 
  \begin{pmatrix}
    \UpperSpinor\eps  \\
    \LowerSpinor\eps
  \end{pmatrix}
\end{equation}
where $\UpperSpinor$ and $\LowerSpinor$ are 2-spinors, i.e.\ column
vectors in $\C^{2}$.
Before one can pass to the limit $\varepsilon \to 0$, the rest energy 
must be subtracted, which for the upper ``positive energy'' component means
multiplication by $e^{it/\varepsilon^{2}}$ and for the lower ``negative energy''
component multiplication by $e^{-it/\varepsilon^{2}}$.

\begin{theorem}\label{NonrelLimitThm} \textbf{(Nonrelativistic limit of DM.)}
Consider the solution of (\ref{DMCscaled}), (\ref{Data}) obtained in 
Theorem \ref{LWPThm}, with data satisfying:
\begin{enumerate}
  \item\label{ConvAssumption} $v_{0} := \lim_{\varepsilon \to 0} \psi_{0}\eps$
  exists in $H^{1}$,
  \item $\Sobdotnorm{\mathbf a_{0}\eps}{1}
  + \varepsilon \twonorm{\mathbf a_{1}\eps}{} = 
  O\left(\frac{1}{\varepsilon^{\fixedpar}}\right)$
  as $\varepsilon \to 0$.
\end{enumerate}
Denote the upper and lower 2-spinors of $v_{0}$ by $v_{0}^{+}$ and 
$v_{0}^{-}$ respectively, and let $(u,v_{+},v_{-})$
be the solution of the Schr\"odinger-Poisson system\footnote{This system
is globally well posed for $L^{2}$ data.}
\begin{equation}\label{PS}
  \Delta u = n,
  \qquad n = \abs{v_{+}}^{2} + \abs{v_{-}}^{2},
  \qquad
  \left( i\partial_{t} \pm \frac{\Delta}{2} \right) v_{\pm}
  + u v_{\pm} = 0,
\end{equation}
with initial data $v_{\pm} \init = v_{0}^{\pm}$.
Then as $\varepsilon \to 0$,
\begin{subequations}\label{Convergence}
\begin{alignat}{2}
  \label{SpinorConv}
  \psi\eps &= e^{-it/\varepsilon^{2}} 
  \begin{pmatrix}
    v_{+}  \\
    0
  \end{pmatrix}
  + e^{+it/\varepsilon^{2}} 
  \begin{pmatrix}
    0  \\
    v_{-}
  \end{pmatrix}
  + o(1) &\qquad & 
  \text{in} \quad H^{1}, 
  \\
  \label{A0Conv}
  A_{0}\eps &= u + o(1) & \qquad  & \text{in} \quad \dot 
  H^{1},
  \\
  \label{ChargeConv}
  \rho\eps &= n + o(1) &\qquad & \text{in} \quad L^{p}, \,\,
  1 \le p \le 3,
\end{alignat}
\end{subequations}
uniformly in every finite time interval. Moreover, the relativistic 
current density converges as follows: Let
\begin{equation}\label{J0}
\begin{split}
  \J^{0} &=
  \im \innerprod{\nabla v_{+}}{v_{+}}_{\C^{2}}
  - \im \innerprod{\nabla v_{-}}{v_{-}}_{\C^{2}}
  \\
  &\quad
  + \frac{1}{2} \curl \innerprod{\vec \sigma v_{+}}{v_{+}}_{\C^{2}}
  - \frac{1}{2} \curl \innerprod{\vec \sigma v_{-}}{v_{-}}_{\C^{2}}
\end{split}
\end{equation}
where $\innerprod{\nabla v_{\pm}}{v_{\pm}}$ and
$\innerprod{\vec \sigma v_{\pm}}{v_{\pm}}$ are the 
vectors with components $\innerprod{\partial^{j} v_{\pm}}{v_{\pm}}$
and $\innerprod{\sigma^{j} v_{\pm}}{v_{\pm}}$, respectively,
for $j = 1,2,3$. Then
\begin{equation}\label{CurrentConv}
  \J\eps \longrightarrow \J^{0}
  \quad \text{in}
  \quad
  \left[ C^{1}_{c}(\R_{t} \times \R_{x}^{3}) \right]'
  \quad \text{weak $*$}
\end{equation}
as $\varepsilon \to 0$.
\end{theorem}

The first line in r.h.s.\eqref{J0} 
is the conserved current associated to the limiting system \eqref{PS},
whereas the second line consists of the well-known divergence-free additional terms
due to the interaction spin-magnetic field \cite{LL3}.

We can improve the convergence rate to $O(\varepsilon)$ by strengthening 
the initial assumptions. Thus, we shall prove:

\begin{theorem}\label{NonrelLimitThm2}
Strengthen the hypotheses of Theorem \ref{NonrelLimitThm} by assuming
$$
  \Sobnorm{\psi_{0}\eps}{2} = O(1),
  \quad
  \Sobnorm{\nabla \mathbf a_{0}\eps}{1}
  + \varepsilon \Sobnorm{\mathbf a_{1}\eps}{1} = 
  O\left(\frac{1}{\varepsilon^{\fixedpar}}\right)
$$
and
$$
  \psi\eps_{0}  = 
  \begin{pmatrix}
    v^{+}_{0}  \\
    0
  \end{pmatrix}
  + O(\varepsilon)
  \qquad \text{in} \quad H^{1}
$$
as $\varepsilon \to 0$. Moreover, assume $v_{0}^{+} \in H^{5}$.
Then
\begin{equation}\label{ImprovedConvergence}
  \psi\eps  = e^{-it/\varepsilon^{2}} 
  \begin{pmatrix}
    v_{+}  \\
    0
  \end{pmatrix}
  + O(\varepsilon)
  \qquad \text{in} \quad H^{1} \quad \text{as} \quad \varepsilon
  \longrightarrow 0 
\end{equation}
uniformly in every finite time interval. Furthermore, the convergence in (\ref{A0Conv})
and (\ref{ChargeConv}) is also $O(\varepsilon)$.
\end{theorem}

\begin{remark} The hypotheses are not strong enough to guarantee 
strong convergence of the current density $\J\eps$ locally uniformly in 
time. In fact, a simple counterexample is given by the initial datum
$$
  \psi\eps_{0} = 
  \begin{pmatrix}
    v^{+}_{0}  \\
    \varepsilon v^{+}_{0}
  \end{pmatrix}.
$$
Then $\J\eps$ initially has vector components $2 \re 
\innerprod{\sigma^{j} v^{+}_{0}}{v^{+}_{0}}$, which does not agree 
with the weak limit $\J^{0}$ given by \eqref{J0}.
\end{remark}

It is instructive to compare the last theorem to the formal derivation of the
nonrelativistic limit usually reproduced in physics 
textbooks, the basic premise of which is a smallness assumption on the 
lower component $\LowerSpinor\eps$ of the spinor.
The idea is to 
define\footnote{Here we break the symmetry of the signs in 
\eqref{SpinorConv}, i.e.\ between ``electrons'' and ``positrons'',
but this is not important since the lower component
is in any case expected to vanish.}
\begin{equation}\label{UpperLower2}
  \phi\eps = 
  \begin{pmatrix}
    \chi\eps  \\
    \eta\eps
  \end{pmatrix}
  := e^{it/\varepsilon^{2}} \psi\eps.
\end{equation}
Then \eqref{ImprovedConvergence} can be restated
\begin{equation}\label{ImprovedConvergence2}
  \chi\eps  = v_{+} + O(\varepsilon),
  \qquad
  \eta\eps  = O(\varepsilon)
  \qquad \text{in} \quad H^{1} \quad \text{as} \quad \varepsilon
  \longrightarrow 0.
\end{equation}
The Dirac equation \eqref{DiracEq} gives
\begin{equation}\label{UpperLowerDirac}
  iD_{0} \chi\eps = - i \sigma^{j} D_{j} \eta\eps,
  \qquad
  iD_{0} \eta\eps = - i \sigma^{j} D_{j} \chi\eps - 
  \frac{2}{\varepsilon} \eta\eps,
\end{equation}
where we write $D_{0} = \varepsilon \partial_{t} - i \varepsilon 
A_{0}\eps$ and $D_{j} = \partial_{j} - i \varepsilon A_{j}\eps$.
Thus,
\begin{equation}\label{LowerComponentExpansion}
  \eta\eps = - \varepsilon \frac{1}{2} i \sigma^{j} \partial_{j} \chi\eps
  - \varepsilon^{2} \frac{1}{2} \left\{ i \partial_{t} \eta\eps + 
  A_{0}\eps \eta\eps + A_{j}\eps \sigma^{j} \chi\eps \right\},
\end{equation}
and substituting this in the first equation in \eqref{UpperLowerDirac} gives, after some algebra,
\begin{equation}\label{ApproxPauli}
  i \partial_{t} \chi\eps = \frac{1}{2} \left( i \nabla + 
  \varepsilon \A\eps \right)^{2} \chi\eps - A_{0}\eps \chi\eps
  - \frac{1}{2} \varepsilon B_{j}\eps \sigma^{j} \chi\eps
  - \varepsilon r\eps,
\end{equation}
where
\begin{equation}\label{PauliRemainder}
  r\eps = \frac{1}{2} \sigma^{j} D_{j} \left(\partial_{t} \eta\eps - 
  i A_{0}\eps \eta\eps \right)
\end{equation}
and $\B\eps = \curl \A\eps$. Then by formal considerations of 
magnitude, in particular assuming $\partial_{t} \eta\eps = O(1)$, one 
obtains a Schr\"odinger equation in the limit $\varepsilon \to 0$.
It is possible to make this argument rigorous, but it has a
fundamental weakness which limits its usefulness, namely that 
$\partial_{t} \eta\eps$ can be no better than $O(1/\varepsilon)$ unless 
one adds a further constraint on the initial data. In fact, it is 
clear from \eqref{LowerComponentExpansion} that $\partial_{t} \eta\eps = 
O(1)$ in $L^{2}$ initially if and only if the constraint
\begin{equation}\label{InitialConstraint}
  \eta\eps = - \varepsilon \frac{1}{2} i \sigma^{j} \partial_{j} \chi\eps
  + O(\varepsilon^{2})
\end{equation}
holds in $L^{2}$ at time $t = 0$, assuming the data \eqref{Data} are $O(1)$.

However, the constraint \eqref{InitialConstraint} is not needed in Theorem 
\ref{NonrelLimitThm2}, the reason being that instead of the simple 
splitting into upper and lower components as in \eqref{UpperLower}, we 
apply the eigenspace projections of the ``free Dirac operator''
$$\mathcal Q\eps = -i \varepsilon \alpha^{k} \partial_{k}
+ \gamma^{0}.$$ As in \cite{BMP} we use the spectral decomposition
$$
  \mathcal Q\eps = \lambda\eps \Pi_{+}\eps - \lambda\eps \Pi_{-}\eps
$$
where
\begin{equation}\label{Eigenvalues}
  \lambda\eps = \sqrt{1 - \varepsilon^{2} \Delta},
  \qquad
  \Pi_{\pm}\eps = \frac{1}{2} \left( I \pm [\lambda\eps]^{-1} 
  \mathcal Q\eps \right).
\end{equation}
Since the positive and negative eigenvalues $\pm \lambda\eps$ 
correspond to positive and negative energies of a free Dirac 
particle, the spectral decomposition is related to electrons and 
positrons (\cite{D}).
The formal limit $\varepsilon \to 0$ of 
$\Pi_{\pm}\eps$ yields the operators
\begin{equation}\label{ZeroOrderProj}
  \Pi_{\pm}^{0} = \frac{1}{2} (I \pm \gamma^{0}),
  \qquad
  \Pi_{+}^{0} =
  \begin{pmatrix}
    I & 0  \\
    0 & 0
  \end{pmatrix},
  \qquad
  \Pi_{-}^{0} =
  \begin{pmatrix}
    0 & 0  \\
    0 & I
  \end{pmatrix}.
\end{equation}

The following basic lemma shows that $\Pi_{\pm}^{0}$ is the leading 
order term in a series expansion of $\Pi_{\pm}\eps$ in powers of 
$\varepsilon$, and moreover that \eqref{InitialConstraint} 
is basically equivalent to $\Pi_{-}\eps \psi\eps = 
O(\varepsilon^{2})$, a condition which resurfaces in the next subsection.

\begin{lemma}\label{DiracProjLemma} For all $s \in \R$,
$\Pi_{\pm}\eps$ is bounded from $H^{s} \to H^{s}$ uniformly 
in $\varepsilon$. Moreover,
\begin{align}
  \label{ProjExpansion1}
  \Pi_{\pm}\eps &= \Pi_{\pm}^{0} \mp \varepsilon \mathcal R_{1}\eps
  \\
  \label{ProjExpansion2}
  &= \Pi_{\pm}^{0} \mp
  i \varepsilon \frac{1}{2} \alpha^{k} \partial_{k} \mp
  \varepsilon^{2} \mathcal R_{2}\eps
\end{align}
where $\mathcal R_{j}\eps$ denotes an operator bounded from 
$H^{s} \to H^{s-j}$ uniformly in $\varepsilon$.
\end{lemma}

\begin{proof} This follows immediately from
\begin{equation}\label{ProjExpansion3}
  \Pi_{\pm}\eps - \Pi_{\pm}^{0} = \mp \frac{1}{2} [\lambda\eps]^{-1} i \varepsilon
  \alpha^{k} \partial_{k}
  \mp \frac{1}{2} \left( 1 - [\lambda\eps]^{-1} \right) \gamma^{0}
\end{equation}
and the fact that the Fourier symbol
of $1 - [\lambda\eps]^{-1}$ satisfies the inequalities
\begin{equation}\label{SymbolBound}
  0 \le 1 - \frac{1}{\sqrt{1 + \varepsilon^{2} \abs{\xi}^{2}}}
  \le \min \left\{ 1, \varepsilon \abs{\xi}, \varepsilon^{2} 
  \abs{\xi}^{2} \right\}
\end{equation}
where $\xi$ is the Fourier variable corresponding to $x$.
\end{proof}
  
Before moving on, we prove that the initial data
assumption \eqref{ConvAssumption} in Theorem \ref{NonrelLimitThm}
implies something stronger, namely the convergence of $\Pi_{\pm}\eps \psi_{0}\eps$.

\begin{lemma}\label{DataLemma}
If
$$
  \lim_{\varepsilon \to 0} \psi_{0}\eps = v_{0} =
  \begin{pmatrix}
    v_{0}^{+}  \\
    v_{0}^{-}
  \end{pmatrix}
$$
exists in $H^{1}$, then
$$
  \lim_{\varepsilon \to 0} \Pi_{+}\eps \psi_{0}\eps = 
  \begin{pmatrix}
    v_{0}^{+}  \\
    0
  \end{pmatrix}
  \quad \text{and} \quad
  \lim_{\varepsilon \to 0} \Pi_{-}\eps \psi_{0}\eps = 
  \begin{pmatrix}
    0 \\
    v_{0}^{-}
  \end{pmatrix}
$$
in $H^{1}$.
\end{lemma}

\begin{proof} It suffices to prove
$\left( \Pi_{\pm}\eps - \Pi_{\pm}^{0} \right) \psi_{0}\eps \to 0$ in $H^{1}$.
But the proof of Lemma \ref{DiracProjLemma} shows that the Fourier symbol
of $\Pi_{\pm}\eps - \Pi_{\pm}^{0}$ is bounded in absolute value 
by $\min\{ 1 , \varepsilon \abs{\xi} \}$.
Thus $\left( \Pi_{\pm}\eps - \Pi_{\pm}^{0} \right) (\psi_{0}\eps - v_{0}) \to 0$
in $H^{1}$, and by Plancherel's theorem and dominated convergence,
$\left( \Pi_{\pm}\eps - \Pi_{\pm}^{0} \right) v_{0} \to 0$ in $H^{1}$.
\end{proof}

\subsection{Semi-nonrelativistic limit}

As in \cite{BMP}, by the ``semi-nonrelativistic limit'' we understand
the approximation of the upper component of the Dirac equation
by the Pauli equation for a 2-spinor, which reads
\begin{subequations}\label{PauliEq}
\begin{equation}
  i \partial_{t} \Pauli\eps = \frac{1}{2} \left( i \nabla + 
  \varepsilon \A\eps \right)^{2} \Pauli\eps - A_{0}\eps \Pauli\eps
  - \frac{1}{2} \varepsilon B_{j}\eps \sigma^{j} \Pauli\eps
\end{equation}
with initial condition
\begin{equation}
  \Pauli\eps \init = \PauliData\eps \in H^{1}.
\end{equation}
\end{subequations}
Note that the naive ``upper and lower components'' approach in 
\eqref{UpperLowerDirac}--\eqref{PauliRemainder} can give at best an 
$O(\varepsilon)$ approximation to the Pauli equation, assuming 
the initial constraint \eqref{InitialConstraint}, which as remarked is 
essentially equivalent to $\Pi_{-}\eps \psi\eps = O(\varepsilon^{2})$.

In contrast, by using the Dirac projections $\Pi_{\pm}\eps$ instead of 
just $\Pi_{\pm}^{0}$, we can prove an $O(\varepsilon^{2})$ approximation,
with the same initial constraint. In fact, we have the following 
result:

\begin{theorem}\label{SeminonrelLimitThm}
Consider the solution of (\ref{DMCscaled}), (\ref{Data}) obtained in 
Theorem \ref{LWPThm}. Define $\chi\eps$ as in (\ref{UpperLower2}) and let
$\Pauli\eps$ be the solution of the Pauli equation (\ref{PauliEq}). 
Assume the initial conditions
\begin{enumerate}
  \item
  $\Sobnorm{\psi_{0}\eps}{5} = O(1), \quad
  \Sobnorm{\nabla \mathbf a_{0}\eps}{4}
  + \varepsilon \Sobnorm{\mathbf a_{1}\eps}{4} = O(1)$,
  \item
  $\Sobnorm{\Pi_{-}\eps \psi_{0}\eps}{1} = O(\varepsilon^{2})$,
\end{enumerate}
as $\varepsilon \to 0$.  Then if
\begin{equation}\label{PauliConv}
  \Sobnorm{\chi\eps - \Pauli\eps}{1} = O(\varepsilon^{2})
\end{equation}
holds at time $t = 0$, it also holds uniformly in every finite time 
interval. For the current density we then have
\begin{equation}\label{SeminonrelCurrent}
  \J\eps = \J_{P}\eps + \frac{1}{2} \curl \innerprod{\vec \sigma 
  \Pauli\eps}{\Pauli\eps}_{\C^{2}}
  + O(\varepsilon)
  \quad \text{in}
  \quad L_{x}^{1},
\end{equation}
uniformly in every finite time interval, where
\begin{equation}\label{PauliCurrent}
  \J_{P}\eps = \im \innerprod{ (\nabla - i\varepsilon \A\eps) \Pauli\eps}{\Pauli\eps}_{\C^{2}}
\end{equation}
is the current density of the Pauli equation.
\end{theorem}

The remainder of this paper is organized as follows: In the following 
section we square the Dirac equation and reinterpret it in terms of 
the projections $\Pi_{\pm}\eps \psi\eps$ of the spinor, and we prove that 
the main bilinear terms can be expressed in terms of null forms. Then in Sect.\ 
\ref{BilinearSection} we discuss the linear and bilinear spacetime 
estimates of Strichartz type that are used in this paper. The proofs of 
those estimates that are not already in the literature can be found in Sect.\ 
\ref{BilinearProofs}. In Sect.\ \ref{FunctionSpaces} we define the 
function spaces that we use, and recall their main properties.
The main estimates for the nonlinear terms are proved in Sect.\ 
\ref{MainEstimates}, which is the heart of the paper. Then in Sects.\ 
\ref{UniformBounds}--\ref{SeminonrelLimit} these estimates are applied 
to prove the main theorems.

To close this section we introduce some notational conventions
which will be in effect throughout:
\begin{itemize}
  \item
  For function spaces we use the following notation.
  If $X$ is a Banach space of functions on $\R_{x}^{3}$, we denote by 
  $L_{t}^{p} X$ the space with norm
  $$
    \norm{u}_{L_{t}^{p}X}
    = \left( \int_{-\infty}^{\infty} \norm{u(t,\cdot)}_{X}^{p} \, dt 
    \right)^{1/p},
  $$
  with the usual modification if $p = \infty$. The localization of this norm 
  to a time slab $S_{T} = [0,T] \times \R^{3}$ is denoted
  $\norm{u}_{L_{t}^{p}X(S_{T})}$.
  \item
  In estimates, we use the notation $\lesssim$ to mean $\le$
  up to multiplication by a positive constant $C$
  independent of $\varepsilon$. Moreover, in estimates over a time slab
  $S_{T}$, $C$ is also understood to be independent of $T$.
  \item
  For exponents, we use the convenient shorthand $p^{+}$ (resp. $p^{-}$)
  for $p + \zeta$ (resp. $p - \zeta$) with $\zeta > 
  0$ sufficiently small, independently of $\varepsilon$.
  The notation $\infty^{-}$ stands for a sufficiently large,
  positive exponent.
  \item
  We denote by $f(x) \mapsto \widehat f(\xi)$
  and $u(t,x) \mapsto \widehat u(\tau,\xi)$ the Fourier 
  transforms on $\R^{3}$ and $\R^{1+3}$, respectively.
  As in \cite{BMS} we split functions $f(x)$ into their
  low ($\abs{\xi} \lesssim 1/\varepsilon$)
  and high ($\abs{\xi} \gtrsim 1/\varepsilon$) frequency parts,
  \begin{equation}\label{LowHigh}
    f = f\low + f\high,
  \end{equation}
  corresponding to a smooth partition of unity in Fourier space.

\end{itemize}

\section{Preliminaries}\label{Preliminaries}

As already mentioned, our approach to the Dirac equation 
is to square it and apply techniques similar to those
used for \emph{KGM} in \cite{BMS}. It is therefore convenient to work with
the ``\emph{KG} splitting''
\begin{equation}\label{KGsplitting}
  \psi_{\pm}\eps
  =
  \begin{pmatrix}
    \UpperSpinor\eps_{\pm}  \\
    \LowerSpinor\eps_{\pm}
  \end{pmatrix}
  := \frac{1}{2} \left\{ \psi\eps \pm 
  \varepsilon^{2} [\lambda\eps]^{-1} \bigl( i \partial_{t} \psi\eps + A_{0}\eps 
  \psi\eps \bigr) \right\}
\end{equation}
as used in \cite{BMS}.
In order to compare this to the Dirac projections \eqref{Eigenvalues}, 
observe that if  $\psi\eps$ solves the Dirac equation \eqref{DiracEq}, then
\begin{equation}\label{DiracvsKGsplitting}
  \psi_{\pm}\eps = \Pi_{\pm}\eps \psi\eps \mp \frac{1}{2}
  \varepsilon^{2} [\lambda\eps]^{-1} (A_{j}\eps \alpha^{j} \psi\eps).
\end{equation}
But using the estimate
\begin{equation}\label{LambdaEst}
  \Sobnorm{[\lambda\eps]^{-1}f}{\sigma}
  \le \varepsilon^{-r} \Sobnorm{f}{\sigma-r}
  \quad \text{for}
  \quad 0 \le r \le 1,
\end{equation}
followed by H\"older's inequality and Sobolev embedding, we see that
\begin{equation}\label{ErrorTerm}
  \varepsilon^{2}
  \Sobnorm{[\lambda\eps]^{-1} (A_{j}\eps \alpha^{j} \psi\eps)}{\sigma}
  \lesssim \varepsilon^{2-\sigma}
  \Sobdotnorm{\A\eps}{1} \Sobnorm{\psi\eps}{1}
  \quad \text{for} \quad 0 \le \sigma \le 1,
\end{equation}
so r.h.s.\eqref{DiracvsKGsplitting} is $O(\varepsilon^{1-\fixedpar})$
in $H^{1}$ at time $t$ if the bound \eqref{EvolutionBound} in
Theorem \ref{LWPThm} holds.
As far as proving Theorem \ref{NonrelLimitThm} is concerned, it is therefore
immaterial whether we use $\psi_{\pm}\eps$ or $\Pi_{\pm}\eps \psi\eps$.

For later use we note the following consequences of \eqref{DiracvsKGsplitting} 
and \eqref{ErrorTerm}. First,
\begin{equation}
  \label{KGsplittingBound}
  \Sobnorm{\psi_{\pm}\eps}{\sigma}
  \lesssim
  \Sobnorm{\psi\eps}{\sigma}
  + \varepsilon^{2-\sigma}\Sobdotnorm{\A\eps}{1} \Sobnorm{\psi\eps}{1}
  \quad \text{for} \quad 0 \le \sigma \le 1,
\end{equation}
using the uniform boundedness of $\Pi_{\pm}\eps$. Second,
\begin{equation}\label{SmallComponentsBound}
  \Lxpnorm{\UpperSpinor_{-}\eps}{2} + 
  \Lxpnorm{\LowerSpinor_{+}\eps}{2}
  \lesssim
  \varepsilon \Sobnorm{\psi\eps}{1}
  + \varepsilon^{2}
  \Sobdotnorm{\A\eps}{1} \Sobnorm{\psi\eps}{1},
\end{equation}
where we used \eqref{ProjExpansion1} and the 
orthogonality of $\Pi_{+}^{0}$ and $\Pi_{-}^{0}$.

\medskip\medskip\medskip
\noindent
Let us now restate the system \eqref{DMCscaled} in terms of the 
splitting \eqref{KGsplitting} of the spinor. 
First we subtract the rest energy, defining
\begin{equation}\label{ModulatedFields}
  \phi_{\pm}\eps
  =
  \begin{pmatrix}
    \chi_{\pm}\eps  \\
    \eta_{\pm}\eps
  \end{pmatrix}
  := e^{\pm it/\varepsilon^{2}} \psi_{\pm}\eps.
\end{equation}
Thus
\begin{equation}\label{PsiExpansion}
  \psi\eps = \psi_{+}\eps + \psi_{-}\eps
  = e^{-it/\varepsilon^{2}} \phi_{+}\eps
  + e^{+it/\varepsilon^{2}} \phi_{-}\eps.
\end{equation}

\begin{lemma}\label{ModifiedDiracLemma}
In terms of the splitting \eqref{PsiExpansion}, defined via 
\eqref{KGsplitting} and \eqref{ModulatedFields}, the Dirac equation
(\ref{DiracEq}) is equivalent to a system of two equations
\begin{equation}\label{ModifiedDirac}
  L_{+}\eps \phi_{+}\eps = - A_{0}\eps \phi_{+}\eps + 
  \frac{1}{2} e^{it/\varepsilon^{2}} R\eps,
  \quad L_{-}\eps \phi_{-}\eps = - A_{0}\eps \phi_{-}\eps -
  \frac{1}{2} e^{- it/\varepsilon^{2}} R\eps,
\end{equation}
provided the constraint \eqref{DiracvsKGsplitting} is satisfied at 
time $t = 0$, or equivalently that the Dirac equation is satisfied at
$t = 0$. Here
\begin{equation}\label{LDef}
  L_{\pm}\eps = i\partial_{t}
  \mp \frac{\lambda\eps - 1}{\varepsilon^{2}}
\end{equation}
and  $R\eps$ is given by
\begin{multline}\label{RDef}
  \lambda\eps R\eps =
  \varepsilon \left\{ 2i \A\eps \cdot \nabla
  + i \dv \A\eps
  + i E_{j}\eps \alpha^{j}
  - B_{j}\eps S^{j} \right\} \psi\eps
  \\
  + \varepsilon^{2} (\A\eps)^{2}\psi\eps
  - [A_{0},\lambda\eps](\psi_{+}\eps - \psi_{-}\eps).
\end{multline}
Further, $[\cdot,\cdot]$ denotes the commutator and
\begin{equation}\label{EandB}
  \E = (E^{1},E^{2},E^{3}) := \nabla A_{0} - \varepsilon \partial_{t} \A,
  \qquad
  \B = (B^{1},B^{2},B^{3}) := \curl \A.
\end{equation}
\end{lemma}

\begin{proof} Squaring the Dirac equation \eqref{DiracEq} yields 
(cf.\ \cite[Sect.\ 70]{D})
\begin{equation}\label{SquaredDirac}
  \left\{ \varepsilon^{2} (i\partial_{t} +  A_{0}\eps)^{2}
  + (\nabla - i\varepsilon \A\eps)^{2}
  - \varepsilon^{-2} - i\varepsilon E_{j}\eps \alpha^{j}
  + \varepsilon B_{j}\eps S^{j} \right\} \psi\eps = 0.
\end{equation}
Applying $i\partial_{t} + A_{0}\eps$ to both sides of
\eqref{KGsplitting} and making use of \eqref{SquaredDirac}
and
$$
  \varepsilon^{2} (i\partial_{t} + A_{0}\eps)\psi\eps
  = \lambda\eps (\psi_{+}\eps - \psi_{-}\eps),
$$
which follows from \eqref{KGsplitting},
one easily obtains \eqref{ModifiedDirac}. Reversing these steps, one 
finds that \eqref{ModifiedDirac} implies the squared Dirac 
equation \eqref{SquaredDirac}. But the latter implies the Dirac 
equation, since we assume that \eqref{DiracvsKGsplitting} holds 
initially, which amounts to saying that the Dirac equation is 
satisfied initially.
\end{proof}

Let us make a brief, heuristic comparison of \eqref{ModifiedDirac} with
the expected limit \eqref{PS}. As it turns out, $R\eps$
vanishes in the limit, so \eqref{ModifiedDirac} tends to
the Schr\"odinger equation in \eqref{PS}.
In fact, the Fourier symbol 
of $(\lambda\eps - 1)/ \varepsilon^{2}$ is
\begin{equation}\label{heps}
  h\subeps(\xi) := \frac{\abs{\xi}^{2}}{1 + \sqrt{1 + \varepsilon^{2}\abs{\xi}^{2}}}
  \sim
  \begin{cases}
    \abs{\xi}^{2}/2 &\quad \text{for $\abs{\xi} \lesssim 1/\varepsilon$},
    \\
    \abs{\xi}/\varepsilon &\quad \text{for $\abs{\xi} \gtrsim 1/\varepsilon$},
  \end{cases}
\end{equation}
so $L_{\pm}$ tends to the Schr\"odinger operator
$i\partial_{t} \pm \Delta/2$ as $\varepsilon \to 0$.
Moreover, the charge and current densities (\ref{ChargeCurrent}) are 
given in terms of the fields \eqref{ModulatedFields} by
\begin{align}
  \label{ChargeExpansion}
  \rho\eps &= \innerprod{\chi_{+}\eps}{\chi_{+}\eps}
  + \innerprod{\chi_{-}\eps}{\chi_{-}\eps}
  + \innerprod{\eta_{+}\eps}{\eta_{+}\eps}
  + \innerprod{\eta_{-}\eps}{\eta_{-}\eps}
  \\
  \notag
  &\qquad\qquad
  + 2 \re \left\{ e^{-2it/\varepsilon^{2}}
  \innerprod{\chi_{+}\eps}{\chi_{-}\eps}
  +  e^{-2it/\varepsilon^{2}} \innerprod{\eta_{+}\eps}{\eta_{-}\eps} 
  \right\},
  \\
  \label{CurrentExpansion}
  \mathbf{J}\eps &= \frac{2}{\varepsilon} \re \Bigl\{ 
  \innerprod{\sigma^{j} \chi_{+}\eps}{\eta_{+}\eps}
  + \innerprod{\sigma^{j} \chi_{-}\eps}{\eta_{-}\eps}
  \\
  \notag
  &\qquad\qquad
  + e^{2it/\varepsilon^{2}} \innerprod{\sigma^{j} \chi_{-}\eps}{\eta_{+}\eps} 
  + e^{-2it/\varepsilon^{2}} \innerprod{\sigma^{j} \chi_{+}\eps}{\eta_{-}\eps} 
  \Bigr\}_{j = 1,2,3}.
\end{align}
We expect [cf.\ \eqref{SmallComponentsBound}] that $\chi_{-}\eps,
\eta_{+}\eps \to 0$. Thus, in
r.h.s.\eqref{ChargeExpansion} only the first and fourth terms are of
importance, and $\Delta A_{0}\eps = \rho\eps$ tends to
the Poisson equation in \eqref{PS}.

For later use we note the estimate
\begin{equation}\label{hEst}
  0 \le \abs{\xi}/\varepsilon - h\subeps(\xi) \lesssim \varepsilon^{-2}.
\end{equation}
This reduces to $r - \alpha(r) \lesssim 1$, where
\begin{equation}\label{AlphaDef}
  \alpha(r) := \frac{r^{2}}{ 1 + \sqrt{1+ r^{2}}}.
\end{equation}
But $r - \alpha(r) = r+1-s = 1 - \frac{1}{r+s}$,
where $s = \sqrt{1+r^{2}}$.

\medskip\medskip\medskip
\noindent
We now turn to the problem of obtaining closed estimates for the 
modified \emph{DM}
system \eqref{ModifiedDirac}, \eqref{A0Eq}, \eqref{AEq}.
A serious obstacle to estimating the bilinear
terms in \eqref{RDef} is the failure
of the endpoint Strichartz estimate for the wave equation in dimension $1+3$.
The salient feature of the Coulomb gauge, however, is that these 
problematic terms can be expressed in terms of the null bilinear forms
\begin{equation}\label{NullForms}
  Q_{0}(u,v) = \partial_{0} u \partial_{0} v
  - \nabla u \cdot \nabla v,
  \quad
  Q_{\alpha\beta}(u,v) = \partial_{\alpha} u \partial_{\beta} v
  - \partial_{\beta} u \partial_{\alpha} v,
\end{equation}
where $\partial_{0}$ denotes $\varepsilon \partial_{t}$
and $0 \le \alpha < \beta \le 3$.
These bilinear forms enjoy better regularity properties than generic products of 
derivatives.

We emphasize that in the following result $\psi$ does not have to
solve the Dirac equation.

\begin{lemma}\label{NullLemma} \textbf{(Null structure.)}
Given a potential $\{A_{\mu}(t,x)\}$
satisfying the Coulomb condition $\dv \A = 0$, let $\E$ and $\B$ be 
defined as in \eqref{EandB}, and consider the bilinear operator
$$
  \psi
  \longrightarrow
  \left\{ 2i \A \cdot \nabla
  + i E_{j} \alpha^{j}
  - B_{j} S^{j} \right\} \psi
$$
appearing in (\ref{RDef}). We have the following identities:
\begin{equation}\label{NullIdentity1}
  2 \A \cdot \nabla \psi
  = - Q_{jk}( \D^{-1} a^{jk}, \psi)
\end{equation}
and
\begin{equation}\label{NullIdentity2}
\begin{split}
  &\left\{i (E_{j} - \partial_{j} A_{0}) \alpha^{j}
  - B_{j} S^{j} \right\} \psi
  \\
  &\qquad =  
  Q_{jk}(\D^{-1} \varepsilon \partial_{t} a^{jk}, U)
  - Q_{jk}(\D^{-1} \partial_{l} a^{jk}, \alpha^{l} U)
  \\
  &\qquad\quad
  + Q_{0}(A_{j},\alpha^{j} U)
  + Q_{0j}(A_{k}, \alpha^{j}\alpha^{k} U)
  - \frac{i}{2} Q_{jk}(A_{m}, \epsilon^{jkl} S_{l}\alpha^{m} U)
\end{split}
\end{equation}
where
$$
  a_{jk} = R_{j} A_{k} - R_{k} A_{j},
  \qquad R_{j} = \D^{-1}\partial_{j}
$$
and $U = U(\psi)$ is the 4-spinor defined by
\begin{equation}\label{UDef}
  \square\subeps U = - i \left( \varepsilon \partial_{t}
  + \alpha^{j} \partial_{j} \right) \psi,
  \qquad
  U \init = 0,
  \qquad i \varepsilon \partial_{t} U \init = \psi_{0}.
\end{equation}
Here $\psi_{0}$ denotes $\psi \init$.
\end{lemma}

\begin{proof} The identity \eqref{NullIdentity1} goes back to the work 
of Klainerman and Machedon \cite{KlMa94} on \emph{KGM},
so we concentrate on the new identity \eqref{NullIdentity2}.
Define
$$
  \partial_{\pm} = \varepsilon \partial_{t} \pm \alpha^{j} \partial_{j}
$$
and observe that
\begin{equation}\label{NullLemmaA}
  (\partial_{-} A_{j}) \alpha^{j} =
  - (E_{j} - \partial_{j} A_{0}) \alpha^{j}
  - i B_{j} S^{j},
\end{equation}
where we used the second identity in \eqref{alphaIdentities}
and the assumption $\dv \A = 0$.
By the first identity in \eqref{alphaIdentities},
$$
  \partial_{+}\partial_{-} = \square\subeps.
$$
Thus \eqref{UDef} implies that $w = \psi - i\partial_{-} U$ satisfies
$\partial_{+} w = 0$ with $w \init = 0$, whence
$$
  \psi = i\partial_{-} U.
$$
Apply \eqref{NullLemmaA} to this and use
$$
  \partial_{+}(\alpha^{j} U) = \alpha^{j} \partial_{-} U + 
  2 \partial^{j} U
  \qquad (\text{by \eqref{alphaIdentities}})
$$
to rewrite l.h.s.\eqref{NullIdentity2} as
$$
  (\partial_{-} A_{j}) \partial_{+}(\alpha^{j} U)
  - 2(\partial_{-} A_{j}) \partial^{j}U.
$$
To the last term we apply the identity \eqref{NullIdentity1}; this we 
can do since $\dv \partial_{\mu} \A = 0$ for $\mu = 0,1,2,3$ by assumption.
To the first term we apply the following general
formula, obtained using the second identity in \eqref{alphaIdentities},
$$
  (\partial_{-} \phi)(\partial_{+} U)
  = Q_{0}(\phi,U) + Q_{0j}(\phi,\alpha^{j} U)
  - \frac{i}{2} Q_{jk}(\phi,\epsilon^{jkl}S_{l}U),
$$
where $\phi$ is a function and $U$ a 4-spinor. This last formula is 
due to Klainerman and Machedon \cite{KlMaPC}.
\end{proof}

\section{Bilinear spacetime estimates}\label{BilinearSection}

The main technical tools used in this paper are spacetime estimates
of Strichartz type for solutions of the free initial value problems
\begin{equation}\label{uvEq}
\begin{split}
  \square\subeps u &= 0, \qquad u \init = f, \qquad \partial_{t} u 
  \init = 0,
  \\
  L_{\pm}\eps v & = 0, \qquad v \init = g,
\end{split}
\end{equation}
on $\R^{1+3}$.
Let us first describe the new $L^{2}$ product estimates that are 
proved in this paper, and then we recall the estimates proved in 
\cite{BMS}.

Let $\mu$ and $\lambda$ be dyadic numbers of the form $2^{j}$, $j 
\in \Z$. Denote by $\LPproj_{\mu}$ the Littlewood-Paley operator given by
$$
  \left( \LPproj_{\mu} f \right)\FT(\xi) = \beta(\xi/\mu) \widehat f(\xi),
$$
where $\beta$ is a bump function supported in $\abs{\xi} \sim 
1$ such that $\sum_{j \in \Z} \beta(\xi/2^{j}) = 1$ for $\xi \neq 0$.
We write $f_{\mu} = \LPproj_{\mu} f$ and similarly for $g, u, v$.
Thus $f = \sum_{\mu} f_{\mu}$ etc.

We shall prove the following:

\begin{theorem}\label{DyadicThm} The solutions $u, v$ of (\ref{uvEq})
satisfy the following dyadic spacetime estimates:
\begin{enumerate}
  \item\label{DiagonalLow}
  $\twonorm{\LPproj_{\mu}(u_{\lambda} v_{\lambda})}{_{t,x}}
  \lesssim \varepsilon^{1/2} \mu \twonorm{f_{\lambda}}{}
  \twonorm{g_{\lambda}}{}
  \quad \text{if} \quad \mu \lesssim \lambda \lesssim 
  1/\varepsilon$.
  \item\label{DiagonalHigh}
  $\twonorm{\LPproj_{\mu}(u_{\lambda} v_{\lambda})}{_{t,x}}
  \lesssim  \varepsilon^{1/2} \mu^{1/2} \lambda^{1/2}
  \twonorm{f_{\lambda}}{}
  \twonorm{g_{\lambda}}{}
  \quad \text{if} \quad \mu \lesssim \lambda,\,\,
  \lambda \gg 1/\varepsilon$.
  \item\label{OffDiagonal}
  $\twonorm{u_{\mu} v_{\lambda}}{_{t,x}}
  \lesssim \varepsilon^{1/2} \min( \mu, \lambda )
  \twonorm{f_{\mu}}{}
  \twonorm{g_{\lambda}}{}
  \quad \text{for all} \quad \mu, \lambda$.
\end{enumerate}
\end{theorem}

See \cite[Thm.\ 12.1]{FK} for the analogous estimates in the case where
$u$ and $v$ both solve the wave equation.

By decomposing the product $uv$ into dyadic pieces, then applying Theorem 
\ref{DyadicThm} and finally exploiting the orthogonality properties in Fourier
space to sum up, one obtains the following corollary. (The complete argument 
can be found in \cite[Sect.\ 12]{FK}.)

\begin{corollary}\label{DyadicCor1}
The solutions $u, v$ of (\ref{uvEq}) satisfy
$$
  \bigtwonorm{\D^{-\sigma} (uv)}{_{t,x}} \le C_{s_{1},s_{2}} \,
  \varepsilon^{1/2}
  \Sobdotnorm{f}{s_{1}} \Sobdotnorm{g}{s_{2}}
$$
provided that
$$
  s_{1}, s_{2} < 1,
  \qquad
  \sigma < \frac{1}{2},
  \qquad
  s_{1}+s_{2}+\sigma = 1.
$$
\end{corollary}

Estimates of this type for the case where
$u$ and $v$ both solve the free wave equation were first investigated by 
Klainerman and Machedon. The case $(s_{1},s_{2},\sigma)
= (0,1,0)$ is excluded, a fact related to the false 
endpoint case of the Strichartz estimates for the wave equation
in $1+3$ dimensions. However, by assuming a little extra regularity 
one can easily sum the dyadic pieces and one obtains
the following nonsharp bilinear estimate.

\begin{corollary}\label{DyadicCor2}
The solutions $u, v$ of (\ref{uvEq}) satisfy
$$
  \twonorm{uv}{_{t,x}} \le C_{\delta} \, \varepsilon^{1/2}
  \twonorm{f}{} \Sobnorm{g}{1+\delta}
$$
for all $\delta > 0$.
\end{corollary}

\begin{proof}
It suffices to prove the sharp estimate
\begin{equation}\label{BilinearEndpoint}
  \twonorm{u v_{\lambda}}{_{t,x}} \lesssim \varepsilon^{1/2}
  \twonorm{f}{} \lambda \twonorm{g_{\lambda}}{}.
\end{equation}
Write $u = \sum_{\mu} u_{\mu}$ and consider the cases 
$\mu \lesssim \lambda$ and $\mu \gg \lambda$. In the first case,
$$
  \twonorm{\left(\sum\nolimits_{\mu \lesssim \lambda} u_{\mu}\right) v_{\lambda}}{}
  \lesssim \sum\nolimits_{\mu \lesssim \lambda} \twonorm{u_{\mu} v_{\lambda}}{}
  \lesssim \varepsilon^{1/2} \left( \sum\nolimits_{\mu \lesssim 
  \lambda} \frac{\mu}{\lambda} \twonorm{f_{\mu}}{} \right)
  \lambda \twonorm{g_{\lambda}}{},
$$
where we used Theorem \ref{DyadicThm}(iii) to get the last inequality. In the second case we 
have, by orthogonality in Fourier space,
$$
  \twonorm{\left(\sum\nolimits_{\mu \gg \lambda} u_{\mu}\right) v_{\lambda}}{}^{2}
  \lesssim \sum\nolimits_{\mu \gg \lambda} \twonorm{u_{\mu} v_{\lambda}}{}^{2},
$$
and by Theorem \ref{DyadicThm}(iii) we dominate this by
$\varepsilon \twonorm{f}{}^{2} \lambda^{2} \twonorm{g_{\lambda}}{}^{2}$.
\end{proof}

Here we could also take $f$ in $H^{1+\delta}$ and $g \in L^{2}$, but 
we shall not need this. However, for null bilinear forms one can get
the sharp result (i.e. $\delta = 0$).
Thus, we recall the following, proved in \cite[Proposition 4]{BMS}:
\begin{equation}\label{NullFormEstA}
  \twonorm{Q_{ij}\bigl(u,v\bigr)}{_{t,x}}
  \lesssim \varepsilon^{1/2} \Sobdotnorm{f}{2} \Sobdotnorm{g}{1}
\end{equation}
where $Q_{ij}$ is given by \eqref{NullForms}.
We remark that this is the analogoue of an estimate for two solutions of the
free wave equation proved by Klainerman and Machedon. 

\medskip\medskip\medskip
\noindent
Since we will prove part \eqref{DiagonalHigh} of Theorem \ref{DyadicThm} 
by a reduction to linear Strichartz estimates, let us recall these 
(for $1+3$ dimensions). We say that a pair $(q,r)$ of Lebesgue exponents is
\emph{wave admissible} if $(q,r) \neq (2,\infty)$ and $1/q + 1/r \le 1/2$,
and \emph{sharp wave admissible} if the last inequality is an equality.

For the free wave $u$ in \eqref{uvEq} one has the well-known estimate
\begin{equation}\label{LinearWaveStr}
  \mixednorm{u}{q}{r} \le C_{q,r} \,
  \varepsilon^{1/q} \Sobdotnorm{f}{s},
\end{equation}
for wave admissible $(q,r)$ and $s = 3/2 - 3/r - 1/q$.
As proved in \cite{KlTa}, this can be improved if the Fourier support 
of $f$ is small. Thus, if $\widehat f$ is supported in a cube with 
side length $\sim \mu$ and at distance $\sim \lambda$ from the origin, 
where $\mu \ll \lambda$, then
\begin{equation}\label{ImprovedWaveStr}
  \mixednorm{u}{q}{r} \le C_{q,r} \, \varepsilon^{1/q}
  \left( \frac{\mu}{\lambda} \right)^{1/2-1/r}
  \Sobdotnorm{f}{s},
\end{equation}
for $q,r,s$ as above.

For $v$ satisfying \eqref{uvEq} we have, as proved in 
\cite[Proposition 1]{BMS},
\begin{equation}\label{LinearStr}
  \mixednorm{v}{q}{r} \le C_{q,r} \left( \Sobdotnorm{g\low}{1/q} + 
  \varepsilon^{1/q} \Sobdotnorm{g\high}{2/q} \right)
\end{equation}
for sharp wave admissible $(q,r)$. (Then one can use Sobolev embedding 
to obtain estimates for all wave admissible pairs.)
In order to prove Theorem \ref{DyadicThm}\eqref{DiagonalHigh} we need 
the analogue of \eqref{ImprovedWaveStr} in this context. Thus, we 
shall prove:

\begin{proposition}\label{ImprovedStrProp}
Let $v$ be as in (\ref{uvEq}), and suppose $\widehat g$ is supported in
a cube with side length $\sim \mu$ and at distance $\sim \lambda$ from the origin, 
where $\mu \ll \lambda$. Then
\begin{equation}\label{ImprovedStr}
  \mixednorm{v}{q}{r} \le C_{q,r} \left( \frac{\mu}{\lambda} 
  \right)^{1/2-1/r}
  \left( \Sobdotnorm{g\low}{1/q} + 
  \varepsilon^{1/q} \Sobdotnorm{g\high}{2/q} \right)
\end{equation}
for sharp wave admissible $(q,r)$.
\end{proposition}

Finally, recalling the basic heuristic that $L_{\pm}\eps$ behaves
like a Schr\"odinger operator at low frequencies,
it is not surprising that we have the following Schr\"odinger type 
estimates, proved in \cite{BMS}. We say that a pair
$(q,r)$ is \emph{Schr\"odinger admissible} if $q,r \ge 2$ and $2/q + 3/r = 3/2$.

\begin{proposition}\label{SchrodingerProp}
Let $(q,r)$ and $(\widetilde q, \widetilde r)$ be any two 
Schr\"odinger admissible pairs. Then for the solution of 
$L\eps_{\pm} v = F$ with data $v \init = f$ we have
$$
  \mixednormlocal{v\low}{q}{r} + \mixednormlocal{v\low}{\infty}{2}
  \lesssim \twonorm{f\low}{} + \mixednormlocal{F\low}{\widetilde 
  q'}{\widetilde r'},
$$
where $\tfrac{1}{\widetilde q} + \tfrac{1}{\widetilde q'} = 1$ and
$\tfrac{1}{\widetilde r} + \tfrac{1}{\widetilde r'} = 1$.
\end{proposition}

\section{Function spaces}\label{FunctionSpaces}

We shall use the following spaces of functions on $\R^{1+3}$
with weighted norms defined in Fourier space:
\begin{itemize}
  \item
  $H^{s,\theta}\subeps$ with norm $\twonorm{ \angles{\xi}^{s} 
  \angles{\abs{\tau} - \varepsilon^{-1} \abs{\xi}}^{\theta} \widehat 
  u(\tau,\xi)}{_{\tau,\xi}}$.
  \item
  $\dot H^{s,\theta}\subeps$ with norm $\twonorm{ \abs{\xi}^{s} 
  \angles{\abs{\tau} - \varepsilon^{-1} \abs{\xi}}^{\theta} \widehat 
  u(\tau,\xi)}{_{\tau,\xi}}$.
  \item
  $\scrH^{s,\theta}\subeps$ with norm $\norm{u}_{ 
  H^{s,\theta}\subeps}
  + \varepsilon \norm{\partial_{t} u}_{H^{s-1,\theta}\subeps}$.
  \item
  $\dot \scrH^{1,\theta}\subeps$ with norm $\norm{u}_{\dot 
  H^{1,\theta}\subeps}
  + \varepsilon \norm{\partial_{t} u}_{H^{0,\theta}\subeps}$.
  \item
  $X^{s,\theta}_{\tau = \pm h\subeps(\xi)}$ with norm
  $\twonorm{ \angles{\xi}^{s} 
  \angles{\tau \mp h\subeps(\xi)}^{\theta} \widehat 
  u(\tau,\xi)}{_{\tau,\xi}}$ and $h\subeps$ as in \eqref{heps}.
\end{itemize}
Here $\angles{\cdot}$ stands for $1 + \abs{\cdot}$.
These spaces are by now standard, and we will recall their main 
properties without proofs. For more details and 
further references to the literature, the reader may consult e.g.\ 
\cite{Tao}, \cite{KlSe}.

It will be convenient to introduce the notation
\begin{equation}\label{Propagators}
\begin{split}
  U\eps(t) &= e^{it( \lambda\eps - 1) / \varepsilon^{2}} = 
  e^{ith\subeps(\D)},
  \\
  S(t) &= e^{it\Delta/2},
  \\
  W\eps(t) &= e^{it \D / \varepsilon }
\end{split}
\end{equation}
for the propagators associated to, respectively, the operators 
$L\eps_{\pm}$ defined in Lemma \ref{ModifiedDiracLemma}, the Schr\"odinger 
operator and the wave operator. 

\medskip\medskip\medskip
\noindent
(i) \emph{Superposition principle.}
A fundamental property of the so-called ``Wave Sobolev space'' 
$H^{s,\theta}\subeps$ is that any function in this space
can be written as a superposition
($H^{s}$-valued integral over the real line) of solutions of the
free wave equation with initial data in $H^{s}$.
(See \cite[Proposition 3.4]{KlSe} for the 
precise statement.) This, in effect, replaces Duhamel's 
principle in the framework of the Wave Sobolev spaces, and it has the
following simple but extremely useful consequence (see \cite{KlSe}):

\begin{TransferPrinciple} Suppose $T$ is a multilinear operator 
$(f_{1}(x),\dots,f_{k}(x)) \mapsto T(f_{1},\dots,f_{k})(x)$ acting in 
$x$-space. If $T$ satisfies an estimate
$$
  \mixednorm{T\bigl(W\eps(\pm t) f_{1},\dots, W\eps(\pm t) f_{k}\bigr)}{q}{r}
  \le C \varepsilon^{1/q} \Sobnorm{f_{1}}{s_{1}} \cdots \Sobnorm{f_{k}}{s_{k}},
$$
for all combinations of signs, then
$$
  \mixednorm{T\left(u_{1},\dots, u_{k} \right)}{q}{r}
  \le C_{\theta} \varepsilon^{1/q} 
  \WaveSobnorm{u_{1}}{s_{1}}{\theta} \cdots \WaveSobnorm{u_{k}}{s_{k}}{\theta}
$$
holds for all $u_{j} \in H\subeps^{s_{j},\theta}$, provided $\theta > 1/2$.
Moreover, the same statement holds with $H^{s}$ and $H\subeps^{s_{j},\theta}$ 
replaced by their homogeneous counterparts.
\end{TransferPrinciple}

The spaces $X^{s,\theta}_{\tau = \pm h\subeps(\xi)}$ are related to the 
equation $L\eps_{\pm} v = 0$ in the same way that the Wave 
Sobolev spaces are related to the free wave equation. Thus, we have a 
superposition principle and hence a transfer principle for these 
spaces as well. To be precise, in the above Transfer Principle, 
one can replace any one of the $W\eps(\pm t)$ by
$U\eps(\pm t)$ and correspondingly $H\subeps^{s_{j},\theta}$
by $X^{s_{j},\theta}_{\tau = \pm h\subeps(\xi)}$.
Applying this to estimates from the previous section, we have, for $\theta > 1/2$,
\begin{alignat}{2}
  \label{WaveEmb}
  &\mixednorm{u\high}{q}{r}
  \lesssim \varepsilon^{1/2 + 1/r} 
  \Xnorm{u\high}{1}{\theta}
  &\quad &\text{for sharp wave adm.\ $(q,r)$},
  \\
  \label{SchrEmb}
  &\mixednorm{u\low}{q}{r}
  \lesssim
  \Xnorm{u\low}{0}{\theta}
  &\quad &\text{for Schr\"odinger adm.\ $(q,r)$},
  \\
  \label{SchrEmb2}
  &\mixednorm{u\low}{2}{\infty}
  \lesssim
  \Xnorm{u\low}{1}{\theta}.
  & &
\end{alignat}
Here \eqref{SchrEmb} follows from Proposition \ref{SchrodingerProp}
with $F = 0$. By Sobolev embedding we reduce \eqref{SchrEmb2}
to the case $(q,r) = (2,6)$ of \eqref{SchrEmb}. Finally,
\eqref{WaveEmb} holds by virtue of \eqref{LinearStr} and 
the trivial estimate
\begin{equation}\label{HighFreqEst}
  \Sobnorm{f\high}{s} \lesssim \varepsilon^{\sigma} 
  \Sobnorm{f\high}{s+\sigma}
  \quad \text{for} \quad \sigma > 0.
\end{equation}

\medskip\medskip\medskip
\noindent
(ii) \emph{Embeddings.} The most basic embeddings are
\begin{equation}\label{BasicEmbedding}
  H^{s,\theta}\subeps,
  \,\, X^{s,\theta}_{\tau = \pm h\subeps(\xi)}
  \hookrightarrow C_{b}(\R;H^{s}),
  \qquad \dot H^{s,\theta}\subeps
  \hookrightarrow C_{b}(\R;\dot H^{s}),
\end{equation}
which hold \emph{uniformly in $\varepsilon$} for any $\theta > 1/2$.
Also uniform in $\varepsilon$ are
\begin{equation}\label{InterpolationEmbeddings}
  \mixed{p}{2} \hookrightarrow H^{0,\theta-1}\subeps,
  \,\, X^{0,\theta-1}_{\tau = \pm h\subeps(\xi)}
  \quad \text{for} \quad
  \frac{1}{\frac{3}{2}-\theta} < p \le 2, \quad
  \frac{1}{2} < \theta < 1.
\end{equation}
In fact, the dual statement $H^{0,1-\theta}\subeps$,
$X^{0,1-\theta}_{\tau = \pm h\subeps(\xi)}
\hookrightarrow \mixed{p'}{2}$ follows by interpolation between the
trivial case $p = 2$, $\theta = 1$ and \eqref{BasicEmbedding}.
We shall also need
\begin{equation}\label{XtoWaveSobEst}
  \bigXnorm{e^{\pm it/\varepsilon^{2}} u}{0}{\theta-1}
  \lesssim \varepsilon^{-2(1-\theta)} \WaveSobnorm{u}{0}{\theta-1}.
\end{equation}
This is obvious if $\widehat u(\tau,\xi)$ is supported in 
$\bigabs{ \abs{\tau} - \abs{\xi}/\varepsilon}
\lesssim \varepsilon^{-2}$; then we can in fact replace the left 
hand side by $\twonorm{u}{}$. On the other hand, if $\widehat u(\tau,\xi)$
is supported in $\bigabs{ \abs{\tau} - \abs{\xi}/\varepsilon }
\gg \varepsilon^{-2}$, then \eqref{XtoWaveSobEst} follows from 
\eqref{hEst}.

\medskip\medskip\medskip
\noindent
(iii) \emph{Time cut-off.}
In view of \eqref{BasicEmbedding}, we can localize to any finite time slab
$$S_{T} = [0,T] \times \R^{3}.$$
The restriction space $H^{s,\theta}\subeps(S_{T})$ is complete when 
equipped with the norm
\begin{equation}\label{RestrictionNorm}
  \norm{u}_{H^{s,\theta}\subeps(S_{T})}
  := \inf \left\{ \WaveSobnorm{v}{s}{\theta} : \text{$v = u$ on 
  $S_{T}$} \right\}.
\end{equation}
Norms on the other restriction spaces are similarly defined.
When $\theta \le 1/2$ the embeddings \eqref{BasicEmbedding} fail, but 
since $H^{s,\theta}\subeps$ etc.\ are spaces of tempered distributions, 
it still makes sense to restrict them to the \emph{interior} of
$S_{T}$, and we will use the same notation $H^{s,\theta}\subeps(S_{T})$
etc.\ for these spaces. Taking the \emph{inf} over all extensions 
produces a seminorm in this case.

The idea behind the following 
``cut-off lemmas'' originates in the work of Bourgain \cite{B} on
the Schr\" odinger and KdV equations, and was developed further by
Kenig-Ponce-Vega \cite{KPV} in their work on KdV
and by Klainerman-Machedon \cite{KlMa} and the last author \cite{Se2}
for the wave equation.
In fact, the argument given in \cite{KPV} applies to $X^{s,\theta}$ 
spaces in general, and in particular proves the following.

\begin{lemma}\label{XCutoffLemma}
Suppose $L\eps_{\pm} v = F$ on the interior of $S_{T}$
with $v \init = f$. Let $\theta > 1/2$. Then for $0 \le T \le 1$,
$$
  \XnormArg{v}{s}{\theta}{S_{T}} \le C_{\theta} \left( \Sobnorm{f}{s} + 
  \XnormArg{F}{s}{\theta-1}{S_{T}} \right)
$$
where $C_{\theta}$ is independent of $T$ and $\varepsilon$.
\end{lemma}

By rescaling $x \to \varepsilon x$ we reduce the next result to the case 
$\varepsilon = 1$, which in turn follows from estimates proved 
in \cite{KlMa}.

\begin{lemma}\label{HomWaveCutoffLemma}
Suppose $\square\subeps u = F$ on the interior of $S_{T}$ with
$(u,\partial_{t} u) \init = (f,g)$. Let $\theta > 1/2$. Then
for $0 \le T \le 1$,
\begin{equation}\label{CutoffA}
  \AnormArg{u}{1}{\theta}{S_{T}}
  \le C_{\theta} \left( \Sobdotnorm{f}{1} + \varepsilon \twonorm{g}{}
  + \frac{1}{\varepsilon} \WaveSobnormArg{F}{0}{\theta-1}{S_{T}} \right)
\end{equation}
where $C_{\theta}$ is independent of $T$ and $\varepsilon$.
Also, for $M \in \N$ large enough,
\begin{equation}\label{CutoffB}
  \SecondWaveSobnormArg{u}{s}{\theta}{S_{T}}
  \le C_{\theta} \varepsilon^{-M} \left( \Sobnorm{f}{s}
  + \Sobnorm{g}{s-1}
  + \SecondWaveSobnormArg{F}{s-1}{\theta-1}{S_{T}} \right).
\end{equation}
\end{lemma}

The last inequality is not sharp w.r.t.\ $\varepsilon$, but it will 
only be used in a situation where powers of $\varepsilon$ are not important.
In order to estimate the Dirac current density
we shall need the following ``integration by 
parts''-version of Lemma \ref{HomWaveCutoffLemma}.

\begin{lemma}\label{IntByPartsCutoffLemma}
Suppose $\square\subeps u = e^{it/\varepsilon^{2}}F$ on the interior 
of $S_{T}$ with vanishing data. Let $1/2 < \theta < 1$. Then
$$
  \AnormArg{u}{1}{\theta}{S_{T}}
  \lesssim
  \varepsilon \twonorm{\partial_{t} F}{(S_{T})}
  + \LpHs{F\ext}{2}{1}
  + \varepsilon \twonorm{\angles{\partial_{t}}^{\theta} F\ext}{_{t,x}}
  + \frac{\twonorm{F\ext}{_{t,x}}}{\varepsilon^{2\theta-1}}
$$
for all $0 \le T \le 1$ and all extensions $F\ext$ of $F$ 
to all of $\R^{1+3}$. Here 
$\angles{\partial_{t}}^{\theta}$ is the multiplier with Fourier symbol
$(1+\abs{\tau})^{\theta}$.
\end{lemma}

\begin{proof}
Let us denote $F\ext$ simply by $F$. Write
$$
  e^{it/\varepsilon^{2}}F
  = (\varepsilon^{2}/i) \, \partial_{t}
  \left[ e^{it/\varepsilon^{2}}F \right]
  - (\varepsilon^{2}/i) e^{it/\varepsilon^{2}} \partial_{t} F
$$
and $u = u_{1} + u_{2}$ accordingly. By \eqref{CutoffA},
$\AnormArg{u_{2}}{1}{\theta}{S_{T}}
\lesssim \varepsilon
\bigtwonorm{\partial_{t} F}{(S_{T})}.$
Now define $G = e^{it/\varepsilon^{2}}F$.
Split $G = G_{1} + G_{2}$ by a partition of unity in Fourier space such that
\begin{alignat*}{2}
  &\widehat G_{1}(\tau,\xi) &\quad  & \text{is supported in} \quad 
  \abs{\tau} \lesssim \abs{\xi}/\varepsilon,
  \\
  &\widehat G_{2}(\tau,\xi) &\quad  & \text{is supported in} \quad 
  \abs{\tau} \gg \abs{\xi}/\varepsilon,
\end{alignat*}
and write $u_{1} = u_{1,1} + u_{1,2}$ accordingly. That is, $\square\subeps 
u_{1,j} = (\varepsilon^{2}/i) \, \partial_{t} G_{j}$ on the interior 
of $S_{T}$ with vanishing data. By \eqref{CutoffA},
$\AnormArg{u_{1,1}}{1}{\theta}{S_{T}}
\lesssim \, \varepsilon
\twonorm{\partial_{t}G_{1}}{}$,
but using Plancherel's theorem and the assumptions on the Fourier 
support,
$$
  \twonorm{\partial_{t}G_{1}}{}
  \lesssim \frac{1}{\varepsilon} \bigtwonorm{\D G_{1}}{}
  \le \frac{1}{\varepsilon} \bigtwonorm{\D G}{}
  = \frac{1}{\varepsilon} \bigtwonorm{\D F}{},
$$
whence $\AnormArg{u_{1,1}}{1}{\theta}{S_{T}} \lesssim
\LpHs{F}{2}{1}$. Finally, to estimate $u_{1,2}$ we
first observe that it has an extension to all of $\R^{1+3}$ defined in 
Fourier space by
$$
  \widehat u_{1,2}(\tau,\xi)
  = \frac{1}{- \varepsilon^{2} \tau^{2} + \abs{\xi^{2}}}
  \left[ (\varepsilon^{2}/i) \, \partial_{t} G_{2} 
  \right]\FT(\tau,\xi).
$$
Thus
$\abs{ \widehat u_{1,2}(\tau,\xi) } \sim \frac{1}{\abs{\tau}}
\bigabs{\widehat G_{2}(\tau,\xi)}$,
and since
$$
  \AnormArg{u_{1,2}}{1}{\theta}{S_{T}}
  \le \AnormArg{u_{1,2}}{1}{\theta}{\R^{1+3}}
  \lesssim
  \twonorm{ (\varepsilon \abs{\tau} + \abs{\xi})
  \angles{ \abs{\tau} - \abs{\xi}/\varepsilon }^{\theta}
  \, \widehat u_{1,2}(\tau,\xi)}{_{\tau,\xi}}
$$
we conclude that $\AnormArg{u_{1,2}}{1}{\theta}{S_{T}}$ is 
dominated by
$$
  \varepsilon \bigtwonorm{ \angles{ \tau }^{\theta}
  \, \widehat G_{2}(\tau,\xi)}{_{\tau,\xi}}
  \lesssim \varepsilon \bigtwonorm{\angles{ \tau + 1/\varepsilon^{2}
  }^{\theta} \, \widehat F(\tau,\xi)}{_{\tau,\xi}}
  \lesssim \varepsilon \twonorm{\angles{\partial_{t}}^{\theta} F}{}
  + \frac{\twonorm{F}{}}{\varepsilon^{2\theta-1}}.
$$
This ends the proof of the lemma.
\end{proof}

\section{Main estimates}\label{MainEstimates}

Here we prove the main \emph{a priori} estimates for the nonlinear 
terms in the modified \emph{DM} system, in terms of the following
spacetime norms.

\begin{definition}\label{NormDef} For $T > 0$ we define
\begin{itemize}
  \item
  $X_{T}\eps = \varepsilon^{\fixedpar} \norm{\A\eps}_{\dot 
  \scrH^{1,\theta}\subeps(S_{T})}$,
  \item
  $Y_{T}\eps = \sum_{\pm}
  \norm{\phi_{\pm}\eps}_{X^{1,\theta}_{\tau
  = \pm h\subeps(\xi)}(S_{T})}$,
  \item
  $Z_{T}\eps = \sum_{\pm}
  \norm{(\phi_{\pm}\eps)\low}_{\mixed{2}{6} \cap 
  \mixed{\infty}{2}(S_{T})}$,
\end{itemize}
for $\theta > 1/2$ sufficiently close to $1/2$, independently of 
$\varepsilon$, but depending on the fixed parameter $\fixedpar$.
In fact, the relevant condition is
\begin{equation}\label{ThetaCondition}
  \Lambda + 1 - 2\theta > 0,
\end{equation}
which we assume from now on.
\end{definition}

We also need the following initial data norms.

\begin{definition}\label{DataNormDef}
For initial data \eqref{Data} we set
\begin{itemize}
  \item $X_{0}\eps = \varepsilon^{\fixedpar}
  \left( \Sobdotnorm{\mathbf a_{0}\eps}{1}
  + \varepsilon \twonorm{\mathbf a_{1}\eps}{} \right)$,
  \item $Y_{0}\eps = \Sobnorm{\psi_{0}\eps}{1}$,
  \item $Z_{0}\eps = \twonorm{\psi_{0}\eps}{}$.
\end{itemize}
\end{definition}

In order to simplify the notation we drop the superscript $\varepsilon$
on the fields $\phi, \psi, A_{\mu}$ etc.\ in the remainder of this 
section. We assume $0 \le T \le 1$ in the estimates that follow,
and we write $$P_{T}\eps = P(X_{T}\eps + Y_{T}\eps)$$ where $P(x) = x + x^{N}$ 
for a sufficiently large $N \in \N$, independent of $\varepsilon$.

\subsection{Estimates for $A_{0}$}\label{A0Estimates}

Split $\psi = \psi\low + \psi\high$ and write
\begin{equation}\label{A0Split}
  A_{0} = A_{0}' + A_{0}''
\end{equation}
where $A_{0}'$ corresponds to ``low-low'' interactions:
$$
  \Delta A_{0}' = \innerprod{\psi\low}{\psi\low}.
$$
Then
\begin{alignat}{2}
  \label{A0EstA1}
  &\mixednormlocal{\Delta A_{0}'}{p}{(3/2)^{+}}
  \lesssim (Z_{T}\eps)^{2}
  &\quad
  \text{for} \quad & 1 \le p < 2,
  \\
  \label{A0EstA2}
  &\mixednormlocal{\Delta A_{0}''}{p}{(3/2)^{+}}
  \lesssim \varepsilon^{1^{-}} (Y_{T}\eps)^{2}
  &\quad
  \text{for} \quad & 1\le p < 2,
  \\
  \label{A0EstB1}
  &\mixednormlocal{\Delta A_{0}'}{2}{r}
  \lesssim (Z_{T}\eps)^{2}
  &\quad
  \text{for} \quad & 1 \le r \le \frac{3}{2},
  \\
  \label{A0EstB2}
  &\mixednormlocal{\Delta A_{0}''}{2}{r}
  \lesssim \varepsilon (Y_{T}\eps)^{2}
  &\quad
  \text{for} \quad & 1 \le r \le \frac{3}{2},
  \\
  \label{A0EstC0}
  &\mixednormlocal{\Delta A_{0}}{2}{r}
  \lesssim (Y_{T}\eps)^{2}
  &\quad
  \text{for} \quad & 2 \le r \le 6,
  \\
  \label{A0EstC}
  &\mixednormlocal{\Delta A_{0}}{\infty}{r}
  \lesssim (Y_{T}\eps)^{2}
  &\quad
  \text{for} \quad & 1\le r \le 3.
\end{alignat}
Here \eqref{A0EstC} follows from H\"older's inequality and 
Sobolev embedding, while \eqref{A0EstC0} reduces to
$$
  \mixednormlocal{\psi}{4}{2r} \lesssim Y_{T}\eps
  \quad \text{for} \quad 2 \le r \le 6.
$$
By Sobolev embedding and the Transfer Principle, the latter reduces
to the $\mixed{4}{4}$ Strichartz estimate \eqref{LinearStr}.
Let us now prove \eqref{A0EstA1} and \eqref{A0EstA2}; the proofs
of \eqref{A0EstB1} and \eqref{A0EstB2} are similar. Write
\begin{equation}\label{A0EstLowLow}
  \Lxpnorm{\innerprod{\psi}{\psi}}{(3/2)^{+}}
  \le
  \Lxpnorm{\psi}{6} \Lxpnorm{\psi}{2^{+}}
  \lesssim
  \Lxpnorm{\psi}{6}^{1^{+}} \Lxpnorm{\psi}{2}^{1^{-}}.
\end{equation}
For $\psi = \psi\low$ the $L_{t}^{2^{-}}$ norm of this is 
clearly dominated by r.h.s.\eqref{A0EstA1}.
On the other hand, if at least one $\psi\high$ is present, then
we dominate by r.h.s.\eqref{A0EstA2}
using the $\dot H^{1} \hookrightarrow L_{x}^{6}$ Sobolev embedding and
the estimate \eqref{HighFreqEst}.

We will also need the embeddings
\begin{equation}\label{EllipticEmbeddings}
\begin{split}
  \Lxpnorm{f}{\infty}
  &\lesssim \Lxpnorm{\Delta f}{(3/2)^{-}} + \Lxpnorm{\Delta 
  f}{(3/2)^{+}},
  \\
  \Lxpnorm{\nabla f}{\infty}
  &\lesssim \Lxpnorm{\Delta f}{3^{-}} + \Lxpnorm{\Delta 
  f}{3^{+}}.
\end{split}
\end{equation}

\subsection{Estimates for the remainder term}\label{RemainderEst}

For the remainder term $R\eps$ given by \eqref{RDef} we shall prove
(cf.\ Lemma \ref{XCutoffLemma})
\begin{align}
  \label{REstA}
  \bigXnormArg{e^{\pm it/\varepsilon^{2}} R\eps}{1}{\theta-1}{S_{T}}
  &\lesssim
  \varepsilon^{(1/2-\fixedpar)^{-}} P_{T}\eps,
  \\
  \label{REstC}
  \twonorm{R\eps}{(S_{T})}
  &\lesssim
  \varepsilon P_{T}\eps.
\end{align}
Using \eqref{RDef}, \eqref{LambdaEst}, \eqref{A0Split}
and \eqref{InterpolationEmbeddings}
we dominate l.h.s.\eqref{REstA} by a sum of terms
\begin{align*}
  N_{1} &=\bigXnormArg{e^{\pm it/\varepsilon^{2}}
  \left\{ i(E_{j}-\partial_{j} A_{0}) \alpha^{j} \psi
  - B_{j} S^{j} \psi \right\}}{0}{\theta-1}{S_{T}},
  \\
  N_{2} &= \twonorm{\A \cdot \nabla \psi}{(S_{T})},
  \\
  N_{3} &= \bigLpHslocal{\varepsilon (\partial_{j} A_{0})
  \alpha^{j} \psi}{2}{1},
  \\
  N_{4} &= \frac{1}{\varepsilon}
  \bigtwonorm{ [A_{0}',\lambda\eps-1]\psi_{\pm}}{(S_{T})},
  \\
  N_{5} &= \frac{1}{\varepsilon}
  \bigmixednormlocal{[A_{0}'',\lambda\eps-1]\psi_{\pm}}{2^{-}}{2},
  \\
  N_{6} &=\twonorm{ \varepsilon \, (\A)^{2} \psi }{(S_{T})}.
\end{align*}
All these terms appear also in the \emph{KGM} case (see \cite{BMS}),
with the notable exception of $N_{1}$.
The latter is however the most interesting (and difficult) term,
so we consider it first. Write
$\left\{ i(E_{j}-\partial_{j} A_{0}) \alpha^{j}
- B_{j} S^{j} \right\} \psi
= \sum I_{\mu}$
where
$$
  I_{\mu} = \left\{ i(E_{j}-\partial_{j} A_{0}) \alpha^{j}
  - B_{j} S^{j} \right\} \LPproj_{\mu} \psi
$$
and the sum is over all dyadic numbers $\mu$
of the form $2^{j}$, $j \in \Z$.
Here $\LPproj_{\mu}$ is the Littlewood-Paley operator
defined in Sect.\ \ref{BilinearSection}.
We split into the cases
\begin{enumerate}
  \item\label{LowCase}
  $\mu \le 1/\varepsilon$,
  \item\label{HighCase}
  $\mu > 1/\varepsilon$.
\end{enumerate}

\medskip
\noindent
\emph{Case (\ref{LowCase}).} By \eqref{InterpolationEmbeddings},
we can reduce to proving
\begin{equation}\label{LowCaseEst}
  \twonorm{\sum\nolimits_{\mu \le 1/\varepsilon} 
  I_{\mu}}{(S_{T})}\
  \lesssim \varepsilon^{(1/2-\fixedpar)^{-}}
   X_{T}\eps Y_{T}\eps,
\end{equation}
but this follows from Corollary \ref{DyadicCor2} 
via the Transfer Principle.

\medskip
\noindent
\emph{Case (\ref{HighCase}).}
Using \eqref{InterpolationEmbeddings} we write
\begin{equation}\label{ImuEst}
  \bigXnormArg{e^{\pm it/\varepsilon^{2}}
  I_{\mu}}{0}{\theta-1}{S_{T}} 
  \lesssim \twonorm{I_{\mu}}{(S_{T})}^{1-\sigma}
  \bigXnormArg{e^{\pm it/\varepsilon^{2}}
  I_{\mu}}{0}{\theta-1}{S_{T}}^{\sigma}
\end{equation}
where $0 < \sigma \ll 1$ will be chosen later.
Proceeding as in case \eqref{LowCase}, but using the sharp estimate 
\eqref{BilinearEndpoint}, we obtain
\begin{equation}\label{ImuEstA}
  \twonorm{I_{\mu}}{(S_{T})}
  \lesssim \varepsilon^{1/2-\fixedpar}
  X_{T}\eps Y_{T}\eps.
\end{equation}
We claim there exist $\zeta > 0$ and $M \in \N$,
both independent of $\varepsilon$ and $\mu$, such that
\begin{equation}\label{ImuEstB}
  \bigXnormArg{e^{\pm it/\varepsilon^{2}}
  I_{\mu}}{0}{\theta-1}{S_{T}}
  \lesssim \varepsilon^{-M} \mu^{-\zeta} P_{T}\eps.
\end{equation}
Granting this for the moment, we see that by choosing
$\sigma$ sufficiently small in \eqref{ImuEst}, depending on $M$,
we get
$$
  \sum_{\mu > 1/\varepsilon}
  \bigXnormArg{e^{\pm it/\varepsilon^{2}}
  I_{\mu}}{0}{\theta-1}{S_{T}}
  \lesssim \varepsilon^{(1/2-\fixedpar)^{-}} P_{T}\eps
$$
as desired. Let us prove the claim.
On account of Lemma \ref{NullLemma},
\begin{align*}
  I_{\mu} &=
  Q_{jk}(\D^{-1} \varepsilon \partial_{t} a^{jk}, \LPproj_{\mu} U)
  - Q_{jk}(\D^{-1} \partial_{l} a^{jk}, \alpha^{l} \LPproj_{\mu} U)
  \\
  &\quad
  + Q_{0}(A_{j},\alpha^{j} \LPproj_{\mu} U)
  + Q_{0j}(A_{k}, \alpha^{j}\alpha^{k} \LPproj_{\mu} U)
  - \frac{i}{2} Q_{jk}(A_{m}, \epsilon^{jkl} S_{l}\alpha^{m} \LPproj_{\mu} U)
\end{align*}
where $a^{jk}, U$ are as in Lemma \ref{NullLemma}.
But since $\psi$ solves the Dirac equation,
\begin{equation}\label{BoxU}
  \square\subeps U = - \varepsilon^{-1} \gamma^{0} \psi
  +\varepsilon A_{j} \alpha^{j} \psi + \varepsilon A_{0} \psi
\end{equation}
Now we appeal to the following null form estimate.

\begin{theorem}\label{NullFormThm}
Let $1/2 < \theta < 1$. Then
$$
  \WaveSobnorm{Q(u,v)}{0}{\theta-1}
  \le C_{\theta} \Anorm{u}{1}{\theta}
  \SecondWaveSobnorm{v}{(1+\theta)^{+}}{1}
$$
holds on $\R^{1+3}$ for all null forms $Q$ in (\ref{NullForms}).
Moreover, if $Q = Q_{ij}$, then the norm
$\Anorm{u}{1}{\theta}$ in the r.h.s.\ can be replaced by
$\WaveSobdotnorm{u}{1}{\theta}$.
\end{theorem}

By a standard procedure we reduce this to well-known bilinear estimates 
for the homogeneous wave equation; the proof can be found in Sect.\ 
\ref{NullProof}.

Applying this estimate, and recalling \eqref{XtoWaveSobEst},
we reduce \eqref{ImuEstB} to proving
\begin{equation}\label{ImuEstC}
  \SecondWaveSobnormArg{\LPproj_{\mu} U}{2-2\zeta}{1}{S_{T}}
  \le C_{\zeta} \varepsilon^{-M} \mu^{-\zeta} P_{T}\eps
\end{equation}
for $\zeta > 0$ such that $1 + \theta < 2 - 2\zeta$.
Clearly, it suffices to show
$$
  \SecondWaveSobnormArg{U}{2-\zeta}{1}{S_{T}}
  \le C_{\zeta} \varepsilon^{-M} P_{T}\eps,
$$
but using \eqref{CutoffB} and \eqref{BoxU} we reduce this to
\begin{align}
  \label{NonsharpDiracEst}
  \LpHslocal{A_{j} \alpha^{j} \psi}{2}{1-\zeta}
  &\le C_{\zeta} \varepsilon^{1/2-\fixedpar} X_{T}\eps Y_{T}\eps,
  \\
  \label{A0ProductEst}
  \LpHslocal{A_{0} \psi}{2}{1}
  &\lesssim (Y_{T})^{3}.
\end{align}
The former follows from
Corollary \ref{DyadicCor1} and the Transfer Principle,
while the latter reduces to \eqref{A0EstC} using
Leibniz' rule, H\"older's inequality and \eqref{EllipticEmbeddings}.
This concludes the estimate for $N_{1}$.

\medskip\medskip\medskip
\noindent
It remains to estimate the terms $N_{2},\dots,N_{6}$.
Use Lemma \ref{NullLemma} and \eqref{NullFormEstA} via the
Transfer Principle to see that
$$
  N_{2} \lesssim \varepsilon^{1/2-\fixedpar} X_{T}\eps Y_{T}\eps.
$$
Next, by Leibniz' rule, H\"older's inequality, 
\eqref{EllipticEmbeddings} and \eqref{A0EstC0},
$$
  N_{3} \lesssim \varepsilon (Y_{T}\eps)^{3}.
$$
To the term $N_{4}$ we apply we apply the commutator estimate
$$
  \Lxpnorm{ [ \Delta^{-1}(fg), \lambda\eps-1] h}{2}
  \le \varepsilon^{2} \, C_{\rho}
  \Sobnorm{f}{1+\rho} \Sobnorm{g}{1+\rho} \Sobnorm{h}{1}
  \quad (\text{for all $\rho > 0$})
$$
proved in \cite[Lemma 9]{BMS}. Thus
$$
  N_{4} \lesssim \varepsilon^{1^{-}} (Y_{T}\eps)^{3}.
$$
In $N_{5}$ we simply expand the commutator and apply the estimate
\begin{equation}\label{MBound}
  \twonorm{(\lambda\eps-1)f}{} \lesssim \varepsilon
  \Sobdotnorm{f}{1},
\end{equation}
which follows from \eqref{heps}. Thus
$$
  N_{5} \lesssim \left( 
  \mixednormlocal{A_{0}''}{2^{-}}{\infty}
  + \mixednormlocal{\nabla A_{0}''}{2^{-}}{3} \right) 
  \LpHslocal{\psi_{\pm}}{\infty}{1},
$$
so in view of \eqref{EllipticEmbeddings}, \eqref{A0EstA2} and \eqref{A0EstB2},
$N_{5}$ satisfies the same bound as $N_{4}$.
Finally, by H\"older's inequality and the $\dot H^{1} \hookrightarrow 
L_{x}^{6}$ Sobolev embedding,
$$
  N_{6} \lesssim \varepsilon^{1-2\fixedpar} (X_{T}\eps)^{2} Y_{T}\eps.
$$
This concludes the proof of \eqref{REstA}.

\medskip\medskip\medskip
\noindent
Now consider the estimate \eqref{REstC}. Using \eqref{LambdaEst} with 
$r = 0^{+}$ we get
$$
  \twonorm{R\eps}{(S_{T})}
  \lesssim \varepsilon^{1^{-}}
  \left( \widetilde N_{1} + N_{2} + N_{4} + \widetilde N_{5} + N_{6} 
  \right) + N_{3},
$$
where $N_{2}, N_{3}, N_{4}$ and $N_{6}$ are as before, whereas
\begin{align*}
  \widetilde N_{1} &=
  \LpHslocal{i(E_{j}-\partial_{j} A_{0}) \alpha^{j} \psi
  - B_{j} S^{j} \psi}{2}{0^{-}},
  \\
  \widetilde N_{5} &= \frac{1}{\varepsilon}
  \bigLpHslocal{[A_{0}'',\lambda\eps-1]\psi_{\pm}}{2}{0^{-}}.
\end{align*}
Write $\widetilde N_{1} \le \widetilde N_{1,1} + \widetilde N_{1,2}$
corresponding to
$\psi = \psi_{\low} + \psi\high$. For the low 
frequency case we apply the nonsharp bilinear Strichartz estimate 
in Corollary \ref{DyadicCor2} via the Transfer Principle, to get
$$
  \widetilde N_{1,1} \lesssim \varepsilon^{(1/2-\fixedpar)^{-}} X_{T}\eps 
  Y_{T}\eps.
$$
By Sobolev embedding and H\"older's inequality,
$$
  \widetilde N_{1,2} \lesssim \left( \mixednorm{\nabla \A}{\infty}{2}
  + \varepsilon \mixednorm{\partial_{t} \A}{\infty}{2} \right)
  \mixednorm{\psi\high}{2^{+}}{\infty^{-}},
$$
where the pair $(2^{+},\infty^{-})$ is chosen to be sharp wave 
admissible. Applying the Strichartz estimate \eqref{WaveEmb}
we then obtain the same estimate for $\widetilde N_{1,2}$ as for
$\widetilde N_{1,1}$. Next,
$$
  \widetilde N_{5}
  \lesssim
  \left( \mixednormlocal{A_{0}''}{2}{\infty^{-}}
  + \mixednormlocal{\nabla A_{0}''}{2}{3} \right) 
  \LpHslocal{\psi_{\pm}}{\infty}{1}
  \lesssim \varepsilon (Y_{T}\eps)^{3},
$$
where we used \eqref{A0EstB2} to get the last inequality.
This ends the proof of \eqref{REstC}.

\subsection{Estimates for the current density}\label{CurrentEstimates}

Split $\psi = \psi\low + \psi\high$ and write
$$
  \J = \J' + \J''
$$
where $\J'$ corresponds to ``low-low'' interactions:
$$
  \varepsilon \J' = \left\{ \innerprod{\alpha^{k}\psi\low}{\psi\low} \right\}_{k=1,2,3}.
$$
By \eqref{WaveEmb} and \eqref{SchrEmb2},
\begin{gather*}
  \twonorm{u\high v\low}{} \lesssim \mixednorm{u\high}{\infty}{2}
  \mixednorm{v\low}{2}{\infty}
  \lesssim \varepsilon \Xnorm{u\high}{1}{\theta}
  \Xnorm{v\low}{1}{\theta},
  \\
  \twonorm{u\high v\high}{} \lesssim \Lpnorm{u\high}{4}
  \Lpnorm{v\high}{4} \lesssim \varepsilon^{3/2} \Xnorm{u\high}{1}{\theta}
  \Xnorm{v\high}{1}{\theta},
\end{gather*}
whence
\begin{equation}\label{JhighEst}
  \twonorm{\J''}{(S_{T})}
  \lesssim (Y_{T}\eps)^{2}.
\end{equation}
In order to estimate $\J'$ we expand it as in
\eqref{CurrentExpansion}. Thus, we write
$$
  \J' = (\J')_{1} + (\J')_{2}
$$
where
$$
  \varepsilon (\J')_{1} = 
  2 \re \left\{ e^{-2it/\varepsilon^{2}} \innerprod{\sigma^{j}
  (\chi_{+})\low}{(\eta_{-})\low} 
  \right\}_{j = 1,2,3}
$$
whereas $(\J')_{2}$ consists of products containing at least one of the fields
$(\chi_{-})\low$ or $(\eta_{+})\low$, which we expect to be small.
The latter we estimate, just to take one of these terms,
$$
  \twonorm{ \innerprod{\sigma^{j} (\chi_{+})\low}
  {(\eta_{+})\low} }{(S_{T})}
  \le \mixednormlocal{(\chi_{+})\low}{2}{\infty}
  \mixednormlocal{(\eta_{+})\low}{\infty}{2}.
$$
Thus, using \eqref{SchrEmb2} and \eqref{SmallComponentsBound},
\begin{equation}\label{JlowlowSmall}
  \twonorm{ (\J')_{2} }{(S_{T})}
  \lesssim
  (Y_{T}\eps)^{2} + \varepsilon^{1-\fixedpar} P_{T}\eps.
\end{equation}
To that part of $\A$ which corresponds to $(\J')_{1}$ we are going to apply
Lemma \ref{IntByPartsCutoffLemma}. Hence we want to estimate
\begin{equation}\label{MDef}
  M_{T}\eps := \varepsilon \twonorm{\partial_{t} \mathbf F}{(S_{T})}
  + \LpHs{\mathbf F\ext}{2}{1}
  + \varepsilon \twonorm{\angles{\partial_{t}}^{\theta} \mathbf F\ext}{_{t,x}}
  + \frac{\twonorm{\mathbf F\ext}{_{t,x}}}{\varepsilon^{2\theta-1}}
\end{equation}
where $\mathbf F\ext$ is an extension of
\begin{equation}\label{FDef}
  \mathbf F = \left\{ \innerprod{\sigma^{j} 
  (\chi_{+})\low}{(\eta_{-})\low} \right\}_{j = 1,2,3}
\end{equation}
from $S_{T}$ to all of $\R^{1+3}$. To choose this extension, let
$$
  \phi'_{\pm}
  =
  \begin{pmatrix}
    \chi'_{\pm}  \\
    \eta'_{\pm}
  \end{pmatrix}
  \in X^{1,\theta}_{\tau = \pm h\subeps(\xi)}
$$
be arbitrary extensions of $\phi_{\pm}$ and define $\mathbf F\ext$ by
\eqref{FDef} with $\chi_{+}$ and $\eta_{-}$ replaced by their 
respective extensions. From now on we denote $\mathbf F\ext$ simply 
by $\mathbf F$. We claim that
\begin{equation}\label{MEst}
  M_{T}\eps \lesssim
  (Y_{T})^{2}
  + \varepsilon^{(1/2)^{-}} P_{T}\eps
  + \varepsilon^{1-2\theta}
  \bigXnormPlus{\phi'_{+}}{1}{\theta}
  \bigXnormMinus{\phi'_{-}}{1}{\theta}
\end{equation}
If this holds, then taking the \emph{inf} over all extensions 
yields
\begin{equation}\label{MEst2}
  M_{T}\eps \lesssim \varepsilon^{1-2\theta} P_{T}\eps.
\end{equation}
Let us prove \eqref{MEst}. First,
\begin{equation}\label{FEstA}
\begin{split}
  \varepsilon \twonorm{\partial_{t} \mathbf F}{(S_{T})}
  &\le \frac{1}{\varepsilon} \mixednormlocal{(\lambda\eps - 1) 
  \chi_{+}}{\infty}{2} \mixednormlocal{(\eta_{-})\low}{2}{\infty}
  \\
  &\quad
  + \varepsilon \twonorm{L_{+}\eps \chi_{+}}{(S_{T})}
  \inftynorm{(\eta_{-})\low}{(S_{T})}
  + (\dots)
  \\
  &\lesssim
  (Y_{T})^{2} + \varepsilon^{(1/2)^{-}} (Y_{T})^{4}
  + \varepsilon^{(3/2)^{-}} P_{T}\eps,
\end{split}
\end{equation}
where $(\dots)$ stands for symmetric terms.
Here we used \eqref{MBound}, \eqref{SchrEmb2} and
\begin{align*}
  \twonorm{L_{+}\eps \chi_{+}}{(S_{T})}
  &\lesssim (Y_{T}\eps)^{3} + \varepsilon P_{T}\eps,
  \\
  \inftynorm{(\eta_{-})\low}{(S_{T})}
  &\lesssim \varepsilon^{-(1/2)^{+}} Y_{T}\eps.
\end{align*}
The former was obtained from
\eqref{ModifiedDirac}, \eqref{REstC} and \eqref{A0EstC},
while the latter follows from \eqref{BasicEmbedding} and Sobolev embedding.
Second, we write
\begin{equation}\label{FEstB}
\begin{split}
  \LpHs{\mathbf F}{2}{1}
  &\lesssim \LpHs{\chi'_{+}}{\infty}{1} \mixednorm{(\eta'_{-})\low}{2}{\infty}
  + \mixednorm{(\chi'_{+})\low}{2}{\infty}
  \LpHs{\eta'_{-}}{\infty}{1}
  \\
  &\lesssim \bigXnormPlus{\phi'_{+}}{1}{\theta}
  \bigXnormMinus{\phi'_{-}}{1}{\theta},
\end{split}
\end{equation}
where we used \eqref{SchrEmb2} and \eqref{BasicEmbedding}.
This estimate can of course
also be used for the term $\twonorm{\mathbf F}{_{t,x}}$.
Third,
\begin{equation}\label{FEstC}
  \varepsilon \twonorm{\angles{\partial_{t}}^{\theta} \mathbf F\ext}{_{t,x}}
  \lesssim \varepsilon^{1-\theta} \bigXnormPlus{\phi'_{+}}{1}{\theta}
  \bigXnormMinus{\phi'_{-}}{1}{\theta},
\end{equation}
where we used the following:

\begin{lemma}
$
  \twonorm{\angles{\partial_{t}}^{\theta} (u\low v\low)}{_{t,x}}
  \lesssim \varepsilon^{-\theta}
  \XnormPlus{u\low}{1}{\theta} \XnormMinus{v\low}{1}{\theta}.
$
\end{lemma}

\begin{proof} To simplify the notation, let us write $u,v$
instead of $u\low, v\low$ here.
W.l.o.g.\ we assume $\widehat u(\tau,\xi), \widehat 
v(\tau,\xi) \ge 0$. Then using Plancherel's theorem,
$$
  \twonorm{\angles{\partial_{t}}^{\theta} (uv)}{}
  \lesssim \twonorm{u T_{-}^{\theta} v}{}
  + \varepsilon^{-\theta} \bigtwonorm{u \D^{\theta} v}{}
  + \twonorm{v T_{+}^{\theta} u}{}
  + \varepsilon^{-\theta} \bigtwonorm{v \D^{\theta} u}{}.
$$
Here $T_{\pm}^{\theta}$ is the multiplier with symbol
$\angles{ \tau \mp h\subeps(\xi) }^{\theta}$, and we used \eqref{heps}.
Write
$$
  \twonorm{u T_{-}^{\theta} v}{}
  \le \mixednorm{u}{\infty}{3} \mixednorm{T_{-}^{\theta} v}{2}{6}
  \lesssim \XnormPlus{u}{1}{\theta} \XnormMinus{v}{1}{\theta}
$$
using \eqref{BasicEmbedding} and Sobolev embedding. Next,
$$
  \bigtwonorm{u \D^{\theta} v}{}
  \le \mixednorm{u}{\infty}{3} \bigmixednorm{\D^{\theta} v}{2}{6}
  \lesssim \XnormPlus{u}{1}{\theta} \XnormMinus{v}{1}{\theta}
$$
where \eqref{SchrEmb} with $(q,r) = (2,6)$ was used.
\end{proof}

Finally, combining \eqref{FEstA}, \eqref{FEstB} and \eqref{FEstC}, we get 
\eqref{MEst}.

\section{Iteration scheme and local existence}\label{LocalExistence}

For fixed $\varepsilon$ we shall prove the following local existence theorem:

\begin{theorem}\label{SeLWPThm} For fixed 
$\varepsilon$, the Dirac-Maxwell-Coulomb system 
\eqref{DMCscaled} is locally well posed for initial data in the space
\eqref{Data}. The existence time $T > 0$ only depends on 
$\varepsilon$ and the size of the norms of the data, and the solution
is in the space
\begin{equation}\label{LWPspace}
  \psi\eps \in H^{1,\theta}\subeps(S_{T}),
  \qquad
  \A\eps \in \dot \scrH^{1,\theta}\subeps(S_{T}),
  \qquad
  A_{0}\eps \in C([0,T];\dot H^{1}),
\end{equation}
for all $1/2 < \theta < 1$. Moreover,
the solution is unique in this regularity class, and we have
$$
  \phi_{\pm}\eps \in X^{1,\theta}_{\tau = \pm h\subeps(\xi)}(S_{T}),
$$
where $\phi_{\pm}\eps$ is defined by \eqref{KGsplitting} and 
\eqref{ModulatedFields}.
\end{theorem}

We shall prove this by Picard iteration.
In order to simplify the notation we drop the superscript $\varepsilon$ on the fields 
$\psi, A_{\mu}$ etc.\ and introduce instead a superscript $(m)$ to denote the
$m$-th iterate of a field. For \eqref{DMCscaled} we use the iteration scheme
\begin{subequations}\label{DMCiterates}
\begin{align}
  \label{DiracIterate}
  \left\{ i \varepsilon \partial_{t}
  + i \alpha^{j} \partial_{j}
  - (1/\varepsilon) \gamma^{0} \right\} \psi\iterate{m+1}
  &= - \varepsilon A_{j}\iterate{m} \alpha^{j} \psi\iterate{m}
  - \varepsilon A_{0}\iterate{m} \psi\iterate{m},
  \\
  \label{A0Iterate}
  \Delta A_{0}\iterate{m} &= \rho\iterate{m},
  \\
  \label{AIterate}
  \square\subeps \A\iterate{m+1} &= \varepsilon \Proj \J\iterate{m},
\end{align}
\end{subequations}
with initial data as in \eqref{Data},
where $\rho\iterate{m}$ and $\J\iterate{m}$ are given by 
\eqref{ChargeCurrent} with $\psi$ replaced by its iterate $\psi\iterate{m}$.
Note that $A_{0}$ is not really iterated; \eqref{A0Iterate} simply 
defines $A_{0}\iterate{m}$ in terms of $\psi\iterate{m}$.
Observe also that for all $m$,
$$
  \dv \A\iterate{m} = 0,
$$
since $w =  \A\iterate{m} - \Proj \A\iterate{m}$ satisfies $\square 
w = 0$ with vanishing initial data.

By convention we start the iteration at $m = -1$ and set all
iterates identically equal to zero there. Then the iterates 
$\psi\iterate{0}, \A\iterate{0}$ are just solutions of the free Dirac
and wave equations with data \eqref{Data}.
Define (cf.\ \eqref{KGsplitting} and \eqref{ModulatedFields})
\begin{subequations}\label{IterationScheme}
\begin{align}
  \label{IteratesKGsplitting}
  \psi_{\pm}\iterate{m+1} &=
  \frac{1}{2} \left\{ \psi\iterate{m+1} \pm 
  \varepsilon^{2} [\lambda\eps]^{-1} \left( i \partial_{t} 
  \psi\iterate{m+1} + A_{0}\iterate{m} \psi\iterate{m} \right) \right\},
  \\
  \label{IteratesModulatedFields}
  \phi_{\pm}\iterate{m}
  &=
  \begin{pmatrix}
    \chi_{\pm}\iterate{m}  \\
    \eta_{\pm}\iterate{m}
  \end{pmatrix}
  := e^{\pm it/\varepsilon^{2}} \psi_{\pm}\iterate{m}.
\end{align}
\end{subequations}
Proceeding as in the proof of Lemma \ref{ModifiedDiracLemma} one finds
\begin{equation}\label{ModifiedIterateDirac}
  L_{\pm}\eps \phi_{\pm}\iterate{m+1}
  = - A_{0}\iterate{m} \phi_{\pm}\iterate{m}
  \pm \frac{1}{2} e^{\pm it/\varepsilon^{2}} R\iterate{m},
\end{equation}
where
\begin{align}
  \notag
  \lambda\eps R\iterate{m} &=
  \varepsilon \mathcal B\iterate{m}
  + \varepsilon^{2} \mathcal C\iterate{m}
  - \left[A_{0}\iterate{m},\lambda\eps\right]
  \left(\psi_{+}\iterate{m} - \psi_{-}\iterate{m}\right),
  \\
  \label{RIterateDef}
  \mathcal B\iterate{m} &=
  \left\{ 2i \A\iterate{m} \cdot \nabla
  + i E_{j}\iterate{m} \alpha^{j}
  - B_{j}\iterate{m} S^{j} \right\} \psi\iterate{m},
  \\
  \notag
  \mathcal C\iterate{m}
  &= \left\{
  A_{j}\iterate{m} A_{k}\iterate{m-1} \alpha^{j} \alpha^{k}
  + A_{j}\iterate{m} A_{0}\iterate{m-1} \alpha^{j}
  - A_{0}\iterate{m} A_{j}\iterate{m-1} \alpha^{j}
  \right\} \psi\iterate{m-1},
\end{align}
and $E_{j}\iterate{m}, B_{j}\iterate{m}$ are given by \eqref{EandB} 
with $A_{\mu}$ replaced by $A_{\mu}\iterate{m}$.

\medskip\medskip\medskip
\noindent
We now turn to the proof of Theorem \ref{SeLWPThm}. By standard 
arguments, this reduces to proving closed estimates for the iterates 
in the space \eqref{LWPspace}. Set
$$
  B\iterate{m}_{T} = \bigWaveSobnormArg{\psi\iterate{m}}{1}{\theta}{S_{T}}
  + \bigAnormArg{\A\iterate{m}}{1}{\theta}{S_T},
$$
and denote by $B_{0}$ the norm of the data \eqref{Data}.
Then it suffices to prove
\begin{equation}
  B\iterate{m+1}_{T} \le CP(B_{0}) + C T^{\delta} P \left( B\iterate{m}_{T} + 
  B\iterate{m-1}_{T} \right),
\end{equation}
for some constants $C, \delta > 0$ and a polynomial $P$ with $P(0) = 
0$. Here $C$ and $P$ may depend on $\varepsilon$, but since the 
latter is fixed we do not indicate this explicitly. In what follows, 
$C, \delta$ and $P$ may change from line to line.
(Observe also that since all the nonlinear terms in 
\emph{DM} are in fact multilinear, the same arguments then give estimates for 
a difference of two iterates.) By Lemma \ref{HomWaveCutoffLemma},
\begin{equation}\label{SimpleAEst}
  \bigAnormArg{\A\iterate{m+1}}{1}{\theta}{S_T}
  \le CB_{0} + C T^{\delta} \bignorm{\psi\iterate{m}}_{L^{4}(S_{T})}^{2},
\end{equation}
where the $T^{\delta}$ comes from H\"older's inequality in time.
Now apply the Strichartz estimate \eqref{LinearStr}
via the Transfer Principle to see that 
$\bignorm{\psi\iterate{m}}_{L^{4}(S_{T})} \lesssim B\iterate{m}_{T}$.
In order to estimate $\psi\iterate{m+1}$ we use the splitting \eqref{PsiExpansion} 
and the embedding
\begin{equation}\label{XintoHEmbeddings}
  \WaveSobnorm{u}{s}{\theta} \lesssim \varepsilon^{-2\theta} 
  \Xnorm{u}{s}{\theta},
\end{equation}
which holds in view of \eqref{hEst}. Thus, we write
$$
  \bigWaveSobnormArg{\psi\iterate{m+1}}{1}{\theta}{S_{T}}
  \le C \sum_{\pm} \bigXnormArg{\phi_{\pm}\iterate{m+1}}{1}{\theta}{S_{T}}.
$$
Using Lemma \ref{XCutoffLemma} and \eqref{InterpolationEmbeddings}, we bound
$\bigXnormArg{\phi_{\pm}\iterate{m+1}}{1}{\theta}{S_{T}}$ by
\begin{equation}\label{SimplePhiEst}
  C P(B_{0}) + CT^{\delta} \bigenergylocal{A_{0}\iterate{m} 
  \phi_{\pm}\iterate{m}}
  + C \bigXnormArg{e^{\pm it/\varepsilon^{2}} 
  R\iterate{m}}{0}{\theta-1}{S_{T}}.
\end{equation}
The second term is trivial to bound, since
\begin{equation}\label{A0SimpleEst}
  \bigmixednormlocal{\Delta A_{0}\iterate{m}}{\infty}{r} \lesssim
  \bigWaveSobnormArg{\psi\iterate{m}}{1}{\theta}{S_{T}}^{2}
  \quad \text{for} \quad 1 \le r \le 3
\end{equation}
and
\begin{equation}\label{IterateKGsplittingBound}
  \bigSobnorm{\psi_{\pm}\iterate{m}(t)}{1}
  \lesssim
  \bigSobnorm{\psi\iterate{m}(t)}{1}
  + \varepsilon \bigSobdotnorm{\A\iterate{m-1}(t)}{1} 
  \bigSobnorm{\psi\iterate{m-1}(t)}{1}.
\end{equation}
The latter is just the analogue of \eqref{KGsplittingBound} for the iterates.

For the third term in \eqref{SimplePhiEst} we can apply 
the estimates proved in Sect.\ \ref{RemainderEst}. In fact, we claim 
that the proof of \eqref{REstA} gives
$$
  \bigXnormArg{e^{\pm it/\varepsilon^{2}} 
  R\iterate{m}}{0}{\theta-1}{S_{T}}
  \le C T^{\delta} P \left( B\iterate{m}_{T} + 
  B\iterate{m-1}_{T} \right).
$$
To see this, consider one by one the terms $N_{1},\dots,N_{6}$ 
in Sect.\ \ref{RemainderEst}. For $N_{1}$ and $N_{2}$ we only have to 
observe that the bilinear estimates in Corollaries \ref{DyadicCor1} and
\ref{DyadicCor2} as well as the null form estimate \eqref{NullFormEstA}
are valid also in the case where both $u$ and $v$ solve the homogeneous
wave equation, so we can apply the Transfer Principle for the 
$H^{s,\theta}\subeps$ spaces instead of $X^{s,\theta}_{\tau = \pm 
h\subeps(\xi)}$.
Note also that \eqref{BoxU} must be replaced by
\begin{equation}\label{BoxUIterate}
  \square\subeps U = - \varepsilon^{-1} \gamma^{0} \psi\iterate{m}
  +\varepsilon A_{j}\iterate{m-1} \alpha^{j} \psi\iterate{m-1}
  + \varepsilon A_{0}\iterate{m-1} \psi\iterate{m-1}.
\end{equation}
The estimate for $N_{6}$ requires no change. Finally, the estimates 
for the terms involving $A_{0}$ can be simplified, since we do not 
care about powers of $\varepsilon$ here. Thus, in $N_{3}$ we can 
replace $H^{1}$ by $L_{x}^{2}$, by giving up the $\varepsilon$, and then 
the estimate reduces to \eqref{A0SimpleEst}. Finally, the commutator 
terms $N_{4}$ and $N_{5}$ are replaced by a single term, since we do not need 
to split $A_{0}$ according to \eqref{A0Split}. We simply expand the 
commutator and proceed as in the estimate for $N_{5}$, reducing to 
\eqref{A0SimpleEst} and \eqref{IterateKGsplittingBound}. This 
concludes the proof of Theorem \ref{SeLWPThm}.

\section{Uniform $H^{1}$ bounds and long time 
existence}\label{UniformBounds}

We shall prove:

\begin{theorem}\label{UniformBoundsTheorem}
Consider the solution $(\psi\eps,A_{\mu}\eps)$ of (\ref{DMCscaled}), 
(\ref{Data}) from Theorem \ref{SeLWPThm}, existing up to a time
$T\subeps > 0$ and belonging to the space (\ref{LWPspace}) over
this time interval. There exist
\begin{enumerate}
  \item a time $T^{*} > 0$ depending only on
  $\sup_{\varepsilon > 0} \twonorm{\psi_{0}\eps}{}$,
  \item constants $C, M, \varepsilon_{0} > 0$
  independent of $\varepsilon$,
\end{enumerate}
such that if
\begin{equation}\label{XYInitialBound}
  X_{0}\eps + Y_{0}\eps \le B
  \quad \text{for all} \quad \varepsilon
\end{equation}
then
\begin{equation}\label{XYBound}
  X_{T}\eps + Y_{T}\eps \le C B
  \quad \text{for} \quad
  \varepsilon < \frac{\varepsilon_{0}}{1 + (CB)^{M}}
  \quad \text{and} \quad 0 \le T \le \min(T^{*},T\subeps).
\end{equation}
Moreover, there is a polynomial $P$ with $P(0) = 0$, independent of 
$\varepsilon$, such that
\begin{equation}\label{ZBound}
  Z_{T}\eps \le C Z_{0}\eps + \varepsilon P(X_{0}\eps + Y_{0}\eps),
\end{equation}
for $T,\varepsilon$ as in \eqref{XYBound}
\end{theorem}

We claim that this result, together with Theorem \ref{SeLWPThm}, 
implies Theorem \ref{LWPThm}. To see this, first observe that
the bound in \eqref{XYBound} implies, on account of 
\eqref{BasicEmbedding},
$$
  \Sobnorm{\psi\eps(T)}{1}
  + \varepsilon^{\fixedpar}
  \left\{ \Sobdotnorm{\A\eps(T)}{1}
  + \varepsilon \twonorm{\partial_{t}\A\eps(T)}{} \right\}
  \le C' B
$$
for some constant $C' > C$ independent of $\varepsilon$.
Thus Theorem \ref{SeLWPThm} implies $T\subeps \ge T^{*}$
for $\varepsilon$ as in \eqref{XYBound}.
In view of the conservation of charge \eqref{ChargeConservation}
for the Dirac equation, we can iterate this argument any 
number of times, obtaining
$$
  T\subeps \ge N T^{*}
  \quad \text{and} \quad
  X_{NT^{*}}\eps + Y_{NT^{*}}\eps \le (C')^{N} B
  \quad \text{for} \quad
  \varepsilon < \frac{\varepsilon_{0}}{1 + (C')^{MN} B^{M}}
$$
for all $N \in \N$. This proves Theorem \ref{LWPThm}.

\medskip\medskip\medskip
\noindent
We shall prove Theorem \ref{UniformBoundsTheorem} using the
iteration scheme from Sect.\ \ref{LocalExistence}. In order to simplify the notation
we drop the superscript $\varepsilon$ on the fields 
$\psi, A_{\mu}$ etc.\ as well as on the $XYZ$-norms in Definitions 
\ref{NormDef} and \ref{DataNormDef} in the remainder of this section,
and introduce instead a superscript $(m)$ to denote the $m$-th iterate of a 
field. We denote by $X_{T}\iterate{m}$ etc.\  the norms in Definition 
\ref{NormDef} with the respective fields replaced by their $m$-th 
iterate.
Then we have:

\begin{proposition}\label{XYZProp}
There exist $C, \gamma, \delta > 0$ and
a polynomial $P$ with $P(0) = 0$, all independent of $\varepsilon$,
such that the estimates
\begin{subequations}\label{XYZEstimates}
\begin{align}
  \label{XEstimate}
  X_{T}\iterate{m+1}
  &\le C X_{0}
  + \varepsilon^{\gamma} P_{T}\iterate{m},
  \\
  \label{YEstimate}
  Y_{T}\iterate{m+1} &\le C Y_{0}
  + C T^{\delta} \left[Z_{T}\iterate{m}\right]^{2} 
  Y_{T}\iterate{m} + \varepsilon^{\gamma} P_{T}\iterate{m},
  \\
  \label{ZEstimate}
  Z_{T}\iterate{m+1} &\le C Z_{0}
  + C T^{\delta} \left[Z_{T}\iterate{m}\right]^{2} 
  Z_{T}\iterate{m} + \varepsilon P_{T}\iterate{m},
\end{align}
\end{subequations}
hold for $T \le 1$ and $m \ge -1$, where
\begin{equation}\label{PTDef}
  P_{T}\iterate{m} :=
  \begin{cases}
  P(X_{0} + Y_{0}) \quad &\text{for} \quad m = -1,
  \\
  P\bigl( X_{T}\iterate{m} + X_{T}\iterate{m-1} + Y_{T}\iterate{m}
  + Y_{T}\iterate{m-1}\bigr) \quad &\text{for} \quad m \ge 0.
  \end{cases}
\end{equation}
In fact, these estimates hold for (recall (\ref{ThetaCondition}))
\begin{equation}\label{GammaCondition}
  \gamma \le \Lambda + 1 - 2\theta.
\end{equation}
\end{proposition}

The proof is deferred to the end of this section.

\begin{corollary}\label{XYZCorollary}
There exist $C,\delta > 0$ and a polynomial $Q$, all 
independent of $\varepsilon$, such 
that if $\gamma > 0$ is sufficiently small depending on $\fixedpar$,
and $T, \varepsilon > 0$ are taken so small that
\begin{equation}\label{Smallness}
  2 C T^{\delta} \bigl[ 2 C \twonorm{\psi_{0}}{} + 1 \bigr]^{2}
  \le 1,
  \qquad 2 \varepsilon^{\gamma/2} Q( X_{0} + Y_{0} ) \le 1,
\end{equation}
then
\begin{subequations}\label{XYZIterateBounds}
\begin{align}
  \label{XIterateBound}
  X_{T}\iterate{m}
  &\le C X_{0} + \varepsilon^{\gamma/2} ( X_{0} + Y_{0} ),
  \\
  \label{YIterateBound}
  Y_{T}\iterate{m} &\le 2 C Y_{0}
  + \varepsilon^{\gamma/2} ( X_{0} + Y_{0} ),
  \\
  \label{ZIterateBound}
  Z_{T}\iterate{m} &\le 2 C Z_{0}
  + \varepsilon^{1-\gamma/2} ( X_{0} + Y_{0} ),
\end{align}
\end{subequations}
for $m \ge 0$.
\end{corollary}

\begin{proof} This is a simple induction.
Since $P(0) = 0$ in Proposition \ref{XYZProp}, there is a polynomial 
$Q(r)$ such that
$$
  P(4[C+1]r) \le r Q(r) \quad \text{for} \quad r \ge 0.
$$
Then
\begin{equation}\label{PBound}
  P_{T}\iterate{m} \le Q( X_{0} + Y_{0} ) \cdot
  ( X_{0} + Y_{0} )
\end{equation}
holds for $m = -1$, in view of the definition \eqref{PTDef}.
Hence \eqref{XYZIterateBounds} for $m = 0$ follows from
\eqref{XEstimate}--\eqref{ZEstimate}
and the fact that the iterates at $m = -1$ all vanish.
Now assume \eqref{XYZIterateBounds} holds for $0 \le m \le m_{0}$.
Then \eqref{PBound} holds for such $m$, and using \eqref{Smallness} and
\eqref{XYZEstimates} we obtain \eqref{XYZIterateBounds} for $m = m_{0} + 1$.
\end{proof}

We are now in a position to prove Theorem \ref{UniformBoundsTheorem}.
Indeed, from the proof of Theorem \ref{SeLWPThm} we know that the 
iterates $\phi_{\pm}\iterate{m}$ converge in the $Y$-norms. We can 
therefore pass to the limit $m \to \infty$ in Corollary \ref{XYZCorollary}.
Thus, from \eqref{XIterateBound}, \eqref{YIterateBound} we get 
\eqref{XYBound}, and from \eqref{ZIterateBound} we get
\begin{equation}\label{FinalZBound}
  Z_{T} \le 2 C Z_{0} + 1.
\end{equation}
Substituting the latter into the second term in the r.h.s.\ of 
\eqref{ZEstimate} in the limit $m \to \infty$,
we then obtain \eqref{ZBound}. This proves Theorem \ref{UniformBoundsTheorem}.

\begin{proof}[Proof of Proposition \ref{XYZProp}.]
By the estimates in Sect.\ \ref{RemainderEst},
\begin{align}
  \label{REstAA}
  \bigXnormArg{e^{\pm it/\varepsilon^{2}} R\iterate{m}}{1}{\theta-1}{S_{T}}
  &\lesssim
  \varepsilon^{(1/2-\fixedpar)^{-}} P_{T}\iterate{m},
  \\
  \label{REstCC}
  \bigtwonorm{R\iterate{m}}{(S_{T})}
  &\lesssim
  \varepsilon P_{T}\iterate{m},
\end{align}
the only difference being that \eqref{BoxU} must be replaced by
\eqref{BoxUIterate}.
Then \eqref{XEstimate} follows from the equation \eqref{AIterate}
by applying Lemmas \ref{HomWaveCutoffLemma} and \ref{IntByPartsCutoffLemma},
the embedding \eqref{InterpolationEmbeddings}, and the estimates proved in
Sect.\ \ref{CurrentEstimates}. However, instead of the estimate
\eqref{SmallComponentsBound}, which was used to prove 
\eqref{JlowlowSmall}, we use the analogous estimate for the iterates:
$$
  \bigLxpnorm{\Pi_{-}^{0} \psi_{+}\iterate{m}}{2} + 
  \bigLxpnorm{\Pi_{+}^{0} \psi_{-}\iterate{m}}{2}
  \lesssim
  \varepsilon \bigSobnorm{\psi\iterate{m}}{1}
  + \varepsilon^{2}
  \bigSobdotnorm{\A\iterate{m-1}}{1} \bigSobnorm{\psi\iterate{m-1}}{1}.
$$
Next, applying Lemma \ref{XCutoffLemma} to the 
equation \eqref{ModifiedIterateDirac} and using the embedding
\eqref{InterpolationEmbeddings} and the estimate \eqref{REstAA},
as well as \eqref{KGsplittingBound} at $t = 0$,
we reduce \eqref{YEstimate} to proving
$$
  \bigLpHslocal{A_{0}\iterate{m}\phi_{\pm}\iterate{m}}{2^{-}}{1}
  \lesssim T^{\delta}
  \bigl[Z_{T}\iterate{m}\bigr]^{2} Y_{T}\iterate{m}
  + \varepsilon^{1^{-}} \bigl[ Y_{T}\iterate{m} \bigr]^{3}.
$$
But this follows from Leibniz' rule, H\"older's inequality,
\eqref{EllipticEmbeddings} and \eqref{A0EstA1}--\eqref{A0EstB2}, in 
view of \eqref{A0Iterate}. The factor
$T^{\delta}$ comes from applying H\"older's inequality in time.
Finally, consider \eqref{ZEstimate}. Apply Proposition \ref{SchrodingerProp}
to \eqref{ModifiedIterateDirac} and use \eqref{KGsplittingBound} at $t = 0$ to get
$$
  Z_{T}\iterate{m+1}
  \lesssim Z_{0} + \varepsilon^{2-\fixedpar} X_{T}\iterate{m} Y_{T}\iterate{m}
  + \sum_{\pm} \bigmixednormlocal{A_{0}\iterate{m} \phi_{\pm}\iterate{m}}{1^{+}}{2^{-}}
  + \bigmixednormlocal{R\iterate{m}}{1}{2}.
$$
The last term is covered by \eqref{REstCC}.
On account of \eqref{A0EstB1} and \eqref{A0EstB2},
$$
  \bigmixednormlocal{A_{0}\iterate{m} \phi_{\pm}\iterate{m}}{1^{+}}{2^{-}}
  \lesssim T^{\delta} \bigl[Z_{T}\iterate{m}\bigr]^{2}
  \bigmixednormlocal{\phi_{\pm}\iterate{m}}{\infty}{2}
  + \varepsilon \bigl[ Y_{T}\iterate{m} \bigr]^{3}.
$$
Then \eqref{ZEstimate} follows, in view of
$$
  \bigmixednormlocal{\phi_{\pm}\iterate{m}}{\infty}{2}
  \lesssim Z_{T}\iterate{m}
  + \varepsilon Y_{T}\iterate{m},
$$
which holds by \eqref{HighFreqEst}.
\end{proof}

\section{Higher order bounds}\label{HigherBoundsSection}

Here we prove bounds for higher order derivatives. For $m = 
0,1,2,\dots$, set (cf.\ Definitions \ref{NormDef} and 
\ref{DataNormDef})
\begin{itemize}
  \item
  $X_{T}\eps[m] = \varepsilon^{\fixedpar}
  \sum_{\abs{\alpha} \le m}
  \norm{\partial_{x}^{\alpha} \A\eps}_{\dot 
  \scrH^{1,\theta}\subeps(S_{T})}$,
  \item
  $Y_{T}\eps[m] = \sum_{\abs{\alpha} \le m} \sum_{\pm}
  \norm{\partial_{x}^{\alpha} \phi_{\pm}\eps}_{X^{1,\theta}_{\tau
  = \pm h\subeps(\xi)}(S_{T})}$,
  \item
  $X_{0}\eps[m] = \varepsilon^{\fixedpar}
  \left(
  \Sobnorm{\nabla \mathbf a_{0}\eps}{m}
  + \varepsilon \Sobnorm{\mathbf a_{1}\eps}{m} \right)$,
  \item
  $Y_{0}\eps[m] = \Sobnorm{\psi_{0}\eps}{m+1}$.
\end{itemize}
The local well-posedness of \emph{DM} in these norms for $m = 0$
was established in Sect.\ \ref{LocalExistence}, and a standard 
argument shows that higher regularity persists, i.e.\ if $ X_{0}\eps[m] + 
Y_{0}\eps[m]$ is finite for some $m \ge 1$, then
$X_{T}\eps[m] + Y_{T}\eps[m]$ is also finite in the interval of existence
$0 \le T \le T\subeps$. Here we concentrate on proving bounds which 
are uniform in $\varepsilon$.
Thus, we shall prove:

\begin{proposition}\label{HigherOrderProp}
If
\begin{equation}\label{HigherOrderDataBound}
  X_{0}\eps[m] + Y_{0}\eps[m] = O(1)
\end{equation}
then
\begin{equation}\label{HigherOrderBound}
  X_{T}\eps[m] + Y_{T}\eps[m] = O(1)
\end{equation}
for $0 \le T \le T\subeps$, where $T\subeps$ is the existence time 
from Theorem \ref{LWPThm}.
\end{proposition}

We claim there exist $C,\delta,\gamma > 0$ and polynomials $Q, 
P_{m}$---all independent of $\varepsilon$---such that for $0 \le T \le 
1$ and $m \ge 1$,
\begin{align}
  \label{XmEst}
  X_{T}\eps[m] &\le C X_{0}\eps[m]
  + \varepsilon^{\gamma} Q\left(X_{0}\eps + Y_{0}\eps \right)
  \cdot \left\{ X_{T}\eps[m] + Y_{T}\eps[m] \right\}
  \\ \notag
  & \quad
  + P_{m} \bigl( X_{T}\eps[m-1] + Y_{T}\eps[m-1] \bigr),
  \\
  \label{YmEst}
  Y_{T}\eps[m] &\le C \left\{ Y_{0}\eps[m] + \varepsilon X_{0}\eps[m] 
  Y_{0}\eps[m] \right\}
  + C T^{\delta} \left\{ 1 + 2C \twonorm{\psi_{0}}{} \right\}^{2} Y_{T}\eps[m]
  \\ \notag
  & \quad
  + \varepsilon^{\gamma} Q\left(X_{0}\eps + Y_{0}\eps \right)
  \cdot \left\{ X_{T}\eps[m] + Y_{T}\eps[m] \right\}
  \\ \notag
  & \quad
  + P_{m} \bigl( X_{T}\eps[m-1] + Y_{T}\eps[m-1] \bigr).
\end{align}
Granting this for the moment, let us prove Proposition 
\ref{HigherOrderProp} by induction on $m$. The case $m = 0$ of 
\eqref{HigherOrderBound} was proved in Sect.\ \ref{UniformBounds}.
Adding up the inequalities \eqref{XmEst} and \eqref{YmEst}, we see 
that if \eqref{HigherOrderBound} holds for $m-1$, then it also holds 
for $m$, provided $T, \varepsilon > 0$ are so small that 
\eqref{Smallness} holds. Arguing as in the paragraph following 
Theorem \ref{UniformBoundsTheorem}, we iterate this argument to cover the 
full time interval $[0,T\subeps]$.

To prove the claim, we apply $\sum_{\abs{\alpha} \le m} 
\partial_{x}^{\alpha}$ to the system, and imitate the proof of the 
estimates in Proposition \ref{XYZProp} for $m = 0$. We single out the 
top order terms where $m$ derivatives fall on one of the fields 
$\A,\psi$ or $\phi_{\pm}$; these are estimated exactly like in the 
case $m =0$. All other terms are lumped together and yield the term 
$$
  P_{m} \bigl( X_{T}\eps[m-1] + Y_{T}\eps[m-1] \bigr).
$$
We skip the straightforward but tedious details of this argument.

\section{Estimates for the small component}\label{SmallComponentSection}

In this section we prove that if the ``positron part''
$\Pi_{-}\eps \psi\eps$ is small initially, then it stays small 
uniformly in every finite time interval, where ``small'' 
means either $O(\varepsilon)$ or $O(\varepsilon^{2})$.
Here is the precise result:

\begin{proposition}\label{SmallnessProp}
\begin{enumerate}
  \item
  Assume \eqref{HigherOrderDataBound} holds for some $m \ge 0$.
  Then if
  \begin{equation}\label{Small1}
    \Sobnorm{\Pi_{-}\eps \psi\eps}{m} = O(\varepsilon)
  \end{equation}
  holds at time $t = 0$, it also holds uniformly in every finite time 
  interval.
  \item
  Now replace \eqref{HigherOrderDataBound} by the stronger condition
  \begin{equation}\label{StrongDataBound}
    \Sobnorm{\psi_{0}\eps}{m+1} = O(1),
    \quad
    \Sobnorm{\nabla \mathbf a_{0}\eps}{m}
    + \varepsilon \Sobnorm{\mathbf a_{1}\eps}{m} = O(1),
  \end{equation}
  as $\varepsilon \to 0$. Then if
  \begin{equation}\label{Small2}
    \Sobnorm{\Pi_{-}\eps \psi\eps}{m-1} = O(\varepsilon^{2})
  \end{equation}
  holds at time $t = 0$, it also holds uniformly in every finite time 
  interval.
\end{enumerate}
\end{proposition}

Let us interpret this result in terms of $\eta\eps$, the lower
component of $e^{it/\varepsilon^{2}} \psi\eps$, as in 
\eqref{UpperLower2}. We claim that \eqref{Small1} is equivalent to
\begin{equation}\label{Small1'}
  \Sobnorm{\eta\eps}{m} = O(\varepsilon)
\end{equation}
while \eqref{Small2} is equivalent to \eqref{Small1'} and
\begin{equation}\label{Small2'}
  \Sobnorm{\partial_{t} \eta\eps}{m-1} = O(1).
\end{equation}
The equivalence of \eqref{Small1} and \eqref{Small1'} follows 
from \eqref{ProjExpansion1}, since
\begin{equation}\label{HighOrderPsiBound}
  \Sobnorm{\psi\eps}{m+1} = O(1)
\end{equation}
on account of Proposition \ref{HigherOrderProp}. To prove the rest of 
the claim, note that by \eqref{ProjExpansion2}, \eqref{Small2} is equivalent to
\begin{equation}
  \label{Small2''}
  \sigma^{j} \partial_{j} \eta\eps = O(\varepsilon),
  \qquad
  \eta\eps + i\varepsilon \frac{1}{2} \sigma^{j} \partial_{j} \chi\eps
  = O(\varepsilon^{2})
  \qquad \text{in} \quad H^{m-1}
\end{equation}
where $\chi\eps$ is the upper component of
$e^{it/\varepsilon^{2}} \psi\eps$, as in \eqref{UpperLower2}.
But by the second equation in \eqref{UpperLowerDirac},
\begin{equation}\label{LowerDiracApprox}
  i \varepsilon^{2} \partial_{t} \eta\eps
  = \eta\eps + i\varepsilon \frac{1}{2} \sigma^{j} \partial_{j} 
  \chi\eps + O(\varepsilon^{2}) \qquad \text{in} \quad H^{m-1}
\end{equation}
where we used the fact, proved below, that if \eqref{Small1} and
\eqref{StrongDataBound} hold initially, then
\begin{equation}\label{ImprovedABound}
  \Sobnorm{\nabla \A\eps}{m}
  + \varepsilon \Sobnorm{\mathbf \partial_{t} \A\eps}{m}
  = O(1)
\end{equation}
uniformly in every finite time interval.

\begin{proof}[Proof of Proposition \ref{SmallnessProp}(i).]
In view of \eqref{DiracvsKGsplitting}--\eqref{ErrorTerm}, we can 
replace $\Pi_{-}\eps \psi\eps$ by $\phi_{-}$ in \eqref{Small1}, and 
by \eqref{HighFreqEst} it suffices to consider the low frequency part
$(\phi_{-})\low$. Set (cf.\ Definition \ref{NormDef})
$$
  \Tilde Z_{T}\eps[m] = \sum_{\abs{\alpha} \le m}
  \norm{\partial_{x}^{\alpha}(\phi_{-}\eps)\low}_{\mixed{2}{6} \cap 
  \mixed{\infty}{2}(S_{T})}.
$$
Then recalling Proposition \ref{HigherOrderProp} and using induction on 
$m$, it suffices to prove that there exist constants $C,\delta > 0$
and polynomials $P_{m}$---all independent of $\varepsilon$---such that
for $0 \le T \le 1$ and $m \ge 0$,
\begin{multline*}
  \Tilde Z_{T}\eps[m] \le C \Sobnorm{\phi_{-}(t=0)}{m}
  + C T^{\delta} \left\{ C \twonorm{\psi_{0}\eps}{} + 1 \right\}^{2} \Tilde 
  Z_{T}\eps[m]
  \\
  + \left\{ \varepsilon + \Tilde Z_{T}\eps[m-1] \right\} 
  P_{m} \bigl( X_{T}\eps[m] + Y_{T}\eps[m] \bigr),
\end{multline*}
where by convention $\Tilde Z_{T}\eps[-1] = 0$. But this estimate follows by a 
straightforward modification of the proof of the estimate for the 
$Z$-norm in Proposition \ref{XYZProp}, taking into account the bound 
\eqref{FinalZBound}.

Let us now prove \eqref{ImprovedABound}, assuming it holds initially. 
In view of Lemma \ref{HomWaveCutoffLemma}, this reduces to proving
$\LpHslocal{\J\eps}{2}{m} = O(1)$.
Split $\J = \J' + \J''$ as in Sect.\ \ref{CurrentEstimates}. To 
estimate $\J''$, we proceed as in the proof of \eqref{JhighEst}, taking into 
account the higher order bound \eqref{HighOrderPsiBound}.
On the other hand, since $\J'$ has vector components $\frac{2}{\varepsilon} \re 
\innerprod{\sigma^{j} (\chi\eps)\low}{(\eta\eps)\low}$, we have
$$
  \twonorm{\partial_{x}^{\alpha} \J'}{} \le \frac{1}{\varepsilon}
  \sum_{\beta+\gamma = \alpha} c_{\alpha \beta} 
  \mixednorm{\partial_{x}^{\beta} (\chi\eps)\low}{2}{\infty}
  \mixednorm{\partial_{x}^{\gamma} (\eta\eps)\low}{\infty}{2},
$$
and the r.h.s.\ is $O(1)$ for $\abs{\alpha} \le m$ on account of \eqref{Small1'},
\eqref{SchrEmb2} and \eqref{HigherOrderBound}.
\end{proof}

\begin{proof}[Proof of Proposition \ref{SmallnessProp}(ii).] Here we 
break with our earlier notation, writing
$$
  \psi_{\pm}\eps = \Pi_{\pm}\eps \psi\eps,
  \qquad
  \phi\eps = e^{it/\varepsilon^{2}} \psi \eps,
  \qquad
  \phi_{\pm}\eps = e^{it/\varepsilon^{2}} \psi_{\pm}\eps.
$$
Then from the Dirac equation,
\begin{align}
  \label{UpperDirac2}
  i\partial_{t}\phi_{+}\eps - \frac{\lambda\eps - 1}{\varepsilon^{2}} 
  \phi_{+}\eps + \Pi_{+}\eps \left( \mathcal A\eps \phi\eps \right) &= 0,
  \\
  \label{LowerDirac2}
  i\partial_{t}\phi_{-}\eps + \frac{\lambda\eps + 1}{\varepsilon^{2}} 
  \phi_{-}\eps + \Pi_{-}\eps \left( \mathcal A\eps \phi\eps \right) &= 0,
\end{align}
where $\mathcal A\eps = A_{j}\eps \sigma^{j} + A_{0}\eps$. Thus,
\begin{equation}\label{PhiMinus}
  \phi_{-}\eps = \left( \lambda\eps + 1 \right)^{-1} \varepsilon^{2} 
  \left\{ - i \partial_{t} \phi_{-}\eps - \Pi_{-}\eps \left( \mathcal 
  A\eps \phi\eps \right) \right\},
\end{equation}
so we reduce \eqref{Small2} to proving
\begin{gather}
  \label{Small2A}
  \Sobnorm{\partial_{t} \phi_{-}\eps}{m-1} = O(1),
  \\
  \label{Small2B}
  \Sobnorm{\Pi_{-}\eps \left( \mathcal 
  A\eps \phi\eps \right)}{m-1} = O(1).
\end{gather}
The latter follows readily from \eqref{ImprovedABound} and
\eqref{HighOrderPsiBound}, since $\Pi_{-}\eps$ is uniformly 
bounded. For later use we also note that \eqref{UpperDirac2} implies
\begin{equation}\label{Small2APlus}
  \Sobnorm{\partial_{t} \phi_{+}\eps}{m-1} = O(1),
\end{equation}
since the symbol of $\frac{\lambda\eps - 1}{\varepsilon^{2}}$ is 
bounded by $\abs{\xi}^{2}$.

To prove \eqref{Small2A} we proceed as in \cite[Sect.\ 4]{BMP}.
Consider first the case $m = 1$. Take a time derivative of 
\eqref{LowerDirac2}, then take the imaginary part of its inner 
product with $\partial_{t} \phi_{-}\eps$ and integrate in $x$.
Making use of the 
self-adjointness of $\lambda\eps$, $\Pi_{-}\eps$ and $\mathcal 
A\eps$, and the fact that $(\Pi_{-}\eps)^{2} = \Pi_{-}\eps$, we then 
obtain
\begin{align*}
  \frac{1}{2} \frac{d}{dt} \twonorm{\partial_{t} \phi_{-}\eps}{}^{2}
  &\le \twonorm{\mathcal A\eps \partial_{t} \phi_{+}\eps}{}
  \twonorm{\partial_{t} \phi_{-}\eps}{}
  + \twonorm{\partial_{t} \mathcal A\eps \cdot \phi\eps}{}
  \twonorm{\partial_{t} \phi_{-}\eps}{}
  \\
  &\lesssim
  \left( \Sobnorm{\nabla \mathcal A\eps}{1} \twonorm{\partial_{t} \phi_{+}\eps}{}
  + \twonorm{\nabla \partial_{t} \mathcal A\eps}{} \Sobnorm{\phi\eps}{1}
  \right) \twonorm{\partial_{t} \phi_{-}\eps}{}.
\end{align*}
Thus, dividing by $\twonorm{\partial_{t} \phi_{-}\eps}{}$ and 
integrating in time, and using the fact that \eqref{Small2A} holds at 
time $t = 0$ (this follows from \eqref{LowerDirac2} and the initial 
assumptions), we reduce \eqref{Small2A} for $m = 1$ to proving that the terms 
inside the parentheses in the last inequality above are all $O(1)$ 
locally uniformly in time. But this follows from the bounds 
\eqref{HighOrderPsiBound}, \eqref{ImprovedABound} and \eqref{Small2APlus}.
Here we use also the fact that $\partial_{t} A_{0}\eps$ enjoys the 
same bounds as $\nabla A_{0}\eps$, in view of the equation
$$
  \Delta \partial_{t} A_{0}\eps = - \dv \J\eps
$$
which follows from \eqref{A0Eq} and the conservation law 
$\partial_{t} \rho\eps + \dv \J\eps = 0$.

Following \cite{BMP} we now proceed by induction on $m$, starting at 
$m = 1$. Thus, we apply $\partial_{t} \partial_{x}^{\alpha}$, where
$\abs{\alpha} \le m-1$,  to the equation \eqref{LowerDirac2},
and we take the imaginary part of its inner product with $\partial_{t} 
\partial_{x}^{\alpha} \phi_{-}\eps$ and integrate in $x$.
Then by a straightforward modification of the argument for $m = 1$, we
reduce \eqref{Small2A} to the $O(1)$ bounds \eqref{HighOrderPsiBound},
\eqref{ImprovedABound} and \eqref{Small2APlus}, as well as
\eqref{Small2A} at the previous induction step. We omit the details. 
\end{proof}

\section{Nonrelativistic limit}

We first prove Theorem \ref{NonrelLimitThm}, then we discuss the 
modifications needed to prove Theorem \ref{NonrelLimitThm2}.

\begin{proof}[Proof of (\ref{SpinorConv}).]
This can be restated:
\begin{equation}\label{NonrelLimitConv2}
  e^{it/\varepsilon^{2}} \Pi_{+}^{0} \psi\eps
  \to
  \begin{pmatrix}
    v_{+}  \\
    0
  \end{pmatrix},
  \quad
  e^{-it/\varepsilon^{2}} \Pi_{-}^{0} \psi\eps
  \to 
  v_{-} =
  \begin{pmatrix}
    0  \\
    v_{-}
  \end{pmatrix}
  \quad
  \text{in $H^{1}$ as $\varepsilon \to 0$}
\end{equation}
locally uniformly in time. We claim it suffices to prove
\begin{equation}\label{NonrelLimitConv3}
  e^{it/\varepsilon^{2}} \Pi_{+}\eps \psi\eps
  \to
  \begin{pmatrix}
    v_{+}  \\
    0
  \end{pmatrix},
  \quad
  e^{-it/\varepsilon^{2}} \Pi_{-}\eps \psi\eps
  \to 
  v_{-} =
  \begin{pmatrix}
    0  \\
    v_{-}
  \end{pmatrix}
  \quad
  \text{in $H^{1}$ as $\varepsilon \to 0$.}
\end{equation}
To prove the claim, write $\Pi_{\pm}\eps \psi\eps = \Pi_{\pm}^{0} \psi\eps 
\pm r\eps$.
By the orthogonality between $\Pi_{+}^{0}$ and $\Pi_{-}^{0}$, we get
$$
  \Pi_{+}^{0} r\eps
  = - \Pi_{+}^{0} \Pi_{-}\eps \psi\eps,
  \qquad
  \Pi_{-}^{0} r\eps
  = \Pi_{-}^{0} \Pi_{+}\eps \psi\eps.
$$
But if \eqref{NonrelLimitConv3} holds, then the right hand sides 
converge to zero in $H^{1}$. Thus $r\eps = o(1)$ in $H^{1}$ and 
we have proved that \eqref{NonrelLimitConv3} implies \eqref{NonrelLimitConv2}.
In the remainder of the proof we skip the superscript $\varepsilon$ 
on the fields, to simplify the notation.

Using \eqref{DiracvsKGsplitting} and \eqref{ErrorTerm} we reduce 
\eqref{NonrelLimitConv3} to proving
\begin{equation}\label{Conv4}
  \phi_{+}
  \longrightarrow
  \begin{pmatrix}
    v_{+}  \\
    0
  \end{pmatrix},
  \quad
  \phi_{-}
  \longrightarrow
  \begin{pmatrix}
    0  \\
    v_{-}
  \end{pmatrix}
  \quad
  \text{in} \quad H^{1}
  \quad \text{as} \quad \varepsilon \longrightarrow 0,
\end{equation}
uniformly in any given time interval $[0,T]$. By (the proof of) Theorem 
\ref{LWPThm}, the solution exists in this time interval for all 
sufficiently small $\varepsilon > 0$, and
\begin{equation}\label{Bounds}
  X_{T}\eps + Y_{T}\eps = O(1),
  \qquad
  \bigXnormArg{e^{\pm it/\varepsilon^{2}} R\eps}{1}{\theta-1}{S_{T}}
  = o(1)
\end{equation}
as $\varepsilon \to 0$.
Note that \eqref{Conv4} holds at time $t = 0$, by Lemma \ref{DataLemma}.
Thus, it suffices to prove that there exist
$K,\delta > 0$, depending on $T$ and $X_{T}\eps + Y_{T}\eps$,
but independent of $\varepsilon$, such that for every time interval
$I = [t_{0},t_{1}] \subset [0,T]$,
\begin{equation}\label{Bootstrap}
  f(I) \le K f( \{t_{0}\} )
  + K \abs{I}^{\delta} f(I) + o(1)
\end{equation}
as $\varepsilon \to 0$, where
\begin{equation}\label{fDef2}
  f(I) = \norm{\phi_{+} - \begin{pmatrix}
    v_{+}  \\
    0
  \end{pmatrix} }_{L_{t}^{\infty}H^{1}(I\times\R^{3})}
  +
  \norm{\phi_{-} - \begin{pmatrix}
    0  \\
    v_{-}
  \end{pmatrix} }_{L_{t}^{\infty}H^{1}(I\times\R^{3})}.
\end{equation}

W.l.o.g.\ we assume $I = [0,T]$, and we only estimate the first term in
\eqref{fDef2}. Write
\begin{align*}
  \phi_{+}(t) &= U\eps(t)\phi_{0}^{+} + \int_{0}^{t} U\eps(t-s) \left[ 
  L_{+}\eps \phi_{+} (s) \right] \, ds,
  \\
  v_{+}(t) &= S(t) v_{0}^{+} - \int_{0}^{t} S(t-s) \left[ 
  (uv_{+})(s) \right] \, ds,
\end{align*}
where $\phi_{0}^{+}, v_{0}^{+}$ are the data of 
$\phi_{+}, v_{+}$ and $U\eps(t), S(t)$ are given by \eqref{Propagators}.
Thus
\begin{equation}\label{DifferenceExpansion}
\begin{split}
  \phi_{+}(t) - \begin{pmatrix}
    v_{+}  \\
    0
  \end{pmatrix}(t)
  &=
  U\eps(t) \left[ \phi_{0}^{+} - \begin{pmatrix}
    v_{0}^{+} \\
    0
  \end{pmatrix}
  \right]
  +
  \left[ U\eps(t) - S(t) \right] \begin{pmatrix}
    v_{0}^{+}  \\
    0
  \end{pmatrix}
  \\
  &\quad
  +
  \int_{0}^{t} U\eps(t-s) \left[
  \begin{pmatrix}
    uv_{+}  \\
    0
  \end{pmatrix}(s)
  + L_{+}\eps \phi_{+}(s) \right] \, ds
  \\
  &\quad
  +
  \int_{0}^{t} \left[ S(t-s) - U\eps(t-s) \right] \begin{pmatrix}
    uv_{+}  \\
    0
  \end{pmatrix}(s)
  \, ds
  \\
  &= I_{1} + I_{2} + I_{3} + I_{4}.
\end{split}
\end{equation}
Clearly,
$$
  \LpHslocal{I_{1}}{\infty}{1} \lesssim
  \Sobnorm{\phi_{0}^{+} - \begin{pmatrix}
    v_{0}^{+}  \\
    0
  \end{pmatrix} }{1}.
$$
As in \cite[Sect.\ 5]{BMS},
$$
  \LpHslocal{I_{j}}{\infty}{1} = o(1) \quad \text{for} \quad j = 2,4,
$$
using the dominated convergence theorem and the fact that
\begin{equation}\label{uBound}
  \Lxpnorm{\nabla u}{3} + \Lxpnorm{u}{\infty} \lesssim
  \Sobnorm{v_{+}}{1} + \Sobnorm{v_{-}}{1}
  < \infty
\end{equation}
uniformly in every finite time interval.
It remains to consider $I_{3}$. By Lemma \ref{XCutoffLemma}
and the embeddings \eqref{BasicEmbedding} and \eqref{InterpolationEmbeddings},
\begin{equation}\label{I3bound1}
  \LpHslocal{I_{3}}{\infty}{1} \lesssim
  \LpHslocal{\begin{pmatrix}
    uv_{+}  \\
    0
  \end{pmatrix}
  - A_{0} \phi_{+} }{2}{1}
  + \bigXnormArg{e^{\pm it/\varepsilon^{2}} R\eps}{1}{\theta-1}{S_{T}}.
\end{equation}
The second term on the r.h.s.\ is $o(1)$ by \eqref{Bounds},
and the first term is bounded by
$$
  T^{1/2} \left(
  \LpHslocal{
  u \left\{ \begin{pmatrix}
    v_{+}  \\
    0
  \end{pmatrix}
  - \phi_{+} \right\} }{\infty}{1}
  +
  \LpHslocal{
  (u - A_{0}) \phi_{+} }{\infty}{1}
  \right).
$$
But using Leibniz' rule, H\"older's inequality and
Sobolev embedding, it is easy to see that 
the terms inside the parentheses are dominated by $K f(I)$,
where $K$ depends on the size of $X_{T}\eps + Y_{T}\eps$ and \eqref{uBound}.
\end{proof}

\begin{proof}[Proof of (\ref{A0Conv}) and (\ref{ChargeConv}).]
Using Sobolev embedding we reduce \eqref{A0Conv} to \eqref{ChargeConv}.
To prove the latter, observe that \eqref{Conv4} implies
\begin{equation}\label{Conv5}
  \chi_{+} \to v_{+},
  \quad \chi_{-} \to 0,
  \quad \eta_{+} \to 0,
  \quad \eta_{-} \to v_{-}
  \quad
  \text{in $H^{1}$ as $\varepsilon \to 0$,}
\end{equation}
locally uniformly in time.
Thus \eqref{ChargeConv} follows immediately from \eqref{ChargeExpansion} 
using H\"older's inequality and Sobolev embedding. 
\end{proof}

\begin{proof}[Proof of (\ref{CurrentConv}).]
Multiply \eqref{CurrentExpansion} by a $C^{1}$ compactly supported
test function $G(t,x)$ and integrate in $t,x$.
W.l.o.g.\ assume $G$ is real-valued.
The integrals corresponding to the last two 
terms in r.h.s.\eqref{CurrentExpansion} are $O(\varepsilon)$ in 
absolute value. To see this, integrate by parts in time and use
\begin{equation}\label{SobIneq}
  \abs{\int fg \, dx} \le \Sobnorm{f}{-1} \Sobnorm{g}{1}
\end{equation}
and the bound, locally uniform in time,
\begin{equation}\label{TimeDerBound}
  \Sobnorm{\partial_{t} \phi_{\pm}}{-1} = O(1).
\end{equation}
The latter is easily reduced to the 
uniform bounds for $X_{T}\eps + Y_{T}\eps$, using Lemma
\ref{ModifiedDiracLemma}, Sobolev embedding and H\"older's inequality.

Next, fix $1 \le j \le 3$ and consider
$$
  I_{\pm} := \frac{2}{\varepsilon} \re \int
  \innerprod{\sigma^{j} \chi_{\pm}}{\eta_{\pm}}
  G \, dt \, dx.
$$
In view of \eqref{Conv5} and \eqref{SmallComponentsBound},
\begin{equation}\label{Iminus1}
  I_{-} = \frac{2}{\varepsilon} \re \int
  \innerprod{\sigma^{j} \chi_{-}}{v_{-}}
  G \, dt \, dx
  + o(1).
\end{equation}
By \eqref{DiracvsKGsplitting},\eqref{ErrorTerm} and \eqref{Bounds},
\begin{equation}\label{Iminus2}
  \frac{1}{\varepsilon} \begin{pmatrix}
    \chi_{-} \\
    0
  \end{pmatrix}
  = e^{-it/\varepsilon^{2}} \frac{1}{\varepsilon} \Pi_{+}^{0} 
  \Pi_{-}\eps \psi + O(\varepsilon^{2-\fixedpar})
  \quad \text{in} \quad L_{x}^{2}
\end{equation}
locally uniformly in time. But by \eqref{ProjExpansion3},
$$
  \frac{1}{\varepsilon} \Pi_{+}^{0} \Pi_{-}\eps \psi
  = \frac{i}{2} [\lambda\eps]^{-1} 
  \begin{pmatrix}
    \sigma^{k} \partial_{k} \LowerSpinor  \\
    0
  \end{pmatrix}
  + \frac{1}{2\varepsilon} \left( 1 - [\lambda\eps]^{-1} \right)
  \begin{pmatrix}
    \UpperSpinor  \\
    0
  \end{pmatrix}.
$$
In view of \eqref{SpinorConv} and the bound \eqref{SymbolBound}, it follows that
\begin{multline}\label{Iminus3}
  e^{-it/\varepsilon^{2}} \frac{1}{\varepsilon} \Pi_{+}^{0} 
  \Pi_{-}\eps \psi
  \\
  = \frac{i}{2} [\lambda\eps]^{-1} 
  \begin{pmatrix}
    \sigma^{k} \partial_{k} v_{-} \\
    0
  \end{pmatrix}
  + \frac{1}{2\varepsilon} \left( 1 - [\lambda\eps]^{-1} \right)
  \begin{pmatrix}
    e^{-2it/\varepsilon^{2}} v_{+} \\
    0
  \end{pmatrix}
  + o(1)
\end{multline}
in $L_{x}^{2}$. Moreover, by dominated convergence,
\begin{equation}\label{Iminus4}
  [\lambda\eps]^{-1} \sigma^{k} \partial_{k} v_{-}
  = \sigma^{k} \partial_{k} v_{-} + o(1)
  \quad \text{in} \quad  H^{-1}.
\end{equation}
Using \eqref{Iminus1}--\eqref{Iminus4} and either H\"older's 
inequality or \eqref{SobIneq}, we conclude that
$$
  I_{-} = \re \int
  i \innerprod{\sigma^{j} \sigma^{k} \partial_{k} v_{-}}{v_{-}}
  G \, dt \, dx
  + I_{-}'
  + o(1)
$$
where
$$
  I_{-}'
  = \frac{1}{\varepsilon}
  \re \int e^{-2it/\varepsilon^{2}}
  \innerprod{\sigma^{j} \left( 1 - [\lambda\eps]^{-1} \right) v_{+}}{v_{-}}
  G \, dt \, dx.
$$
But the latter is $O(\varepsilon)$ in absolute value (integrate by 
parts in time and use the analogue of \eqref{TimeDerBound} for 
$v_{\pm}$). Using \eqref{alphaIdentities} we finally conclude that
$$
  I_{-} = \int
  \left\{
  - \im \innerprod{\partial_{j} v_{-}}{v_{-}}
  - \frac{1}{2} \epsilon^{jkl} \partial_{k}
   \innerprod{\sigma_{l} v_{-}}{v_{-}}
  \right\}
  G \, dt \, dx + o(1).
$$
A similar calculation can be done for $I_{+}$, and this proves 
\eqref{CurrentConv}.
\end{proof}

Next, we prove Theorem \ref{NonrelLimitThm2}. By hypothesis, 
\eqref{Small1'}, or equivalently \eqref{Small1}, holds initially
and therefore also uniformly in every finite time interval, by 
Proposition \ref{SmallnessProp}. Next observe that since
\eqref{ImprovedConvergence} holds initially, we have
$$
  \Pi_{\pm}\eps \psi\eps = \psi_{\pm}\eps + O(\varepsilon)
  \qquad \text{in} \quad H^{1}
$$
locally uniformly in time.
In fact, this follows from 
\eqref{DiracvsKGsplitting}--\eqref{ErrorTerm}, since 
\eqref{HigherOrderBound} holds with $m = 1$. We conclude that
it suffices to prove \eqref{ImprovedConvergence} with $\psi\eps$ 
replaced by $\psi_{+}\eps$. We proceed as in the proof of Theorem 
\ref{NonrelLimitThm}, but now the remainder term in \eqref{Bootstrap} must 
be improved from $o(1)$ to $O(\varepsilon)$, and $f(I)$ is given by the 
first term in r.h.s.\eqref{fDef2}. Again we reduce to estimating 
the terms $I_{1},\dots,I_{4}$ as given by \eqref{DifferenceExpansion}.

The term $I_{1}$ is estimated exactly as before, but is now 
$O(\varepsilon)$ since \eqref{ImprovedConvergence} is assumed to hold
initially. Using the fact that $U\eps(t) - S(t) = \varepsilon^{2} \mathcal 
R_{4}\eps(t)$, where $\mathcal R_{4}\eps(t)$ is bounded from $H^{s+4} \to 
H^{s}$ uniformly in $\varepsilon$ and $0 \le t \le T$,
and the assumption that the initial datum of $v_{+}$ is in $H^{5}$, we 
find that
$$
  \LpHslocal{I_{j}}{\infty}{1} = O(\varepsilon^{2})
  \quad \text{for} \quad j = 2,4.
$$
For $I_{3}$ we use again \eqref{I3bound1}, but now the last term is 
$O(\varepsilon)$, as follows from the proof of \eqref{REstC} 
taking into account the fact that \eqref{HigherOrderBound} holds with $m = 1$.
The first term in r.h.s.\eqref{I3bound1} is estimated exactly as before.
This proves \eqref{ImprovedConvergence}, and then it follows 
immediately that \eqref{A0Conv} and \eqref{ChargeConv} are also 
improved to $O(\varepsilon)$. 

\section{Semi-nonrelativistic limit}\label{SeminonrelLimit}

Here we prove Theorem \ref{SeminonrelLimitThm}.
The initial assumptions (i), (ii) imply, as proved in Sects.\ 
\ref{HigherBoundsSection} and \ref{SmallComponentSection}, that
\eqref{HigherOrderBound} holds with $m = 4$, while
\eqref{Small1}--\eqref{ImprovedABound} hold with $m = 2$.
We write
$$
  \phi\eps = 
  \begin{pmatrix}
    \chi\eps  \\
    \eta\eps
  \end{pmatrix}
  := e^{it/\varepsilon^{2}} \psi\eps,
  \qquad
  \phi_{\pm}\eps := e^{it/\varepsilon^{2}} \psi_{\pm}\eps,
$$
with $\psi_{\pm}\eps$ defined as in \eqref{KGsplitting}. Also, we 
denote by $\chi_{\pm}\eps$ the upper component of $\phi_{\pm}\eps$.
Observe that
$$
  \Pi_{\pm}\eps \psi\eps = \psi_{\pm}\eps + O(\varepsilon^{2})
  \qquad \text{in} \quad H^{1},
$$
in view of \eqref{DiracvsKGsplitting}--\eqref{ErrorTerm} and the 
bounds \eqref{HighOrderPsiBound} and \eqref{ImprovedABound} for $m = 2$.
On account of \eqref{Small2}, we may therefore replace $\chi\eps$ in 
\eqref{PauliConv} by $\chi_{+}\eps$. By Lemma \ref{ModifiedDiracLemma},
$$
  \left( i\partial_{t} - \frac{\lambda\eps - 1}{\varepsilon^{2}} 
  \right) \chi_{+}\eps + A_{0}\eps \chi_{+}\eps
  = \Tilde R\eps,
$$
where
\begin{align*}
  \Tilde R\eps &=
  \frac{1}{2} \varepsilon i \A\eps \cdot \nabla \chi\eps
  - \frac{1}{2} \varepsilon B_{j}\eps \sigma^{j} \chi\eps
  + \frac{1}{2} \varepsilon^{2} \left(\A\eps\right)^{2} \chi\eps
  \\
  &\quad
  - \frac{1}{2} \left( 1 - [\lambda\eps]^{-1} \right) \varepsilon
  \left\{ 2i \A\eps \cdot \nabla \chi\eps - B_{j}\eps \sigma^{j} \chi\eps 
  \right\}
  + \frac{1}{2} [\lambda\eps]^{-1} \varepsilon \left\{ i E_{j}\eps \sigma^{j} \eta\eps 
  \right\}
  \\
  &\quad - \frac{1}{2} \left( 1 - [\lambda\eps]^{-1} \right) \varepsilon^{2} \left\{ 
  \left(\A\eps\right)^{2} \chi\eps \right\}
  - \frac{1}{2} \varepsilon^{2} [\lambda\eps]^{-1} \left[ A_{0}\eps, 
  \frac{\lambda\eps - 1}{\varepsilon^{2}} \right] \left( \chi_{+}\eps - \chi_{-}\eps \right).
\end{align*}
Recalling the bound \eqref{SymbolBound} on the symbol of $1 - 
[\lambda\eps]^{-1}$ and using the fact that \eqref{HigherOrderBound}, 
\eqref{Small1'} and \eqref{ImprovedABound} hold with $m = 2$, we 
conclude that
$$
  \Tilde R\eps =
  \frac{1}{2} \varepsilon i \A\eps \cdot \nabla \chi\eps
  - \frac{1}{2} \varepsilon B_{j}\eps \sigma^{j} \chi\eps
  + \frac{1}{2} \varepsilon^{2} \left(\A\eps\right)^{2} \chi\eps
  + O(\varepsilon^{2}) \quad \text{in} \quad H^{1},
$$
locally unformly in time.

Then, since \eqref{HigherOrderBound} holds with $m = 4$ and
$\frac{\lambda\eps - 1}{\varepsilon^{2}} = \frac{\Delta}{2}
+ \varepsilon^{2} \mathcal R_{4}$, where $\mathcal R_{4}\eps$ is bounded from
$H^{s+4} \to H^{s}$ uniformly in $\varepsilon$, we further
conclude that
\begin{equation}\label{ApproxPauli2}
  i \partial_{t} \chi_{+}\eps = \frac{1}{2} \left( i \nabla + 
  \varepsilon \A\eps \right)^{2} \chi_{+}\eps - A_{0}\eps \chi_{+}\eps
  - \frac{1}{2} \varepsilon B_{j}\eps \sigma^{j} \chi_{+}\eps
  + \varepsilon^{2} r\eps
\end{equation}
where $r\eps = O(1)$ in $H^{1}$ locally unformly in time.
Comparing \eqref{ApproxPauli2} to the Pauli equation \eqref{PauliEq} via
the energy inequality for the self-adjoint ``Pauli operator'',
$$
  P\eps =  \frac{1}{2} \left( i \nabla + 
  \varepsilon \A\eps \right)^{2}
  - \frac{1}{2} \varepsilon B_{j}\eps \sigma^{j},
$$
one finds that
$$
  f(I) \le f( \{t_{0}\} )
  + K \abs{I} f(I) + O(\varepsilon^{2})
$$
as $\varepsilon \to 0$, where
$$
  f(I) = \norm{ \chi_{+}\eps - \Pauli\eps }_{L_{t}^{\infty} 
  H^{1} (I \times \R^{3}) }
$$
for time intervals $I = [t_{0},t_{1}] \subset [0,T]$, and where 
$K$ depends on $T$ but not on $\varepsilon$. In fact, $K$ depends on 
the $O(1)$ bounds in \eqref{HighOrderPsiBound} and \eqref{ImprovedABound},
which hold for $m = 2$ as we recall. We conclude that
$f([0,T]) = O(\varepsilon^{2})$, and this proves \eqref{PauliConv}.

Observe that 
\eqref{InitialConstraint} holds in $H^{1}$ locally uniformly in time, 
in view of \eqref{LowerComponentExpansion} and the fact that 
\eqref{Small1'}, \eqref{Small2'} and \eqref{ImprovedABound} hold for 
$m = 2$. Substituting \eqref{InitialConstraint} into
$$
  \J\eps = \varepsilon^{-1} \left\{ 2 \re
  \innerprod{\sigma^{k}\chi\eps}{\eta\eps}_{\C^{2}} \right\}_{k=1,2,3}
$$
and using \eqref{PauliConv} yields \eqref{SeminonrelCurrent}.

\section{Proofs of the spacetime estimates}\label{BilinearProofs}

Here we prove Theorem \ref{DyadicThm} and Proposition \ref{ImprovedStrProp}.

\begin{proof}[Proof of Proposition \ref{ImprovedStrProp}]
Let $Q$ be a cube with side length $\sim \mu$ centered at $\xi_{0}$, 
where $\abs{\xi_{0}} \sim \lambda$, and let $\chi_{Q}(\xi)$ be a 
smooth cut-off function equal to $1$ on $Q$. For example, 
we can take
\begin{equation}\label{CubeCutoff}
  \chi_{Q}(\xi) := \eta \left( \frac{\xi-\xi_{0}}{\mu} \right),
\end{equation}
where $\eta$ is a smooth bump function equal to $1$ on a neighborhood of 
the origin. Then by the $TT^{*}$ method, we reduce 
\eqref{ImprovedStr} to the decay estimate
\begin{equation}\label{KernelEst}
  \abs{K_{\varepsilon,Q}(t,x)} \lesssim
  \begin{cases}
    \mu \abs{t}^{-1} &\quad \text{for $\lambda \lesssim 1/\varepsilon$},
    \\
    \varepsilon \mu \lambda \abs{t}^{-1} &\quad \text{for $\lambda \gg 
    1/\varepsilon$,}
  \end{cases}
\end{equation}
for the convolution kernel
$$K_{\varepsilon,Q}(t,x) := \int_{\R^{3}} e^{ix\cdot \xi} 
e^{ith\subeps(\xi)} \chi_{Q}(\xi) \, d\xi,$$ with $h\subeps$ given by 
\eqref{heps}. In view of the scaling identity $$K_{\varepsilon,Q}(t,x)
= \varepsilon^{-3} K_{1,\varepsilon Q}(\varepsilon^{-2} t, 
\varepsilon^{-1} x),$$ it suffices to prove \eqref{KernelEst} for 
$\varepsilon = 1$. To simplify the notation we write $K_{Q}$ 
instead of $K_{1,Q}$. Thus,
\begin{equation}\label{KernelPolar}
  K_{Q}(t,x) = \int_{0}^{\infty} \int_{S^{2}}  e^{i r x \cdot \omega}
  e^{it \alpha(r)} \chi_{Q}(r\omega) r^{2} \, d\sigma(\omega) \, dr
\end{equation}
where $\sigma$ is surface measure on $S^{2}$ and $\alpha$ is given 
by \eqref{AlphaDef}. Note that
\begin{equation}\label{AlphaDer}
  \alpha'(r) = \frac{r}{\sqrt{1+r^{2}}}
  \quad \text{and} \quad
  \alpha''(r) = \frac{1}{(1+r^{2})^{3/2}}.
\end{equation}

We split the problem into the following cases:
\begin{enumerate}
  \item\label{case1}
  $\lambda \lesssim 1$ and $\abs{x} \gtrsim \lambda\abs{t}$, 
  \item\label{case2}
  $\lambda \lesssim 1$ and $\abs{x} \ll \lambda \abs{t}$, 
  \item\label{case3}
  $\lambda \gg 1$ and $\abs{x} \gtrsim \abs{t}$, 
  \item\label{case4}
  $\lambda \gg 1$ and $\abs{x} \ll \abs{t}$.
\end{enumerate}
Rewrite \eqref{KernelPolar} as $K_{Q}(t,x) = \int_{0}^{\infty}
e^{it \alpha(r)} a(r,x) r^{2} \, dr$ where
$$
  a(r,x) := \int_{S^{2}}  e^{i r x \cdot \omega}
  \chi_{Q}(r\omega) \, d\sigma(\omega).
$$
We shall need the following:

\begin{lemma}\label{aDecayLemma}
$\abs{a(r,x)} \lesssim \mathbb(r\abs{x})^{-1} \chi_{I}(r)$, where 
$\chi_{I}$ is the characteristic function of an
interval $I$ of length $\sim \mu$ and centered at a distance 
$\sim \lambda$ from the origin.
\end{lemma}

\begin{proof}
The statement about the $r$-support of $a(r,x)$ is obvious, and the 
decay statement follows from the fact that
$$
  \abs{ \int_{S^{2}}  e^{i x \cdot \omega}
  \gamma(\omega) \, d\sigma(\omega) } \lesssim 1/\abs{x}
$$
for all smooth functions $\gamma$ such that $\abs{\gamma} \le 1$.
But this fact is easily proved by passing to spherical coordinates 
and rescaling.
\end{proof}

Thus $$\abs{K_{Q}(t,x)} \lesssim \int_{I} \bigl(r/\abs{x}\bigr) \, dr \sim \mu 
\lambda / \abs{x},$$ and this covers the cases \eqref{case1} and 
\eqref{case3} above.

To handle the remaining cases we write \eqref{KernelPolar} as
$K_{Q}(t,x) = \int_{S^{2}} b(\omega) \, d\sigma(\omega)$, where
$$
  b(\omega) = \int_{0}^{\infty} \frac{d}{dr} \left[ e^{i (t 
  \alpha(r) + r x\cdot \omega)} \right]
  \frac{\chi_{Q}(r\omega) r^{2}}{i \bigl( t \alpha'(r) + 
  x \cdot \omega \bigr)} \, dr.
$$
Integrate by parts and write
$$
  - \frac{d}{dr} \left[ \frac{\chi_{Q}(r\omega) r^{2}}
  {i \bigl( t \alpha'(r) +   x \cdot \omega \bigr)} \right]
  =  \frac{\chi_{Q}(r\omega) r^{2} t \alpha''(r)}
  {i \bigl( t \alpha'(r) + x \cdot \omega \bigr)^{2}}
  - \frac{\tfrac{d}{dr} \bigl[ \chi_{Q}(r\omega) r^{2} 
  \bigr]}{ i \bigl( t \alpha'(r) + 
  x \cdot \omega \bigr)}.
$$
Correspondingly we split $b = b_{1} + b_{2}$. Observe that the 
$r$-support of $\chi_{Q}(r\omega)$ is contained in an
interval $I$ of length $\sim \mu$ and centered at a distance 
$\lambda$ from the origin, while the $\omega$-support is contained in 
a set given by
\begin{equation}\label{Angle}
  \angle(\omega,\omega_{0}) \lesssim \mu/\lambda
\end{equation}
for some $\omega_{0} \in S^{2}$. Moreover, in view of \eqref{CubeCutoff} 
we have
\begin{equation}\label{CutoffDer}
  \abs{ \tfrac{d}{dr} \chi_{Q}(r\omega) } \lesssim 1/\mu.
\end{equation}

Now consider case \eqref{case4}. Then on account of 
\eqref{AlphaDer} we have $\alpha'(r) \sim 1$ and $\alpha''(r) \sim 
\lambda^{-3}$ for $r \in I$, so
$\abs{t \alpha'(r) +   x \cdot \omega} \gtrsim \abs{t}$. Thus
$$
  \abs{b_{1}(\omega)} \lesssim \bigl( 1 / \lambda \abs{t} \bigr) 
  \int_{I} \, dr \lesssim \mu / \lambda \abs{t},
$$
which is more than good enough. Next, using \eqref{CutoffDer} we have
$$
  \abs{b_{2}(\omega)} \lesssim \bigl( 1 / \abs{t} \bigr) 
  \int_{I} \left( r + r^{2}/\mu \right) \, dr
  \lesssim (\lambda \mu + \lambda^{2}) / \abs{t}
  \lesssim \lambda^{2} / \abs{t}.
$$
But integrating this over the region \eqref{Angle} on $S^{2}$ gives 
us a bound $\mu^{2} / \abs{t}$, which again is more than good enough.

Finally, consider case \eqref{case2}. Then $\alpha'(r) \sim \lambda$ and
$\alpha''(r) \sim 1$ for $r \in I$,
so $\abs{t \alpha'(r) +   x \cdot \omega} \gtrsim \lambda \abs{t}$.
Thus
$$
  \abs{b_{1}(\omega)} \lesssim \bigl( 1 / \abs{t} \bigr) 
  \int_{I} \, dr \lesssim \mu / \abs{t}
$$
and
$$
  \abs{b_{2}(\omega)} \lesssim \bigl( 1 / \lambda \abs{t} \bigr) 
  \int_{I} \left( r + r^{2}/\mu \right) \, dr
  \lesssim (\mu + \lambda) / \abs{t}
  \lesssim \lambda / \abs{t}.
$$
Taking into account \eqref{Angle} we thus get the desired bound, and 
this concludes the proof of Proposition \ref{ImprovedStrProp}.
\end{proof}

\begin{proof}[Proof of Theorem \ref{DyadicThm}(\ref{DiagonalHigh})]
If $\mu \sim \lambda$, this reduces to part \eqref{OffDiagonal}
of the theorem, so we may assume $\mu \ll \lambda$ (and 
$\lambda \gg 1/\varepsilon$). But then 
by an orthogonality argument (see, e.g., the proof of the analogous 
estimate in Theorem 12.1 of \cite{FK}) we reduce to proving 
$$
  \twonorm{u v}{_{t,x}}
  \lesssim  \varepsilon^{1/2} \mu^{1/2} \lambda^{1/2}
  \twonorm{f}{}
  \twonorm{g}{}
$$
in the case where the Fourier transforms of $f, g$ are
supported in (diametrically opposite) cubes with side length
$\sim \mu$ and at distance $\sim \lambda$ from the origin. But this
follows from H\"older's inequality and the estimates 
\eqref{ImprovedWaveStr} and \eqref{ImprovedStr} with $(q,r) = (4,4)$.
\end{proof}

\begin{proof}[Proof of Theorem \ref{DyadicThm}(\ref{DiagonalLow})]
If $\mu \sim \lambda$, this reduces to part \eqref{OffDiagonal}
of the theorem, so we may assume $\mu \ll \lambda \lesssim 1/\varepsilon$.
By orthogonality, we reduce to proving 
\begin{equation}\label{DiagonalLowA}
  \twonorm{u v}{_{t,x}}
  \lesssim \varepsilon^{1/2} \mu
  \twonorm{f}{}
  \twonorm{g}{}
\end{equation}
in the case where $\widehat f, \widehat g$ are
supported in opposite cubes $Q, -Q$ with side length
$\sim \mu$ and at distance $\sim \lambda$ from the origin.
By rescaling $t \to t/\varepsilon$ we further reduce to proving 
\eqref{DiagonalLowA} without the $\varepsilon^{1/2}$ in the right hand 
side, and with $u, v$ given by
\begin{equation}\label{uvRedefined}
  [u(t)]\FT(\xi) = e^{it\abs{\xi}} \widehat f(\xi),
  \qquad
  [v(t)]\FT(\xi) = e^{\pm it \varepsilon^{-1} \alpha(\varepsilon 
  \abs{\xi})} \widehat g(\xi).
\end{equation}
Here $\alpha$ is given by \eqref{AlphaDef}. Then by a standard 
Cauchy-Schwarz argument, see e.g.\ \cite[Sect.\ 3.4]{BMS},
we finally reduce to proving that
\begin{equation}\label{DeltaEstA}
  \int \chi_{\{\eta : \eta \in Q \} \cap
  \{\eta : \eta-\xi \in Q\}}(\eta)
  \delta \bigl( \tau - \abs{\eta} \pm
  k(\abs{\xi-\eta}) \bigr) \, d\eta \lesssim \mu^{2}
\end{equation}
where
\begin{equation}\label{kDef}
  k(\rho) := \varepsilon^{-1} \alpha(\varepsilon \rho).
\end{equation}
(Here and in what follows we use the notation $\chi_{A}$ for the 
characteristic function of a set $A$.)
Then in view of \eqref{AlphaDer} there is an absolute constant $c_{0}$ such that
\begin{equation}\label{kDerBound}
  \abs{k'(\rho)} \le c_{0} < 1
  \quad \text{for all} \quad \rho \lesssim 1/\varepsilon, \,\,
  0 < \varepsilon < 1.
\end{equation}

Denote by $I_{\pm}(\tau,\xi)$ the integral in \eqref{DeltaEstA}. In 
polar coordinates $\eta = r\omega$, $r > 0$, $\omega \in S^{2}$, we 
have $I_{\pm}(\tau,\xi) = \int_{S^{2}} a_{\pm}(\tau,\xi;\omega) \, 
d\sigma(\omega)$, where
$$
  a_{\pm}(\tau,\xi;\omega) := \int_{0}^{\infty} \chi_{\{r : r\omega \in Q \}
  \cap \{r : r\omega-\xi \in Q\}}(r)
  \delta \bigl( \tau - r \pm k(\abs{\xi-r\omega}) \bigr) r^{2} \, dr.
$$
Observe that the $\omega$-support of $a$ is contained in 
a set given by \eqref{Angle}, so it suffices to prove that $a_{\pm} \lesssim 
\lambda^{2}$. Observe also that in the integral defining $a_{\pm}$, the 
variable $r$ is restricted to an interval $I$ of length $\sim \mu$ and centered
at a distance $\lambda$ from the origin.

We shall use the following fact: If $f : \R \to \R$ is differentiable
with $\abs{f'(r)} > 0$, and $f$ has a zero at $r_{0}$, then
\begin{equation}\label{DeltaIdentity}
  \delta\bigl( f(r) \bigr) \, dr = \frac{\delta(r-r_{0}) \, dr}{\abs{f'(r_{0})}}.
\end{equation}
Take
\begin{equation}\label{fDef}
  f(r) := \tau - r \pm k(\abs{\xi-r\omega}),
\end{equation}
for fixed $\tau,\xi,\omega$. Then for $r$ such that $r\omega - \xi \in Q$,
\begin{equation}\label{fDer}
  \abs{f'(r)} = 1 \mp k'(\abs{\xi-r\omega}) 
  \frac{(\xi-r\omega)\cdot\omega}{\abs{\xi-r\omega}}
  \ge 1 - c_{0} \gtrsim 1,
\end{equation}
where we used \eqref{kDerBound} and the assumption $\lambda \lesssim 
1/\varepsilon$. On account of \eqref{DeltaIdentity} 
and \eqref{fDer}, we then get $a_{\pm} \lesssim \lambda^{2}$ as desired.
This concludes the proof of part \eqref{DiagonalLow} of Theorem \ref{DyadicThm}.
\end{proof}

\begin{proof}[Proof of Theorem \ref{DyadicThm}(\ref{OffDiagonal})]
This reduces to proving
\begin{equation}\label{DeltaEstB}
  \int \chi_{\{ \eta: \abs{\eta} \sim \mu \} \cap
  \{ \eta : \abs{\xi-\eta} \sim \lambda \}}(\eta)
  \delta \bigl( \tau - \abs{\eta} \pm 
  k(\abs{\xi-\eta}) \bigr) \, d\eta \lesssim
  \bigl[\min(\mu,\lambda)\bigr]^{2},
\end{equation}
for $k$ defined by \eqref{kDef}. Let us denote the above integral 
by $I_{\pm}(\tau,\xi)$. Passing to polar coordinates we have
$I_{\pm}(\tau,\xi) = \int_{S^{2}} a_{\pm}(\tau,\xi;\omega) \, d\sigma(\omega)$,
where now
$$
  a_{\pm}(\tau,\xi;\omega) := \int_{0}^{\infty} \chi_{\{r : r \sim \mu \}
  \cap \{r : \abs{r\omega-\xi} \sim \lambda \}}(r)
  \delta \bigl( \tau - r \pm k(\abs{\xi-r\omega}) \bigr) r^{2} \, dr.
$$

We split into the cases
{
\renewcommand{\theenumi}{\alph{enumi}}
\begin{enumerate}
  \item\label{FirstCase}
  $\lambda \lesssim 1/\varepsilon$,
  \item\label{SecondCase}
  $\lambda \gg 1/\varepsilon$.
\end{enumerate}

\medskip
\noindent
\emph{Case (\ref{FirstCase}).}
Then in view of \eqref{DeltaIdentity} and \eqref{fDer} with $f(r)$ 
given by \eqref{fDef}, we have $a_{\pm} \lesssim \mu^{2}$. Now
integrate over $S^{2}$, taking into account the fact that on the
support of $a_{\pm}$,
\begin{equation}\label{aSupportAngle}
  \angle(\omega,\xi) = \angle(\eta,\xi) \lesssim \lambda/\mu
  \quad \text{if} \quad \mu \gg \lambda.
\end{equation}

\medskip
\noindent
\emph{Case (\ref{SecondCase}).}
}
By rotational symmetry we may assume $\xi = (\abs{\xi},0,0)$. Now
parametrize the sphere $S^{2}$ by
$$
  (y,\theta) \mapsto \omega = \left(y,\sqrt{1-y^{2}}\, \vec n(\theta) 
  \right),
  \qquad \vec n(\theta) = (\cos \theta, \sin \theta).
$$
Then surface measure $d\sigma(\omega)$ on $S^{2}$ becomes
$dy \, d\theta$. Again we use \eqref{DeltaIdentity} with $f(r)$ given 
by \eqref{fDef}. Observe that $f$ depends implicitly on $y$ but not 
on $\theta$. Denote by $A = A(\tau,\xi)$ the set of $y \in (-1,1)$ such that
$f(r)$ given by \eqref{fDef} has a zero $r_{0} = r_{0}(y) > 0$.
Since $\abs{f'(r)} > 0$, the implicit function theorem guarantees 
that $A$ is open and $r_{0} : A \to (0,\infty)$ is a smooth function.
Differentiating $f\bigl(r_{0}(y)\bigr) = 0$ gives
\begin{equation}\label{fDerExpression}
  0 = f'(r_{0}) r_{0}'(y)
  \mp k'(\abs{\xi-r_{0}\omega}) 
  \frac{r_{0}\abs{\xi}}{\abs{\xi-r_{0}\omega}},
\end{equation}
where we used $\xi \cdot \partial_{y} \omega = \xi_{1} = \abs{\xi}$ and $\omega \cdot 
\partial_{y} \omega = 0$.

Let us suppress the subscript and write 
$r(y)$ instead of $r_{0}(y)$ from now on. Solving 
\eqref{fDerExpression} for $r'(y)$ and using the fact that $f'(r) < 
0$, we see that $\partial r / \partial y$ is either strictly negative 
or strictly positive, depending on whether we have the $+$ sign or 
the $-$ sign in \eqref{DeltaEstB}. The function $r(y)$ is therefore a 
change of variables.

With this information in hand, we solve \eqref{fDerExpression} for
$f'(r)$ and substitute into \eqref{DeltaIdentity}, thus arriving at the
identity
$$
  \int F(\eta) \delta \bigl( \tau - \abs{\eta} \pm
  k(\abs{\xi-\eta}) \bigr) \, d\eta
  = \int \!\! \int F(r\omega)
  \frac{r \abs{\xi-r\omega}}{\abs{\xi}k'(\abs{\xi-r\omega})}
  \abs{ \frac{\partial r}{\partial y} } \, dy \, d\theta.
$$
Changing variables $y \to r$ finally gives
\begin{equation}\label{IntegralIdentity}
  \int F(\eta) \delta \bigl( \tau - \abs{\eta} \pm
  k(\abs{\xi-\eta}) \bigr) \, d\eta
  = \int \!\! \int F(r\omega)
  \frac{r \abs{\xi-r\omega}}{\abs{\xi}k'(\abs{\xi-r\omega})}
  \, dr \, d\theta,
\end{equation}
where $\omega$ is now a function of $r$ and $\theta$. 
We apply this with
$$
  F(\eta) := \chi_{\{ \eta: \abs{\eta} \sim \mu \} \cap
  \{ \eta : \abs{\xi-\eta} \sim \lambda \}}(\eta).
$$
Since $\lambda \gg 1/\varepsilon$, we see from \eqref{AlphaDer} that
$k'(\abs{\xi-r\omega}) \sim 1$, whence
\begin{equation}\label{IntegrandEst}
  F(r\omega)
  \frac{r \abs{\xi-r\omega}}{\abs{\xi}k'(\abs{\xi-r\omega})}
  \sim
  \frac{\mu\lambda}{\abs{\xi}} F(r\omega).
\end{equation}

We now split into the subcases
\begin{itemize}
  \item[(b1)] $\mu \ll \lambda$,
  \item[(b2)] $\mu \sim \lambda$,
  \item[(b3)] $\mu \gg \lambda$.
\end{itemize}

\medskip
\noindent
\emph{Case (b2).}
In this case we can prove the estimate in Theorem 
\ref{DyadicThm}(\ref{OffDiagonal}) directly, by applying H\"older's 
inequality followed by the linear Strichartz estimate \eqref{LinearStr} 
with $(q,r) = (4,4)$. (This works because we are at high frequency, 
i.e. $\gg 1/\varepsilon$.)

\medskip
\noindent
\emph{Case (b1).}
Then $\abs{\xi} \sim \lambda$, so the desired estimate \eqref{DeltaEstB}
follows readily from \eqref{IntegrandEst} and \eqref{IntegralIdentity}.

\medskip
\noindent
\emph{Case (b3).}
Then $\abs{\xi} \sim \mu$, so \eqref{IntegrandEst} and \eqref{IntegralIdentity} 
imply
$$
  I_{\pm}(\tau,\xi) \lesssim \lambda \int \!\! \int \chi_{ \{r : r \sim \mu\} 
  \cap \{ r : \abs{\xi-r\omega} \sim \lambda \} }(r) \, dr \, d\theta.
$$
Recall that $\omega$ is now a function of $(r,\theta)$. However, 
$\abs{\xi-r\omega}$ is independent of $\theta$, so by a slight abuse 
of notation we will simply write $\omega = \omega(r)$ and integrate 
out $\theta$, leaving us with
$$
  \lambda \int \chi_{ \{r : r \sim \mu\} 
  \cap \{ r : \abs{\xi-r\omega(r)} \sim \lambda \} }(r) \, dr.
$$
Clearly it suffices to prove that the support of the integrand is contained
in an interval of length $\sim \lambda$. Let us assume there is no such
interval, and obtain a contradiction. Fix a point $r_{0}$ in the 
support, and write
$$
  r = r_{0} + \kappa
$$
for a general point $r$ in the support. In view of our assumption, 
$\kappa$ varies on a scale $\gg \lambda$. Thus, if we can show that
\begin{equation}\label{Variation}
  \abs{\xi-r\omega(r)}^{2} = a + \kappa^{2} + O(\lambda \kappa + 
  \lambda^{2}),
\end{equation}
for some constant $a$, it follows that $\abs{\xi-r\omega(r)}$ also varies on a 
scale $\gg \lambda$, and we have the contradiction we 
seek, since $\abs{\xi-r\omega(r)} \sim \lambda$ on the support.

To prove \eqref{Variation}, write
$$
  \abs{\xi-r\omega}^{2} = \abs{\xi}^{2} + (r^{2} - 2 r \abs{\xi}) + 2r 
  (1-\omega_{1}) \abs{\xi}.
$$
On account of \eqref{aSupportAngle} we have $1-\omega_{1} \lesssim 
(\lambda/\mu)^{2}$, so the last term on the right hand side is
$O(\lambda^{2})$. For the second term we calculate
$$
  r^{2} - 2 r \abs{\xi} = (r_{0}^{2} - 2r_{0} \abs{\xi})
  + 2(r_{0} - \abs{\xi}) \kappa + \kappa^{2}.
$$
But
$$
  \bigabs{r_{0} - \abs{\xi}} \le \abs{r_{0}\omega(r_{0}) - \xi}
  \sim \lambda,
$$
so we conclude that \eqref{Variation} holds.
This ends the proof of Theorem \ref{DyadicThm}.
\end{proof}

\section{Proof of Theorem \ref{NullFormThm}}\label{NullProof}

As remarked,  by a standard procedure this reduces to some well-known bilinear estimates 
for the homogeneous wave equation. The first observation is that by 
rescaling $x \to \varepsilon x$ we can reduce to the case 
$\varepsilon = 1$. Thus we suppress the subscript on $H^{s,\theta}$ etc.\ from now on.

Some notation: For $s \in \R$, let
$D^{s}$, $D_{+}^{s}$ and $D_{-}^{s}$ be the Fourier multipliers
$$
  \bigl(D^{s} u\bigr)\FT = \abs{\xi}^{s} \widehat u,
  \quad
  \bigl(D_{+}^{s} u\bigr)\FT = \bigl(\abs{\tau} + \abs{\xi} \bigr)^{s} \widehat u,
  \quad
  \bigl(D_{-}^{s} u\bigr)\FT = \bigabs{\abs{\tau} - \abs{\xi}}^{s} \widehat u.
$$
The notation $u \precsim v$ means $\bigabs{\widehat u} 
\lesssim \widehat v$. We are concerned with bilinear operators 
$B(u,v)$ of the form
$$
  [B(u,v)]\FT(\tau,\xi)
  = \int b(\tau-\lambda,\xi-\eta;\lambda,\eta)
  \widehat u(\tau-\lambda,\xi-\eta)
  \widehat v(\lambda,\eta) \, d\lambda \, d\eta,
$$
where $b(\tau,\xi;\lambda,\eta)$ is the \emph{symbol} of $B$.
The symbols of the null forms $Q_{0}$, $Q_{ij}$ and $Q_{0j}$ are, 
respectively,
\begin{subequations}\label{NullFormSymbols}
\begin{align}
  \label{Q0Symbol}
  q_{0}(\tau,\xi;\lambda,\eta) &= \tau \lambda - \xi \cdot \eta,
  \\
  \label{QijSymbol}
  q_{ij}(\tau,\xi;\lambda,\eta) &= - \xi_{i} \eta_{j} + \xi_{j} 
  \eta_{i},
  \\
  \label{Q0jSymbol}
  q_{0j}(\tau,\xi;\lambda,\eta) &= - \tau \eta_{j} + \lambda \xi_{j}.
\end{align}
\end{subequations}
Since we rely on estimates for the \emph{absolute values} of these symbols, and since all norms 
involved only depend on the absolute value of the Fourier transform, 
we may assume $\widehat u, \widehat v \ge 0$ henceforth.

For $s \in \R$, let $R^{s}$ be the bilinear operator with symbol $r^{s}$, where
$$
  r(\tau,\xi;\lambda,\eta) =
  \begin{cases}
  \abs{\xi} + \abs{\eta} - \abs{\xi+\eta} &\text{if} \quad \tau\lambda \ge
  0, \\
  \abs{\xi+\eta} - \bigl\vert \abs{\xi} - \abs{\eta} \bigr\vert &\text{if}
  \quad \tau\lambda < 0.
  \end{cases}
$$
We shall need the estimate, for $\theta > 1/2$,
\begin{equation}\label{RsqrtEst}
  \bigtwonorm{R^{1/2}(u,v)}{} \lesssim \norm{u}_{H^{0,\theta}}
  \norm{v}_{H^{3/2,\theta}},
\end{equation}
which derives from an estimate for the homogeneous wave equation via the 
Transfer Principle; see \cite{KlSe} for the details. We also need
\begin{equation}\label{RExpansion}
  R^{s}(u,v)
  \precsim
  D_{-}^{s} (u v) + (D_{-}^{s} u) v + u D_{-}^{s}v
  \qquad (s \ge 0).
\end{equation}
This follows easily from the triangle inequality, if one keeps track 
of the signs of $\tau$ and $\lambda$ as in the proof of the following 
lemma, which is more or less standard.

\begin{lemma}\label{NullFormExpansionLemma}
The following estimates hold:
\begin{subequations}
\begin{align}
  \label{QijExpansion1}
  Q_{ij}(u,v) &\precsim 
  R^{1/2}(D u,D^{1/2}v) + R^{1/2}(D^{1/2}u,v)
  \\
  \label{QijExpansion2}
  &\precsim 
  R^{(1/2)^{-}}(Du,D^{(1/2)^{+}}v) + R^{(1/2)^{-}}(D^{(1/2)^{+}}u,D v),
  \\
  \label{Q0jExpansion}
  Q_{0j}(u,v) &\precsim [\text{r.h.s.\eqref{QijExpansion1}}]
  + D u \cdot D_{-} v + D_{-} u \cdot D v,
  \\
  \label{Q0Expansion}
  Q_{0}(u,v) &\precsim [\text{r.h.s.\eqref{QijExpansion1}}]
  + D_{+} u \cdot D_{-} v
  + D_{-} u \cdot D_{+} v.
\end{align}
\end{subequations}
\end{lemma}

\begin{proof} All these statements reduce to estimates on the absolute 
values of the symbols \eqref{NullFormSymbols}. First,
by \cite[Lemma 13.2]{FK} we have
$$
  \abs{q_{ij}(\tau,\xi;\lambda,\eta)}
  \le \abs{\xi \times \eta} \le
  \abs{\xi}^{1/2} \abs{\eta}^{1/2} \abs{\xi+\eta}^{1/2}
  [r(\tau,\xi;\lambda,\eta)]^{1/2},
$$
where $r$ is the symbol of $R$ as defined above.
Then \eqref{QijExpansion1} and \eqref{QijExpansion2} follow,
in view of the fact that
\begin{equation}\label{rEst}
  r(\tau,\xi;\lambda,\eta) \le 2 \min(\abs{\xi},\abs{\eta}).
\end{equation}
To prove \eqref{Q0jExpansion}, write
$$
  q_{0j}(\tau,\xi;\lambda,\eta) = ( \epsilon_{1} \abs{\xi} - \tau) \eta_{j}
  + ( \lambda - \epsilon_{2} \abs{\eta} ) \xi_{j}
  - \epsilon_{1} ( \abs{\xi} \eta_{j} - \epsilon_{1} 
  \epsilon_{2} \abs{\eta} \xi_{j} ),
$$
where $\epsilon_{1}$ and $\epsilon_{2}$ are the signs of $\tau$ and 
$\lambda$, respectively. That is, $\epsilon_{1} \tau = \abs{\tau}$
and $\epsilon_{2} \lambda = \abs{\lambda}$.
Now take absolute values and use the fact (see \cite[Lemma 
13.2]{FK}) that
$$
  \bigabs{ \abs{\xi} \eta_{j} \pm \abs{\eta} \xi_{j} }
  \le \abs{\xi}^{1/2} \abs{\eta}^{1/2} (\abs{\xi} + \abs{\eta})^{1/2}
  [r(\tau,\xi;\lambda,\eta)]^{1/2}
$$
holds for all $\tau,\xi,\lambda,\eta$. (The sign in the left hand 
side is independent of the signs of $\tau, \lambda$.)
This proves \eqref{Q0jExpansion}. The proof of \eqref{Q0Expansion} is 
similar. Write
$$
  q_{0}(\tau,\xi;\lambda,\eta) = ( \tau - \epsilon_{1} \abs{\xi} ) 
  \lambda
  + ( \lambda - \epsilon_{2} \abs{\eta} ) \epsilon_{1} \abs{\xi} 
  + \epsilon_{1} \epsilon_{2} \abs{\eta} \abs{\xi} - \xi \cdot \eta.
$$
Then use (see \cite[Lemma 13.2]{FK})
$$
  \bigabs{\abs{\eta} \abs{\xi} - \xi \cdot \eta} \le 
  (\abs{\xi} + \abs{\eta}) r(\tau,\xi;\lambda,\eta)
$$
and \eqref{rEst}.
\end{proof}

Finally, we need the estimate (here $s_{1}, s_{2}, \theta_{1}, 
\theta_{2} \ge 0$)
\begin{equation}\label{SimpleBilinear}
  \twonorm{uv}{} \lesssim \norm{u}_{H^{s_{1},\theta_{1}}} 
  \norm{u}_{H^{s_{2},\theta_{2}}}
  \quad \text{for} \quad s_{1} + s_{2} > \tfrac{3}{2},
  \quad \theta_{1} + \theta_{2} > \tfrac{1}{2}.
\end{equation}
See \cite[Proposition A.1]{KlSe} for the simple proof of this fact.

\medskip\medskip\medskip
\noindent
We are now ready to prove Theorem \ref{NullFormThm}. By interpolation, we reduce to
\begin{align}
  \label{BilinearA}
  \twonorm{Q(u,v)}{} &\lesssim \norm{u}_{\dot \scrH^{1,\theta}}
  \norm{v}_{\scrH^{2,\theta}} 
  \\
  \label{BilinearB}
  \norm{Q(u,v)}_{H^{0,(-1/2)^{-}}} &\lesssim \norm{u}_{\dot \scrH^{1,\theta}} 
  \norm{v}_{\scrH^{(3/2)^{+},1}} 
\end{align}
where $\norm{u}_{\dot \scrH^{1,\theta}}$ in the right hand side can
be replaced by $\norm{u}_{\dot H^{1,\theta}}$ if $Q = Q_{ij}$.

\begin{proof}[Proof of \eqref{BilinearA}.] First observe that for
the last two terms in the right hand sides of \eqref{Q0jExpansion} and
\eqref{Q0Expansion}, the estimate reduces to special cases of
\eqref{SimpleBilinear}, since we can always replace $D_{-}$ by 
$D_{-}^{\theta}D_{+}^{1-\theta}$. Thus, it only remains to prove the estimate for
the right hand side of \eqref{QijExpansion1}, but this reduces to
\eqref{RsqrtEst}.
\end{proof}

\begin{proof}[Proof of \eqref{BilinearB}.]
First consider $Q_{ij}$. Applying \eqref{RExpansion} to \eqref{QijExpansion2}, we reduce to
\begin{subequations}\label{BilinearBReduction}
\begin{align}
  \label{BilinearB1}
  \twonorm{uv}{} &\lesssim \norm{u}_{H^{0,(1/2)^{+}}} \norm{v}_{H^{1^{+},1}},
  \\
  \label{BilinearB2}
  \twonorm{uv}{} &\lesssim \norm{u}_{\dot H^{(1/2)^{-},(1/2)^{+}}} 
  \norm{v}_{H^{(1/2)^{+},1}},
  \\
  \label{BilinearB3}
  \norm{uv}_{H^{0,(-1/2)^{-}}}
  &\lesssim \twonorm{u}{} \norm{v}_{H^{1^{+},1}},
  \\
  \label{BilinearB4}
  \norm{uv}_{H^{0,(-1/2)^{-}}}
  &\lesssim \norm{u}_{\dot H^{(1/2)^{-},0}} \norm{v}_{H^{(1/2)^{+},1}},
  \\
  \label{BilinearB5}
  \norm{uv}_{H^{0,(-1/2)^{-}}}
  &\lesssim 
  \norm{u}_{H^{0,(1/2)^{+}}} \norm{v}_{H^{1^{+},(1/2)^{+}}},
  \\
  \label{BilinearB6}
  \norm{uv}_{H^{0,(-1/2)^{-}}}
  &\lesssim \norm{u}_{\dot H^{(1/2)^{-},(1/2)^{+}}} 
  \norm{v}_{H^{(1/2)^{+},(1/2)^{+}}}.
\end{align}
\end{subequations}
Via duality and the Transfer Principle, these reduce to the estimates 
in Corollaries \ref{DyadicCor1} and \ref{DyadicCor2}, which are valid 
in the case where $u,v$ are both solutions of the homogeneous wave 
equation, as remarked in Sect.\ \ref{BilinearSection}.

It remains to consider the second and third terms in the right hand 
sides of \eqref{Q0jExpansion} and \eqref{Q0Expansion}. For the second 
term we can apply \eqref{SimpleBilinear} directly, while for the third term
we replace $D_{-}$ by $D_{-}^{(1/2)^{-}}D_{+}^{(1/2)^{+}}$, thus reducing
to \eqref{BilinearB4}.
\end{proof}

\medskip\medskip\medskip
\noindent
{\bf Acknowledgment.} {\em

Financial support by the Austrian START project ``Nonlinear Schr\"odinger and quantum Boltzmann
equations" (FWF Y137-TEC) of N.J.M.\ and by the European network HYKE
(HPRN-CT-2002-00282) as well as by the OeAD (``acciones integradas") is acknowledged.
}

\end{document}